\documentclass[11pt,reqno]{amsart}
\usepackage{amsmath,amssymb}
\usepackage{graphicx,psfrag}
\usepackage{amsthm}
\usepackage{enumerate,placeins,wrapfig,booktabs,chngcntr}
\usepackage{subcaption}

\usepackage[colorlinks=true, linkcolor=blue, citecolor=green]{hyperref}
\newtheorem{thm}{Theorem}[section] 
\newtheorem{conj}[thm]{Conjecture}

\newtheorem{rmk}[thm]{Remark}
\newtheorem{cor}[thm]{Corollary}
\newtheorem{lem}[thm]{Lemma}
\newtheorem{prop}[thm]{Proposition}

\newtheorem{defn}[thm]{Definition}

  \allowdisplaybreaks[4]

\begin{document}

	\title[Homotopic to the Identity]{Partially Hyperbolic Diffeomorphisms Homotopic to the Identity in Dimension Three}

	\author{Ziqiang Feng}
    \address{Beijing International Center for Mathematical Research, Peking University, Beijing, 100871, China.}
    \email{zqfeng@pku.edu.cn}

    \author{Ra\'{u}l Ures}
    \address{Department of Mathematics and Shenzhen International Center for Mathematics, Southern University of Science and Technology, Shenzhen, 518000, China.}
    \email{ures@sustech.edu.cn}

    \thanks{This work was partially supported by National Key R\&D Program of China 2022YFA1005801, NSFC 12071202, and NSFC 12161141002.}

	\begin{abstract}
		We show that any conservative partially hyperbolic diffeomorphism homotopic to the identity is accessible unless the fundamental group of its ambient 3-manifold is virtually solvable. As a consequence, such diffeomorphisms are ergodic, giving an affirmative answer to the Hertz-Hertz-Ures Ergodicity Conjecture in the homotopy class of identity.
	\end{abstract}

    \footnote{Preliminary version}
    \date{\today}
    
	\maketitle
	

	\vspace{.5cm}

	

\tableofcontents

\section{Introduction}


In the late 19th century, Boltzmann introduced \textbf{ergodicity} as a hypothesis in thermodynamics to describe the long-term statistical behavior of gas particles. Since then, ergodic theory has been extensively studied as a powerful framework for formulating the equidistribution of typical orbits in dynamical systems, with applications across various mathematical and physical disciplines. Through the ergodic theorems \cite{Birkhoff1931,Birkhoff1932,Neumann1932}, a dynamical system $f: M\rightarrow M$ is \emph{ergodic} with respect to an invariant measure $\mu$ if for any continuous function $\phi: M\rightarrow \mathbb{R}$, the time average along $\mu$-typical orbits equals the space average:
$$
\lim\limits_{n\rightarrow\infty}\frac{1}{n}\sum\limits_{i=0}^{n-1}\phi(f^i(x)) = \int_M \phi \, d\mu \quad \text{for } \mu\text{-a.e. } x.
$$
Following the work of Oxtoby and Ulam \cite{OU1941} on the genericity of ergodic homeomorphisms, a parallel conjecture emerged for smooth systems, which persisted until the advent of KAM theory \cite{KAM54,Arnold65,Moser62}. The KAM phenomenon demonstrated the existence of smooth integrable systems of elliptic type that are stably non-ergodic, shifting focus to systems exhibiting hyperbolic features as potential candidates for generic ergodicity.

In his seminal work, Hopf \cite{Hopf1939} developed the \emph{Hopf argument}, establishing ergodicity for geodesic flows on compact negatively curved surfaces. Building on this foundation, Anosov and Sinai \cite{Anosov1967,AnosovSinai67} extended the technique to \emph{uniformly hyperbolic systems} (Anosov systems), characterized by two invariant subbundles that exhibit uniform expansion and contraction under the differential.

Beyond uniform hyperbolicity, the Hopf argument might naturally extend to two broader contexts: nonuniformly hyperbolic systems (characterized by non-vanishing Lyapunov exponents and now known as Pesin theory) and partially hyperbolic systems (featuring a dynamical preserving direction where hyperbolic directions dominate). Our work focuses on the partially hyperbolic setting.

The study of ergodic properties in partially hyperbolic systems originated independently with Brin and Pesin \cite{BrinPesin74} for skew products and frame flows, and with Pugh and Shub \cite{PughShub72} for Anosov actions. Contemporary definitions differ slightly from these original formulations, with modern treatments considering more general cases. A diffeomorphism $f: M\rightarrow M$ on a compact Riemannian manifold is called (pointwise) \emph{partially hyperbolic} if there exists a nontrivial $Df$-invariant splitting $TM = E^s \oplus E^c \oplus E^u$ satisfying:
$$
\|Df_x(v^s)\| < 1 < \|Df_x(v^u)\| \quad \text{and} \quad \|Df_x(v^s)\| < \|Df_x(v^c)\| < \|Df_x(v^u)\|
$$
for all $x \in M$ and unit vectors $v^\sigma \in E^\sigma_x$ ($\sigma = s, c, u$). To apply the Hopf argument effectively, we focus on invariant measures that are absolutely continuous with respect to a volume form. A system $f: M\rightarrow M$ preserving such a measure is termed \emph{conservative}.

Following the groundbreaking work of Grayson, Pugh, and Shub \cite{GPS94}, Pugh and Shub formulated their Stable Ergodicity Conjecture: (stable) ergodicity constitutes a generic property among conservative partially hyperbolic diffeomorphisms. Their approach utilized accessibility - the property that any two points can be joined through paths consisting of consecutive arcs tangent to either $E^s$ or $E^u$. They proposed two refined conjectures: (1) Accessibility is both open and dense among partially hyperbolic diffeomorphisms; (2) (Essential) accessibility implies ergodicity for $C^r$ ($r > 1$) conservative partially hyperbolic diffeomorphisms. Substantial progress has been made toward proving the Stable Ergodicity Conjecture, as evidenced by works including \cite{Wilkinson98, BPW00, NT01, BDP02, DW03, Hertz05, 08invent, BW10annals, HHTU11, ACW21, AV20}.

Research \cite{08invent, BW10annals} has demonstrated that ergodicity constitutes a prevalent property among conservative partially hyperbolic diffeomorphisms in 3-manifolds. This naturally leads to the challenge of characterizing this prevalence through a complete identification of obstructions to ergodicity in dimension three - the minimal dimension permitting partial hyperbolicity.

In 2008, Hertz, Hertz, and Ures \cite{2008nil} (see also \cite{2018survey}) formulated the following conjecture:

\begin{conj}[HHU Ergodicity Conjecture]
    Let $f: M^3\rightarrow M^3$ be a $C^r$ ($r > 1$) conservative partially hyperbolic diffeomorphism of a closed 3-manifold that fails to be ergodic. Then there exists an embedded 2-torus tangent to $E^s \oplus E^u$. Moreover, under orientability assumptions, the manifold $M^3$ must be one of:
    \begin{enumerate}
        \item the 3-torus $\mathbb{T}^3$;
        \item the mapping torus of $-id: \mathbb{T}^2 \rightarrow \mathbb{T}^2$;
        \item the mapping torus of an Anosov diffeomorphism on $\mathbb{T}^2$.
    \end{enumerate}
\end{conj}

This paper establishes a complete proof of this conjecture within the homotopy class of the identity:

\begin{thm}\label{ergodic_id}
	Let $f:M^3\rightarrow M^3$ be a $C^r$, $r>1$, conservative partially hyperbolic diffeomorphism homotopic to the identity on a closed 3-manifold. Then, $f$ is a K-system (and thus ergodic) if there is no embedded 2-torus tangent to $E^s\oplus E^u$.
\end{thm}

We note that the forward implication does not hold, as ergodic examples with $su$-tori can be readily constructed. For instance, consider an irrational-time map of a suspension Anosov flow, which remains ergodic on torus bundles over the circle where fibers align with the $E^s \oplus E^u$ distribution (see \cite{FU1}). Furthermore, conservative partially hyperbolic diffeomorphisms with a single $su$-torus and a single open accessibility class on solvmanifolds are essentially accessible and therefore ergodic by \cite{08invent, BW10annals}.

Significant progress has been made in verifying the HHU Ergodicity Conjecture across typical 3-manifold classes. Hertz, Hertz, and Ures \cite{2008nil} provided foundational motivation by demonstrating that all conservative partially hyperbolic diffeomorphisms on nilmanifolds are ergodic, thereby confirming a weak version of the conjecture for manifolds with virtually solvable fundamental group. Gan and Shi \cite{GS20DA}, building on contributions from \cite{HamU14CCM}, fully resolved the conjecture for these manifolds. This key step transfers attention to manifolds with non-solvable fundamental group. Subsequent efforts have addressed more complex cases: Hammerlindl, Hertz, and Ures \cite{2020Seifert} and Fenley-Potrie \cite{FP_hyperbolic} established affirmative results for Seifert manifolds in specific isotopy classes. Notably, Fenley and Potrie have announced a complete proof for Seifert manifolds without isotopy restrictions. They also established ergodicity for all conservative partially hyperbolic diffeomorphisms on hyperbolic 3-manifolds \cite{FP_hyperbolic}.

Though conservative non-ergodic partially hyperbolic diffeomorphisms constitute a meager set, the core challenge lies in precisely characterizing this exceptional class. Adopting a dynamical perspective rather than focusing on manifold topology, Fenley and Potrie demonstrated accessibility (and consequently ergodicity) for discretized Anosov flows \cite{FP_hyperbolic} and collapsed Anosov flows \cite{FP21accessible}. In \cite{FU1}, we proved ergodicity for all partially hyperbolic diffeomorphisms without periodic points, thereby resolving the HHU Ergodicity conjecture in that setting. These techniques were extended by the first author \cite{Feng24} to systems with bounded center dynamics, broadening the class of diffeomorphisms satisfying the conjecture.

Inspired by the Pugh-Shub Stable Ergodicity Conjecture, as shown in \cite{GPS94,08invent, BW10annals}, examining ergodicity can be reduced to a geometric problem studying accessibility with less regularity condition and a weaker constraint than conservative property. 

Our following result establishes accessibility in the identity homotopy class:

\begin{thm}\label{accessible-id}
    Let $f:M^3\rightarrow M^3$ be a $C^1$ partially hyperbolic diffeomorphism homotopic to the identity on a closed 3-manifold. If $\pi_1(M)$ is not virtually solvable and $NW(f) = M$, then $f$ is accessible.
\end{thm}


This result offers a step forward by removing geometric constraints on $M$. By contrast, this theorem was previously established in case where the ambient manifold is a Seifert 3-manifold \cite{2020Seifert}, with independent proofs emerging from \cite{FP_hyperbolic, BFFP1} through distinct approaches. Furthermore, the framework naturally extends to hyperbolic 3-manifolds due to Mostow rigidity, as demonstrated in \cite{FP_hyperbolic, FP-gafa, FP21accessible}. 

Accessibility has emerged as a significant property of independent interest in partially hyperbolic dynamics, with applications extending far beyond ergodicity verification. Notably, it enables generalizations of periodic obstructions for solving cohomological equations in the partially hyperbolic setting \cite{Wilkinson13}. Furthermore, accessibility serves as a powerful tool for studying properties of special ergodic measures, as detailed in \cite{HHTU12,VY13,AVW-flow,CP}. While substantial work has focused on establishing accessibility itself, Hammerlindl and Shi \cite{Ham17CMH,HS21DA} classified accessibility classes for manifolds with virtually solvable fundamental groups. Recent work \cite{FU2} demonstrates accessibility for partially hyperbolic diffeomorphisms without periodic points by removing the non-wandering constraint, building on techniques from \cite{FU1}. We refer the reader to \cite{Wilkinson10-survey,2018survey} for surveys of related developments.

Combining Brin's foundational work \cite{Brin75_transitive} with our results yields the following corollary, which may hold independent interest:

\begin{cor}
    Let $f:M^3\rightarrow M^3$ be a $C^1$ partially hyperbolic diffeomorphism homotopic to the identity on a closed 3-manifold. If $\pi_1(M)$ is not (virtually) solvable and $NW(f)=M$, then $f$ is necessarily transitive.
\end{cor}

Ergodicity and transitivity serve as fundamental concepts in partially hyperbolic dynamics, representing key phenomena from measure-theoretical and topological perspectives respectively. Moreover, partial hyperbolicity naturally arises in both stable ergodicity theory \cite{GPS94} and the study of robust transitivity \cite{DPU99}. The following result constructs a conceptual bridge between transitivity and ergodicity within the homotopy class of the identity:

\begin{thm}\label{ergodic=transitive}
    Let $f:M\rightarrow M$ be a $C^r, r>1,$ conservative partially hyperbolic diffeomorphism homotopic to the identity on a closed 3-manifold. Then, $f$ is transitive if and only if it is ergodic.
\end{thm}

As a consequence of Theorem~\ref{accessible-id}, the following result clarifies a mechanism to examine the accessibility and ergodicity of partially hyperbolic diffeomorphisms that are homotopic to the identity. 

\begin{cor}\label{cor-id-nonDAF}
	Let $f:M\rightarrow M$ be a $C^1$ partially hyperbolic diffeomorphism homotopic to the identity in a closed 3-manifold satisfying either $NW(f)=M$ or dynamical coherence. Assume that any finite iterate of $f$ is not a discretized suspension Anosov flow. Then, $f$ is accessible. In particular, $f$ is a K-system, and thus ergodic, provided that it is $C^r$, $r>1$, and conservative.
\end{cor}

We now outline our proof strategy. Our principal novelty reside in Theorem~\ref{accessible-id}, making its proof our central focus. We address extensions to other results in the concluding section.

To establish Theorem~\ref{accessible-id}, we employ Theorem~\ref{su-foliation} to reduce the problem to analyzing an invariant lamination whose leaves are saturated by stable and unstable manifolds. The non-virtually solvable fundamental group condition eliminates compact lamination leaves (Theorem~\ref{maptori}), while the complementary geometric structure (which may be empty) ensures lamination minimality and trivial foliation extension. In light of \cite{BI08}, we can naturally utilize branching foliations, which may be not true foliations due to the lack of dynamical coherence, to interact with the lamination. Critically, stable/unstable subfoliations within lamination leaves remain uncharacterized. Our core objective involves deriving intrinsic structures from intersections between the lamination and branching subfoliations.

We summarize our technical results for general codimension-one foliations in closed 3-manifolds as follows:

\begin{thm}\label{thm-foliations}
    Let $\mathcal{F}_1$, $\mathcal{F}_2$ be two transverse minimal $\mathbb{R}$-covered foliations by non-compact Gromov hyperbolic leaves in a closed 3-manifold $M \neq \mathbb{T}^3$, and let $\mathcal{G}$ be the intersecting one-dimensional foliation. Denote by $\widetilde{\mathcal{F}}_1$, $\widetilde{\mathcal{F}}_2$, and $\widetilde{\mathcal{G}}$ the associated lifts to the universal cover $\widetilde{M}$. Then the leaves of $\mathcal{F}_1$ constitute the weak-stable foliation of a topologically mixing topological Anosov flow $\phi_t$ if either:
    \begin{enumerate}
        \item $\widetilde{\mathcal{G}}$ has non-dense limit set in each leaf of $\widetilde{\mathcal{F}}_1$; or
        \item $\mathcal{F}_1, \mathcal{F}_2$ are uniform foliations with $\widetilde{\mathcal{G}}$ having dense limit sets and Hausdorff leaf spaces in $\widetilde{\mathcal{F}}_1$-leaves, in which case $\mathcal{G}$-leaves coincide with $\phi_t$-orbits.
    \end{enumerate}
    Consequently, $M$ admits a transitive smooth Anosov flow.
\end{thm}

This result proves more than sufficient for establishing Theorem~\ref{accessible-id}. To apply it, we collapse complementary regions of our given lamination and replace the chosen branching foliation with a well-approximated foliation. The homotopy class constraint ensures satisfaction of Theorem~\ref{thm-foliations}'s hypotheses. This process yields rapid contradictions through existence proofs for closed stable/unstable intersection leaves.

The paper's structure unfolds as follows: After providing some preliminaries in Section~\ref{preliminaries}, we will concentrate on the analysis of Gromov hyperbolic foliation pairs. Section~\ref{section-asymptotic} examines large-scale behavior of intersected subfoliations near Gromov boundaries. Section~\ref{section-distance} relates geodesics to intersected leaves, providing foundational insights for the study in Section~\ref{section-Hausdorff} of Hausdorff leaf space geometries. The case of non-dense limit set in the ideal boundary is treated in Section~\ref{section-degenerate}, while Section~\ref{section-nonHausdorff} addresses non-Hausdorff configurations. Finally, Section~\ref{section-accessible} includes applications of these results to partially hyperbolic systems, containing all necessary proofs.


\section{Preliminaries}\label{preliminaries}

\subsection{Partial hyperbolicity}

Let \( M \) be a compact Riemannian manifold. A diffeomorphism \( f: M \to M \) is \emph{partially hyperbolic} if the tangent bundle splits into three nontrivial \( Df \)-invariant subbundles \( TM = E^s \oplus E^c \oplus E^u \), such that for an adapted metric and for all \( x \in M \) and unit vectors \( v^\sigma \in E^\sigma_x \) (\( \sigma = s, c, u \)):
\begin{equation*}
    \begin{aligned}
        \|Df(x)v^s\| < 1 < \|Df(x)v^u\| \quad
        \text{and} \quad \|Df(x)v^s\| < \|Df(x)v^c\| < \|Df(x)v^u\|.
    \end{aligned}
\end{equation*}

\begin{thm}{\cite{BrinPesin74, HPS77}}
    Let \( f: M \to M \) be a partially hyperbolic diffeomorphism. Then there exist unique invariant foliations \( \mathcal{F}^s \) and \( \mathcal{F}^u \), tangent to \( E^s \) and \( E^u \), respectively.
\end{thm}

This classical result establishes the unique integrability of the strong stable and unstable bundles. The foliations \( \mathcal{F}^s \) and \( \mathcal{F}^u \) are called the \emph{stable} and \emph{unstable foliations}. However, the center bundle \( E^c \) is not always integrable, and its integrability remains a longstanding open problem. A partially hyperbolic diffeomorphism is \emph{dynamically coherent} if there exist invariant foliations \( \mathcal{F}^{cs} \) and \( \mathcal{F}^{cu} \), tangent to \( E^s \oplus E^c \) and \( E^c \oplus E^u \), respectively. Dynamical coherence is not always guaranteed; see \cite{Wilkinson98, BurnsWilkinson08, 2016example, BGHP3} for incoherent examples. The foliations \( \mathcal{F}^{cs} \) and \( \mathcal{F}^{cu} \) are called the \emph{center-stable} and \emph{center-unstable foliations}. Under dynamical coherence, intersecting \( \mathcal{F}^{cs} \) and \( \mathcal{F}^{cu} \) yields a one-dimensional \( f \)-invariant foliation tangent to \( E^c \). If \( E^c \) is integrable, the resulting foliation \( \mathcal{F}^c \) is called the \emph{center foliation}.

Three-dimensional manifolds admitting compact tori tangent to \( E^s \oplus E^c \), \( E^c \oplus E^u \), or \( E^s \oplus E^u \) are classified below. Such tori force the fundamental group to be (virtually) solvable.

\begin{thm} \cite{2011TORI}\label{maptori}
    Let \( f: M^3 \to M^3 \) be a partially hyperbolic diffeomorphism on a closed orientable 3-manifold. If there exists an \( f \)-invariant 2-torus \( T \) tangent to \( E^s \oplus E^u \), \( E^c \oplus E^u \), or \( E^c \oplus E^s \), then \( M^3 \) must be:
    \begin{enumerate}
        \item The 3-torus \( \mathbb{T}^3 \);
        \item The mapping torus of \( -\mathrm{id}: \mathbb{T}^2 \to \mathbb{T}^2 \); or
        \item The mapping torus of a hyperbolic automorphism of \( \mathbb{T}^2 \).
    \end{enumerate}
\end{thm}

A set is \emph{s-saturated} (resp.~\emph{u-saturated}) if it is a union of stable (resp.~unstable) leaves. A set is \emph{su-saturated} if it is both $s$- and $u$-saturated. The \emph{accessibility class} \( AC(x) \) of \( x \in M \) is the minimal $su$-saturated set containing \( x \). Points in the same accessibility class are connected by an \emph{$su$-path} (piecewise tangent to \( E^s \cup E^u \)). The diffeomorphism \( f \) is \emph{accessible} if \( AC(x) = M \) for all \( x \), i.e., any two points are connected by an $su$-path. Let \( \Gamma(f) \) denote the set of non-open accessibility classes; then \( f \) is accessible if and only if \( \Gamma(f) = \emptyset \).

\begin{thm}{\cite{2008nil}}\label{su-foliation}
    Let \( f \) be a partially hyperbolic diffeomorphism on an orientable 3-manifold \( M \) with \( NW(f) = M \). Assume \( E^s \), \( E^c \), \( E^u \) are orientable and \( f \) is not accessible. Then one of the following holds:
    \begin{enumerate}
        \item There exists an incompressible torus tangent to \( E^s \oplus E^u \);
        \item There exists a unique nontrivial \( f \)-invariant minimal lamination \( \Gamma(f) \subsetneq M \) tangent to \( E^s \oplus E^u \), extendable to a foliation without compact leaves. Boundary leaves of \( \Gamma(f) \) are periodic and contain dense periodic points;
        \item There exists an \( f \)-invariant foliation tangent to \( E^s \oplus E^u \) without compact leaves.
    \end{enumerate}
\end{thm}

A fundamental class of partially hyperbolic diffeomorphisms consists of time-one maps of Anosov flows, where the flow direction governs intermediate dynamics via the center bundle. We adopt topological Anosov flows as a simplifying framework:

\begin{thm}\cite{IM90, Paternain93}\cite[Theorem 5.9]{BFP23collapsed}\label{flow-defn}
    A non-singular continuous flow $\phi_t: M\rightarrow M$ is a \emph{topological Anosov flow} if and only if it is expansive and preserves a foliation.
\end{thm}

A partially hyperbolic diffeomorphism $f: M\rightarrow M$ on a closed 3-manifold is called a \emph{discretized Anosov flow} if there exists a topological Anosov flow $\phi_t: M\rightarrow M$ and continuous function $\tau:M\rightarrow \mathbb{R}$ satisfying $f(x)=\phi_{\tau(x)}(x)$ for all $x\in M$. When $\phi_t$ is orbit-equivalent to a suspension flow over an Anosov diffeomorphism on $\mathbb{T}^2$, we call $f$ a \emph{discretized suspension Anosov flow}. Detailed treatments appear in \cite{FP_hyperbolic,BFFP1,BFP23collapsed,Martinchich23}.

The following result establishes accessibility for discretized Anosov flows:
\begin{thm}\cite[Theorem C]{FP_hyperbolic}\label{DAF}
    Let $f: M\rightarrow M$ be a partially hyperbolic diffeomorphism on a closed 3-manifold whose fundamental group is not (virtually) solvable. If $f$ is a discretized Anosov flow, then it is accessible.
\end{thm}

Discretized Anosov flows play pivotal roles in 3-dimensional classification programs, particularly well-understood under recurrent conditions when the ambient manifold possesses virtually solvable fundamental group. The non-wandering property precludes periodic tori tangent to either $E^s\oplus E^c$ or $E^c\oplus E^u$. The following classification result is derived from this obstruction.

\begin{thm}\cite{HP15}\label{sol-nil}
    Let $f: M \rightarrow M$ be a partially hyperbolic diffeomorphism of a closed 3-manifold such that $\pi_1(M)$ is virtually solvable but not virtually nilpotent. If either $f$ is dynamically coherent or $NW(f)=M$, then up to a finite iterate, $f$ is a discretized suspension Anosov flow.
\end{thm}

\subsection{Taut Foliations}

We summarize some results on codimension-one foliations, particularly 2-dimensional foliations in 3-manifolds.  

A \emph{Reeb component} of a foliation is a solid torus whose interior is foliated by planes transverse to the core, with all leaves limiting on the boundary torus (which is itself a leaf). A foliation is \emph{Reebless} if it contains no Reeb components. Let \( \widetilde{\mathcal{F}} \) denote the lift of \( \mathcal{F} \) to the universal cover \( \widetilde{M} \). A foliation \( \mathcal{F} \) is \emph{\(\mathbb{R}\)-covered} if the leaf space of \( \widetilde{\mathcal{F}} \) is homeomorphic to \( \mathbb{R} \) (in particular, it is Hausdorff). A foliation is \emph{taut} if there exists a closed transversal intersecting every leaf. For codimension-one foliations:

\begin{thm}{\cite[Corollary 3.3.8, Corollary 6.3.4]{CC00I}}
    Let \( \mathcal{F} \) be a transversely oriented, codimension-one foliation on a compact 3-manifold \( M \). If \( \mathcal{F} \) has no compact leaves, it is taut.
\end{thm}

The next theorem synthesizes properties of taut foliations, drawing from Novikov's work, \cite{Roussarie71}, \cite{Palmeira}, and \cite{Calegari07book}. Analogous results hold for essential laminations \cite{GabaiOertel89}.

\begin{thm}\label{taut}
    Let \( M \) be a closed 3-manifold not finitely covered by \( \mathbb{S}^2 \times \mathbb{S}^1 \). If \( M \) admits a taut foliation \( \mathcal{F} \), then:
    \begin{itemize}
        \item Every leaf of \( \widetilde{\mathcal{F}} \) is a properly embedded plane dividing \( \widetilde{M} \) into two connected half-spaces.
        \item The universal cover \( \widetilde{M} \) is homeomorphic to \( \mathbb{R}^3 \).
        \item Every closed transversal to \( \mathcal{F} \) is homotopically non-trivial.
        \item Transversals to \( \widetilde{\mathcal{F}} \) intersect each leaf at most once.
    \end{itemize}
\end{thm}

The following result was originally proved in \cite{Rosenberg} for $C^2$ foliations and then generalized in \cite{Gabai90} for essential laminations. We will utilize it to obtain leaves with non-trivial stabilizer. 
\begin{thm}\cite{Gabai90}\label{rosenberg}
    Let $\mathcal{F}$ be a foliation or an essential lamination of a closed 3-manifold $M\neq \mathbb{T}^3$. Then, $\mathcal{F}$ must contain a non-planar leaf. In particular, there exists a leaf $F\in \widetilde{\mathcal{F}}$ and $\gamma\in \pi_1(M)\setminus\{id\}$ such that $\gamma(F)=F$, where $\widetilde{\mathcal{F}}$ is the lifted foliation of $\mathcal{F}$ to the universal cover $\widetilde{M}$.
\end{thm}

The following alternative formulation of the theorem above better suits our needs, bypassing the theory of essential laminations:

\begin{thm}\label{rosenberg-lift}
    Let $\mathcal{F}$ be a Reebless foliation in a closed 3-manifold $M$ with non-abelian $\pi_1(M)$, and $\widetilde{\mathcal{F}}$ its lift to the universal cover $\widetilde{M}$. For every non-empty closed $\pi_1(M)$-invariant $\widetilde{\mathcal{F}}$-saturated set $\Lambda\subset \widetilde{M}$, there exists a leaf $F\in\Lambda$ satisfying $\gamma(F)=F$ for some non-trivial $\gamma\in \pi_1(M)$.
\end{thm}

A foliation $\mathcal{F}$ is called \emph{uniform} if every pair of lifted leaves in $\widetilde{M}$ lies within bounded Hausdorff distance.

\begin{thm}\cite[Theorem 1.1]{FP20minimal}
    A uniform Reebless foliation in a closed 3-manifold is $\mathbb{R}$-covered.
\end{thm}

Note that Reeb components contain toroidal boundary leaves, which cannot be made Gromov-hyperbolic through any metric choice. The theorem above then applies to have the following consequence.
\begin{cor}\label{uniform-Gromov-Rcovered}
    A uniform foliation with Gromov hyperbolic leaves in a closed 3-manifold is Reebless and $\mathbb{R}$-covered.
\end{cor}

\subsection{Branching foliation}

Without assuming dynamical coherence, there may be no foliations tangent to the invariant bundles $E^c\oplus E^s$ nor $E^c\oplus E^u$. See \cite{2016example} for an example on $\mathbb{T}^3$ and \cite{BGHP3} for examples on the unit tangent bundle of a higher genus surface. Instead, we can utilize branching foliations introduced by \cite{BI08} as substitutes for center-stable and center-unstable foliations. We present some properties of branching foliations in this subsection. 

\begin{defn}
	A \emph{branching foliation} in a 3-manifold is a collection $\mathcal{W}$ of $C^1$-immersed surfaces, complete under the induced metric and satisfying the following properties:
	\begin{enumerate}
	    \item for each $x\in M$, there is at least one $L\in \mathcal{W}$ containing $x$;
        \item any leaf has no topological crossing with itself;
        \item any pair of leaves has no topological crossing;
        \item each leaf $L(x)\in\mathcal{W}$ through a point $x$ varies continuously with $x$.
	\end{enumerate}
\end{defn}

Note that a branching foliation may be not a real foliation, which allows different leaves to merge. Although, branching foliations could be approximated well by real foliations.

\begin{defn}
	A branching foliation $\mathcal{W}$ is \emph{well-approximated by foliations} if for every arbitrarily small $\epsilon>0$, there is a foliation $\mathcal{W}_{\epsilon}$ by $C^1$ leaves and a continuous map $h_{\epsilon}: M\rightarrow M$ such that the following properties hold:
	\begin{enumerate}
	    \item the angles between tangent spaces of $\mathcal{W}$ and $\mathcal{W}_{\epsilon}$ have a uniform upper bound $\epsilon$;
        \item the $C^0$-distance between $h_{\epsilon}$ and the identity is bounded by $\epsilon$;
        \item the map $h_{\epsilon}$ restricted to each leaf of $\mathcal{W}_{\epsilon}$ is a local diffeomorphism to a leaf of $\mathcal{W}$;
        \item every leaf $L\in \mathcal{W}$ has a unique associated leaf $L_{\epsilon}\in \mathcal{W}_{\epsilon}$ such that $h_{\epsilon}(L_{\epsilon})=L$.
	\end{enumerate}
\end{defn}

We say that the foliation $\mathcal{W}_{\epsilon}$ in the definition above for each $\epsilon>0$ is a \emph{well-approximated foliation} of the branching foliation $\mathcal{W}$. A priori, the surjective map $h_{\epsilon}$ is only a local diffeomorphism on each leaf of $\mathcal{W}_{\epsilon}$, instead of a global diffeomorphism. Note that $h_{\epsilon}$ lifts to a uniformly bounded diffeomorphism in the universal cover $\widetilde{M}$.


We will use the following result on the existence of branching foliations preserved by a partially hyperbolic diffeomorphism.
\begin{thm}\cite{BI08}
    Let $f: M\rightarrow M$ be a partially hyperbolic diffeomorphism of a closed 3-manifold. Assume that $M$ and $E^{\sigma}, \sigma=s, c,u$ are all oriented and the orientations are preserved by $Df$. Then, there exist $f$-invariant branching foliations $\mathcal{W}^{cs}$ and $\mathcal{W}^{cu}$ tangent to $E^s\oplus E^c$ and $E^c\oplus E^u$, respectively, which are well-approximated by foliations.
\end{thm}

As demonstrated in \cite{BI08}, branching foliations exert comparable influence on transversals to well-approximated foliations. The following lemma adapts arguments from \cite{BI08}, with proof outlined for completeness:

\begin{lem}\label{branch-Reebless}
    Let $f: M\rightarrow M$ be a partially hyperbolic diffeomorphism on a closed 3-manifold with branching foliation $\mathcal{W}^{cs}$ tangent to $E^c\oplus E^s$. Then $\mathcal{W}^{cs}$ admits no Reeb components, and $\mathcal{W}^{cs}_{\epsilon}$ forms a Reebless foliation for sufficiently small $\epsilon>0$. Moreover, every transversal to $\widetilde{\mathcal{W}}^{cs}$ in $\widetilde{M}$ intersects leaves at most once.
\end{lem}

\begin{proof}
    The unstable bundle $E^u$ integrates uniquely to a non-compact leafwise transverse foliation $\mathcal{F}^u$. A Reeb component in $\mathcal{W}^{cs}$ would induce closed leaves in $\mathcal{F}^u$ via \cite[Lemma 2.2]{BI08}, contradicting partial hyperbolicity. This conclusion extends to $\mathcal{W}^{cs}_{\epsilon}$-approximations since they maintain $\mathcal{F}^u$-transversality for small $\epsilon$. The final property follows directly from Theorem~\ref{taut}.
\end{proof}

The following lemma is a consequence of \cite[Lemma B.2]{BFFP2}, see \cite[Section 3]{BFFP2} for the meaning of minimal branching foliations. 

\begin{lem}\label{minimal-branching}
	Let $\mathcal{W}^{cs}$ and $\mathcal{W}^{cu}$ be $f$-invariant branching foliations without compact leaves tangent to $E^s\oplus E^c$ and $E^c\oplus E^u$, respectively. If $NW(f)=M$, then both $\mathcal{W}^{cs}$ and $\mathcal{W}^{cu}$ are minimal.
\end{lem}

We refer the reader to \cite{BFFP2, BFP23collapsed} for a comprehensive account for branching foliations and their leaf spaces.

\section{Asymptotic behaviors on the boundary at infinity}\label{section-asymptotic}

In this section, we establish general results for transverse codimension-one foliations in closed 3-manifolds $M$ independent of dynamical constraints. Specifically, we require $M$ to lack finite covers by $\mathbb{S}^2 \times \mathbb{S}^1$ except when analyzing partially hyperbolic systems. Our focus lies in characterizing asymptotic behaviors of a pair of transverse foliations and their intersections near ideal boundaries.

\subsection{Boundary at infinity and limit points}\label{subsection-topology}

We first describe leafwise ideal boundaries and their global topology. Let $\widetilde{M}$ denote the universal cover of the given closed manifold $M$, and $\widetilde{\mathcal{F}}$ the lift of a foliation $\mathcal{F}$ from $M$ to $\widetilde{M}$.

A 2-dimensional foliation $\mathcal{F}$ of a 3-manifold $M$ comprises \emph{Gromov hyperbolic leaves} if there exists a uniform quasi-isometric map sending $(F, d_F)$ to $(\mathbb{H}^2, d_{\mathbb{H}^2})$ for each $F\in \widetilde{\mathcal{F}}$, where $d_F$ is the induced path metric of a given Riemannian metric lifted to the universal cover $\widetilde{M}$ and $d_{\mathbb{H}^2}$ is the standard hyperbolic metric of constant negative curvature. This definition aligns with standard Gromov hyperbolicity due to Candel's uniformization theorem \cite{Candel93}. For any metric space $(X, d_X)$ quasi-isometric to a hyperbolic space $(\mathbb{H}^n, d_{\mathbb{H}^n})$, a quasi-isometric map extends to a boundary homeomorphism, yielding a canonical compactification $X\cup \partial_{\infty}X$ identified with $\mathbb{H}^n\cup\partial\mathbb{H}^n$ \cite{Gromov87, Thurston88}. The topology of the boundary $\partial_{\infty}X$ is given by the natural topology of $\partial\mathbb{H}^n$ through a homeomorphism. Thus, each Gromov hyperbolic leaf $F\in \widetilde{\mathcal{F}}$ compactifies to $F\cup\partial_{\infty}F$ identified with the hyperbolic Poincar\'e disc $\mathbb{D}^2=\mathbb{H}^2\cup\partial\mathbb{H}^2$. We call $\partial_{\infty}F$ the \emph{ideal boundary} or \emph{ideal circle} of $F$.

The union of all ideal circles $\mathcal{A}_{\mathcal{F}} := \bigcup_{F\in \widetilde{\mathcal{F}}} \partial_\infty F$ forms the \emph{tubulation at infinity}, termed the \emph{cylinder at infinity} when $\mathcal{F}$ is $\mathbb{R}$-covered. By adapting Candel's metric \cite{Candel93}, one induces a natural topology on $\mathcal{A}_{\mathcal{F}}$ through identifications of ideal circles with unit tangent circles at basepoints along connected transversals to $\widetilde{\mathcal{F}}$ \cite{Calegari00,Fenley02}. Leafwise negative curvature ensures bijective correspondence between ideal circles and unit tangent circles in this metric.

To make a unitary topology of tubulations at infinity for two foliations with Gromov hyperbolic leaves that is compatible with the topology of $\widetilde{M}$, we need to fix a metric which might make no sense to have leaf curvature. Instead of having bijective correspondence between each ideal circle and a unit tangent circle, we identify each ideal circle with a topological circle given by monotonely quotienting some connected intervals of a unit tangent circle. See for instance \cite[Proposition 4.3]{Fenley92} and \cite[Section 3]{FP23intersection}. Crucially, the resultant tubulation topology on $\mathcal{A}_{\mathcal{F}}$ for any foliation $\mathcal{F}$ with Gromov hyperbolic leaves remains compatible with the given topology of $\widetilde{M}$. Consequently, the universal covering $\widetilde{M}$ becomes homeomorphic to a solid connected tubulation for any such foliation, which is particularly a solid cylinder for $\mathbb{R}$-covered foliations. Under these considerations, we define a \emph{geodesic} in each leaf $F\in \widetilde{\mathcal{F}}$ as a bi-infinite curve minimizing leafwise path distances between all point pairs, with analogous definitions for \emph{geodesic segments} and \emph{geodesic rays}.

The topology described above remains invariant under Riemannian metric changes on $M$ \cite[Lemma 3.7]{FP23intersection}. This equivalence demonstrates that $\mathcal{A}_{\mathcal{F}}$'s topology coincides with Candel's metric-induced topology \cite{Calegari00,Fenley02} (see \cite[Section 7.2]{Calegari07book}).

For any connected ideal interval $I \subset \partial_\infty F$ and point $x \in F$, we introduce the \emph{wedge} $W_x(I)$ as the maximal union of geodesic rays from $x$ terminating in $I$. The \emph{visual measure} $V_x(I)$ denotes the angular span of $W_x(I)$. When identifying ideal circles with unit tangent circles via surjective maps (potentially involving non-trivial quotients), an ideal interval may produce distinct maximal and minimal wedges, resulting in different visual measures. Our analysis exclusively employs maximal wedges and their associated visual measures.

Given a point $x\in \widetilde{M}$, the visual measure at $x$ alternatively arises from the visual metric on the ideal boundary $\partial_{\infty}F$ using the Gromov product at $x$ with some parameters. The visual metric is well-defined up to H\"{o}lder equivalence on the ideal circle $\partial_{\infty}F$ for each $F\in \widetilde{\mathcal{F}}$, see \cite[Chapter III.H.3]{BH13metric}. For connected $I$, this measure corresponds to the endpoint visual metric at $x$, regardless of the choices of parameters. Moreover, these visual measures vary continuously across points in $\widetilde{M}$ due to geodesic convergence under point sequences. The maximal wedge visual measure exhibits upper semicontinuity relative to ideal intervals, while the minimal version shows lower semicontinuity.

We investigate the ideal limit sets of intersection foliation $\widetilde{\mathcal{G}} = \widetilde{\mathcal{F}}_1 \cap \widetilde{\mathcal{F}}_2$ within ideal boundaries for transverse foliations $\mathcal{F}_1$ and $\mathcal{F}_2$. For any leaf $F \in \widetilde{\mathcal{F}}_i$, let $\widetilde{\mathcal{G}}_F$ denote the restricted subfoliation of $\widetilde{\mathcal{G}}$ on $F$ ($i=1,2$).

\begin{lem}\label{singlelimit}\cite[Theorem 4.1]{FU1}
    Let $\mathcal{F}_1, \mathcal{F}_2$ be transverse 2-dimensional foliations with non-compact Gromov hyperbolic leaves in a closed 3-manifold $M \neq \mathbb{T}^3$, and let $\mathcal{G} = \mathcal{F}_1 \cap \mathcal{F}_2$ be their one-dimensional intersection subfoliation. If $\mathcal{F}_1$ (resp. $\mathcal{F}_2$) is minimal and $\mathbb{R}$-covered, then every ray in $\widetilde{\mathcal{G}}_F$ converges to a unique point in $\partial_\infty F$ for all $F \in \widetilde{\mathcal{F}}_1$ (resp. $F \in \widetilde{\mathcal{F}}_2$).
\end{lem}

The condition $M \neq \mathbb{T}^3$ in \cite[Theorem 4.1]{FU1} ensures non-planar leaves via Rosenberg's theorem \cite{Rosenberg} (see Theorem~\ref{rosenberg}), which asserts that 3-manifolds foliated by planes must be tori. Our framework similarly excludes tori to guarantee non-planar leaf existence through this mechanism.

\subsection{Coherent leaf-wise structures}\label{subsection-coherent}

Given a codimension-one foliation $\mathcal{G}$ of a space $X$ (open or compact), the \emph{leaf space} of $\mathcal{G}$ is the quotient space obtained by collapsing each leaf to a point, endowed with the quotient topology. This leaf space forms a simply connected, separable one-dimensional manifold that need not be Hausdorff \cite{Calegari07book}. For codimension-one foliations, Hausdorff leaf spaces characterize $\mathbb{R}$-covered structures. When $\mathcal{G}$ possesses a non-Hausdorff leaf space, certain leaves correspond to points lacking Hausdorff neighborhoods; these are termed \emph{non-separated leaves}. Equivalently, a leaf $l$ is non-separated if another distinct leaf $l'$ and a sequence $\{r_n\}$ exist where $r_n$ accumulates to both $l$ and $l'$ as $n\to\infty$.

The following result, established for transverse foliations in \cite{FP23_transverse}, adapts to our framework without requiring Gromov hyperbolic leaves:
\begin{prop}\label{intersection}
    Let $\mathcal{F}_1, \mathcal{F}_2$ be transverse taut 2-dimensional foliations in a closed 3-manifold with intersection subfoliation $\mathcal{G} = \mathcal{F}_1 \cap \mathcal{F}_2$. 
    \begin{itemize}
        \item If $\mathcal{F}_2$ is $\mathbb{R}$-covered, then for any $F \in \widetilde{\mathcal{F}}_1$, the leaf space of $\widetilde{\mathcal{G}}_F$ is non-Hausdorff if and only if there exists $E \in \widetilde{\mathcal{F}}_2$ such that $F\cap E$ is not connected. Moreover, non-separated leaves $s_1, s_2 \in \widetilde{\mathcal{G}}_F$ must lie within a common $E \in \widetilde{\mathcal{F}}_2$.
        
        \item If both $\mathcal{F}_1$ and $\mathcal{F}_2$ are $\mathbb{R}$-covered with $\mathcal{F}_1$ minimal, then Hausdorff $\widetilde{\mathcal{G}}$-leaf space for any $F \in \widetilde{\mathcal{F}}_1$ implies all leaves in $\widetilde{\mathcal{F}}_1$ and $\widetilde{\mathcal{F}}_2$ have Hausdorff $\widetilde{\mathcal{G}}$-leaf spaces.
    \end{itemize}
\end{prop}
\begin{proof}
    Assume that a leaf $F\in \widetilde{\mathcal{F}}_1$ has disconnected intersection with some leaf $E\in \widetilde{\mathcal{F}}_2$, yielding distinct leaves $l_1, l_2 \in \widetilde{\mathcal{G}}_F$ within $F \cap E$. If $\widetilde{\mathcal{G}}_F$ has Hausdorff leaf space, a transversal in $F$ intersecting both $l_1$ and $l_2$ would produce a $\widetilde{\mathcal{F}}_2$-transversal joining $E$ at two points - contradicting Theorem~\ref{taut}.
    
    For non-Hausdorff $\widetilde{\mathcal{G}}_F$-leaf spaces, consider non-separated leaves $s_1, s_2 \in \widetilde{\mathcal{G}}_F$ and a sequence of leaves $r_n\in \widetilde{\mathcal{G}}_F$ accumulating simultaneously to $s_1, s_2$. Let $E_n \in \widetilde{\mathcal{F}}_2$ denote leaves containing $r_n$. For large $n\neq m$, the leaves $E_n$ and $E_m$ are distinct. Indeed, if $x_n \in r_n$ and $x_m \in r_m$ shared $E_n$ for $x_n, x_m$ close enough, there would be a $\widetilde{\mathcal{F}}_2$-transversal intersecting $E_n$ twice, violating Theorem~\ref{taut}. Thus, up to a subsequence, distinct $E_n$ accumulate to $\widetilde{\mathcal{F}}_2$-leaves through $s_1$ and $s_2$. The $\mathbb{R}$-covered property of $\mathcal{F}_2$ forces a single $E \in \widetilde{\mathcal{F}}_2$ containing both $s_1$ and $s_2$, making $F \cap E$ non-connected through saturation by $\widetilde{\mathcal{G}}$-leaves.

    When $\widetilde{\mathcal{G}}_F$ has non-Hausdorff leaf space, the corresponding $E \in \widetilde{\mathcal{F}}_2$ necessarily induces a non-Hausdorff $\widetilde{\mathcal{G}}_E$-leaf space. Should $\widetilde{\mathcal{G}}_E$ instead be Hausdorff, repeating the argument would yield a transversal to both $\widetilde{\mathcal{G}}_E$ and $\widetilde{\mathcal{F}}_1$ intersecting $s_1$ and $s_2$ - impossible since these lie in the same $\widetilde{\mathcal{F}}_1$-leaf, contradicting Theorem~\ref{taut}. Symmetrically, non-Hausdorff $\widetilde{\mathcal{G}}_E$ for any $E \in \widetilde{\mathcal{F}}_2$ implies existence of $F \in \widetilde{\mathcal{F}}_1$ with non-Hausdorff $\widetilde{\mathcal{G}}_F$-leaf space.
    
    Let $F \in \widetilde{\mathcal{F}}_1$ be a leaf having disconnected intersection with some $E \in \widetilde{\mathcal{F}}_2$. By $\mathcal{F}_1$-minimality, for any $L \in \widetilde{\mathcal{F}}_1$ there exists a deck transformation $\rho$ with $\rho(L)$ arbitrarily close to $F$ such that $\rho(L) \cap E$ remains disconnected. Specifically, the leaf $F$ divides $\widetilde{M}$ into transverse orientation consistent components $F^+$ and $F^-$. Transversality of $\widetilde{\mathcal{F}}_1$ and $\widetilde{\mathcal{F}}_2$ ensures nearby $\widetilde{\mathcal{F}}_1$-leaves also intersect $E$. The disconnected intersection $F\cap E$ remains valid for nearby leaves of $\widetilde{\mathcal{F}}_1$ within either $F^+$ or $F^-$. Through minimality, choose deck transformations $\rho^\pm$ positioning $\rho^\pm(L)$ in $F^\pm$ near $F$. Selecting $\rho = \rho^+$ or $\rho^-$ then propagates the disconnected intersection property to $\rho(L)\cap E$. It follows that $L$ admits disconnected intersection with some leaf of $\widetilde{\mathcal{F}}_2$.
    
    Thus, non-Hausdorff $\widetilde{\mathcal{G}}_F$-leaf space for one $F \in \widetilde{\mathcal{F}}_1$ necessitates non-Hausdorff leaf spaces for all $\widetilde{\mathcal{F}}_1$-leaves. Combining these arguments completes the proof.
\end{proof}

\begin{rmk}
	For general transverse Reebless foliations $\mathcal{F}_1$ and $\mathcal{F}_2$, non-separated leaves $l_1, l_2$ of $\widetilde{\mathcal{G}} = \widetilde{\mathcal{F}}_1 \cap \widetilde{\mathcal{F}}_2$ within a leaf $F \in \widetilde{\mathcal{F}}_1$ need not remain non-separated in their containing leaf $E \in \widetilde{\mathcal{F}}_2$. In such cases, the $\widetilde{\mathcal{G}}_E$-leaf space between $l_1$ and $l_2$ cannot be homeomorphic to a real interval - it must contain non-separated leaves separating $l_1$ from $l_2$, as no $\widetilde{\mathcal{G}}_E$-transversal in $E$ can intersect both. See \cite[Section 7]{FP23_transverse} for illustrative examples.
\end{rmk}

The \emph{ideal limit set} of $\widetilde{\mathcal{G}}_F$ refers to the union of all limit points of $\widetilde{\mathcal{G}}_F$-leaves in $\partial_\infty F$, which is always a non-empty subset of $\partial_{\infty}F$. We say $F \in \widetilde{\mathcal{F}}_i$ ($i=1,2$) has a \emph{degenerate limit set} if this limit set of $\widetilde{\mathcal{G}}_F$ reduces to a single ideal point in $\partial_{\infty}F$.

\begin{lem}\label{nondense=single}\cite[Corollary 4.14]{FU1}
    Let $\mathcal{F}_1, \mathcal{F}_2$ be transverse 2-dimensional foliations with non-compact leaves in a closed 3-manifold, and $\mathcal{G} = \mathcal{F}_1 \cap \mathcal{F}_2$ be their one-dimensional intersecting subfoliation. Assume that $\mathcal{F}_1$ is a minimal $\mathbb{R}$-covered foliation with Gromov hyperbolic leaves. If the ideal limit set of $\widetilde{\mathcal{G}}_L$ is dense in $\partial_{\infty}L$ for some $L\in\widetilde{\mathcal{F}}_1$, then each leaf $F\in \widetilde{\mathcal{F}}_1$ has dense limit set in $\partial_{\infty}F$. If the ideal limit set of $\widetilde{\mathcal{G}}_L$ is not dense in $\partial_{\infty}L$ for some $L\in\widetilde{\mathcal{F}}_1$, then each leaf $F\in \widetilde{\mathcal{F}}_1$ has degenerate limit set in $\partial_{\infty}F$. 
\end{lem}

\subsection{Transverse continuity}
We establish the continuous variation of ideal points of $\widetilde{\mathcal{G}}_F$ under the topology of $\mathcal{A}_{\mathcal{F}_1}$ as $F$ ranges through $\widetilde{\mathcal{F}}_1$.

Let us first present the notion of non-expanding directions:

\begin{defn}[Contracting and Non-expanding Direction]\label{def_contracting}
    Let $F$ be a leaf of $\widetilde{\mathcal{F}}$ containing a point $x$, and $\gamma$ a geodesic ray in $F$ starting at $x$ with initial tangent vector $v$. We say $\gamma$ (or $v$) is a \emph{contracting direction} if:
    \begin{itemize}
        \item There exists a transversal $\tau$ to $\widetilde{\mathcal{F}}$ through $x$;
        \item For every leaf $L \in \widetilde{\mathcal{F}}$ intersecting $\tau$, the distance $d(\gamma(t), L) \to 0$ as $t \to \infty$.
    \end{itemize}
    
    Similarly, $\gamma$ (or $v$) is called an \emph{$\epsilon$-non-expanding direction} if:
    \begin{itemize}
        \item There exists $\epsilon > 0$ and a transversal $\tau_\epsilon$ to $\widetilde{\mathcal{F}}$ through $x$;
        \item For every leaf $L \in \widetilde{\mathcal{F}}$ intersecting $\tau_\epsilon$, the distance $d(\gamma(t), L) \leq \epsilon$ for all $t \geq 0$.
    \end{itemize}
\end{defn}

Thurston formulates the abundance of contracting and non-expanding directions in hyperbolic leaves of codimension-one foliations. The work of \cite{Calegari00,Fenley02} features stronger hypotheses and conclusions compared to Thurston's original result, specifically guaranteeing density of contracting directions between arbitrary leaf pairs. While Thurston's result holds in arbitrary dimensions, we specialize to 3-manifolds for simplicity.

\begin{thm}[Thurston]\label{Thurston}
    Let $\mathcal{F}$ be a codimension-one foliation with hyperbolic leaves in a closed 3-manifold $M^3$. For any $\epsilon > 0$ and leaf $F \in \widetilde{\mathcal{F}}$, the $\epsilon$-non-expanding directions are dense at every $x \in F$. If no holonomy-invariant transverse measure supported on $\pi(F)$ exists for any $F \in \widetilde{\mathcal{F}}$, then each $\epsilon$-non-expanding direction becomes a contracting direction.
    Furthermore, when $\pi(F)$ is non-compact, the transversal $\tau$ from Definition~\ref{def_contracting} may be chosen with $x$ in its interior.
\end{thm}

We say that a foliation (or a lamination) $\widetilde{\mathcal{F}}$ has \emph{transverse continuity} with respect to a transverse foliation $\widetilde{\mathcal{W}}$ if for any sequence of leaves $L_i\in \widetilde{\mathcal{F}}$ converging to a leaf $L\in \widetilde{\mathcal{F}}$ and any sequence of leaves $l_i\in L_i\cap E$ accumulating to a leaf $l\in L\cap E$ in a common $E\in \widetilde{\mathcal{W}}$, the ideal points $l_i^+, l_i^-\in \partial_{\infty}L_i$ converge to the ideal points $l^+, l^-\in \partial_{\infty}L$, respectively, in the topology of the tubulation at infinity $\mathcal{A}_{\mathcal{F}}$.
We also say that ideal points in $\mathcal{A}_{\mathcal{F}}$ for leaves of $\widetilde{\mathcal{F}}\cap \widetilde{\mathcal{W}}$ are transversely continuous with respect to the foliation $\widetilde{\mathcal{F}}$. 
For the sake of use, we will present another equivalent definition below as an interpretation of transverse continuity.

Given a leaf $F\in \widetilde{\mathcal{F}}$, let $I_F\subset \partial_{\infty}F$ be an ideal interval bounded by two $\epsilon$ non-expanding ideal points for some $\epsilon$. For any geodesic rays $\alpha_F, \beta_F\subset F$ towards these two $\epsilon$ non-expanding ideal points, all points on these two rays lie within $\epsilon$-neighborhoods of a nearby leaf $L\in \widetilde{\mathcal{F}}$. Through transversals in leaves of $\widetilde{\mathcal{W}}$ to $\widetilde{\mathcal{F}}$, there are two curves $\alpha_L, \beta_L\subset L$ defined by $\alpha_F, \beta_F$, respectively, arbitrarily close to geodesic rays in $L$ if $\epsilon$ is chosen small enough. Thus, they uniquely determine two corresponding ideal points in $\partial_{\infty}L$ that bound an ideal interval $I_L\subset \partial_{\infty}L$ corresponding to $I_F$. We call the ideal interval $I_L$ the \emph{projection} of $I_F$ on the leaf $L$.

Analogous to contracting directions, we can construct a collection of embedded curves in $\mathcal{A}_{\mathcal{F}}$ through $\epsilon$-non-expanding directions. Given a geodesic ray $\gamma$ associated with an  $\epsilon$-non-expanding direction in a leaf $L \in \widetilde{\mathcal{F}}$ along a transversal $\tau_{\epsilon}$, the corresponding curve $\gamma_F\subset F$ defined by $\gamma$ uniquely determines an ideal point $\xi_F\in \partial_{\infty}F$ for any nearby leaf $F \in \widetilde{\mathcal{F}}$. When $\mathcal{F}$ is a minimal foliation consisting of non-compact leaves, the leaf $F$ given above can be taken as any leaf of $\widetilde{\mathcal{F}}$. We define an \emph{$\epsilon$-marker} as the union $\bigcup_{F \in \widetilde{\mathcal{F}}} \xi_F$ for each $\epsilon$-non-expanding direction. The $\epsilon$-markers form embedded curves in $\mathcal{A}_{\mathcal{F}}$ with pairwise disjointness or coincidence \cite{Fenley02,Calegari00}. Note that $\epsilon$-markers can also be defined analogously without $\mathbb{R}$-covered or minimal requirements for $\mathcal{F}$ \cite{Calegari07book}. We present the tools here in a way for our convenience.

Using the denseness of $\epsilon$ non-expanding ideal points (see Theorem \ref{Thurston}) and the continuity of $\epsilon$-markers (see \cite{Fenley02}), we have the following equivalent definition of transverse continuity.

\begin{defn}
	Let $l_F$ be a ray in a leaf of $F\cap E$ for any $F\in \widetilde{\mathcal{F}}$ and $E\in \widetilde{\mathcal{W}}$. Consider any sufficiently small $\epsilon>0$ and any ideal interval $I_F\subset \partial_{\infty}F$ containing the ideal point of $l_F$, say $l_F^+$, in the interior so that two endpoints of $I_F$ are $\epsilon$ non-expanding ideal points. We say that the foliation $\widetilde{\mathcal{F}}$ has \emph{transverse continuity} with respect to $\widetilde{\mathcal{W}}$ if for any $L\in \widetilde{\mathcal{F}}$ sufficiently close to $F$, there is a ray $l_L$ of a leaf of $L\cap E$ such that its ideal point $l_L^+$ is contained in the projection $I_L\subset \partial_{\infty}L$.
\end{defn}

We provide a sufficient condition for transverse continuity of a foliation without compact leaves. Notably, the $\mathbb{R}$-covered and minimal properties are not required provided that each ray accumulates at a single ideal point in its underlying boundary.

\begin{prop}\label{transverse-continuity-id}
	Let $\mathcal{F}_1, \mathcal{F}_2$ be two transverse foliations by Gromov hyperbolic leaves and $\mathcal{G}$ be the intersection foliation. Assume that every ray in $\widetilde{\mathcal{G}}_F$ and $\widetilde{\mathcal{G}}_E$ converges to a single ideal point of $\partial_{\infty}F$ and $\partial_{\infty}E$, respectively, for any $F\in \widetilde{\mathcal{F}}_1$ and $E\in \widetilde{\mathcal{F}}_2$. If the foliation $\widetilde{\mathcal{G}}_E$ has Hausdorff leaf space for each leaf $E\in \widetilde{\mathcal{F}}_2$, then the ideal points $l^{\pm}(x)$ associated to $\widetilde{\mathcal{G}}_E$-leaves vary continuously with $x\in E$. Moreover, the foliation $\widetilde{\mathcal{F}}_1$ has transverse continuity with respect to $\widetilde{\mathcal{F}}_2$.
\end{prop}
\begin{proof}
	Suppose that there is a sequence of points $(x_n)_{n\in\mathbb{N}}$ in $E\in \widetilde{\mathcal{F}}_2$ converging to a point $x\in E$ such that the ideal point $l(x_n)^+$ does not converge to $l(x)^+$ under the topology of $\partial_{\infty}E$. It turns out that the ideal limit set of $\widetilde{\mathcal{G}}_E$ is not equal to a single ideal point of $\partial_{\infty}E$. By Lemma \ref{nondense=single}, this ideal limit set must be dense in $\partial_{\infty}E$. Then there is at least one sequence of points $p_n\in l(x_n)$ in the ray from $x_n$ to the ideal point $l(x_n)^+$ converging to a point not contained in $l(x)$. Otherwise, by the denseness of ideal limit set, all points in $l(x_n)$ converge to points in $l(x)$, as well as ideal points. Denote by $p\in E$ the accumulation point of $p_n$. As leaves of $\widetilde{\mathcal{G}}_E$ accumulate in leaves of $\widetilde{\mathcal{G}}_E$, the leaf $l(p)\in \widetilde{\mathcal{G}}_E$ is also accumulated by $l(x_n)$. Then $l(p)$ is not separated from $l(x)$. Thus, the leaf space of $\widetilde{\mathcal{G}}_E$ is not Hausdorff, which is a contradiction. 
	
	Consider any leaf $L\in \widetilde{\mathcal{F}}_1$ and any leaf $l_L\in\widetilde{\mathcal{G}}_L$.
	Let $l_L^{\pm}\in \partial_{\infty}L$ be the ideal points of $l_L$. We consider a parameterization $l_L(t): (-\infty, +\infty)\rightarrow l_L$ such that $l_L(t)$ goes to $l_L^+$ and $l_L^-$ as $t\rightarrow +\infty$ and $t\rightarrow -\infty$, respectively. Let $E\in \widetilde{\mathcal{F}}_2$ be the leaf through $l_L$. By transversality of $\widetilde{\mathcal{F}}_1$ and $\widetilde{\mathcal{F}}_2$, we can pick an $\epsilon>0$ small enough such that if a leaf $F\in \widetilde{\mathcal{F}}_1$ contains a point $x\in F$ of $d(x,z)\leq \epsilon$ for some point $z\in L$, then any leaf $E(x)\in \widetilde{\mathcal{F}}_2$ through $x$ intersects $L$ in a point $y\in E(x)\cap L$ with $d_L(y,z)\leq \delta$ for a small constant $\delta>0$ depending on $\epsilon$. Let $I_L\subset \partial_{\infty}L$ be a small ideal interval containing $l_L^+$ in the interior such that its two endpoints are $\epsilon$ non-expanding ideal points. Suppose that $F_n\in \widetilde{\mathcal{F}}_1, n\in \mathbb{N},$ is a sequence of leaves converging to $L$ such that for the leaf $l_{F_n}\in \widetilde{\mathcal{G}}_{F_n}$ in the intersection $F_n\cap E$, the ideal point $l_{F_n}^+\in\partial_{\infty}F_n$ is not contained in the projection $I_{F_n}\subset\partial_{\infty}F_n$ for each $n\in \mathbb{N}$. The leaf $l_{F_n}\in \widetilde{\mathcal{G}}_{F_n}$ is unique by Proposition~\ref{intersection}, and it converges to the leaf $l_L$ as $n$ goes to infinity. 
	
	As the leaf space of $\widetilde{\mathcal{G}}_E$ is Hausdorff and each leaf is oriented, up to reparameterizations, the point $l_{F_n}(t)$ converges to $l_L(t)$ as $n\rightarrow +\infty$ for each $t\in \left[0, +\infty\right)$. Denote by $\alpha_L, \beta_L\subset L$ two geodesic rays from the point $l_L(0)$ to two endpoints of $I_L$, respectively. Let $\alpha_{F_n}\subset F_n$ be a curve of bounded distance $d_H(\alpha_{F_n}, \alpha_L)\leq \epsilon$ and $\beta_{F_n}\subset F_n$ be a curve of bounded distance $d_H(\beta_{F_n}, \beta_L)\leq \epsilon$ for each $n\in \mathbb{N}$. Since $l_L^+$ is in the interior of $I_L$, for the constant $\delta>0$ given above, there is a large $R_{\delta}>0$ such that for any $t> R_{\delta}$, $l_L(t)$ is $\delta$-away from $\alpha_L$ and $\beta_L$ with respect to the distance $d_L$. Note that $l_{F_n}^+$ is not contained in the ideal interval $I_{F_n}$ for each $n\in \mathbb{N}$. Then, for a sufficiently large $N\in \mathbb{N}$, the ray $l_{F_N}(R_{\delta}, +\infty)$ has to intersect either $\alpha_{F_N}$ or $\beta_{F_N}$ at some moment and then never comes back again. Without loss of generality, we assume that $l_{F_N}$ intersects $\alpha_{F_N}$ in a point $l_{F_N}(T)$ for $T>R_{\delta}$. It implies that $d(l_{F_N}(T), \alpha_L)\leq \epsilon$, and thus there is a point $y\in E\cap L$ of distance $d_L(y, \alpha_L)\leq \delta$ by the definition of $\delta$. Therefore, the point $y$ is not contained in the leaf $l_L$, which implies that $E\cap L$ contains at least two stable leaves. This is a contradiction by Proposition~\ref{intersection}, and thus we finish the proof.
\end{proof}

\section{Bounded distance from geodesics to curves}\label{section-distance}

In this section, we analyze transverse minimal $\mathbb{R}$-covered foliations $\mathcal{F}_1$ and $\mathcal{F}_2$ with non-compact Gromov hyperbolic leaves in a closed 3-manifold $M \neq \mathbb{T}^3$. Let $\mathcal{G} = \mathcal{F}_1 \cap \mathcal{F}_2$ denote their one-dimensional intersection foliation. As established in Lemma~\ref{singlelimit}, every ray in $\widetilde{\mathcal{G}}_F$ accumulates at a unique ideal point in $\partial_\infty F$ for all $F \in \widetilde{\mathcal{F}}_1$ and $F \in \widetilde{\mathcal{F}}_2$.

Under the topology of $\mathcal{A}_{\mathcal{F}_1}$, geodesics between fixed ideal points in $\partial_\infty F$ remain within bounded Hausdorff distance for each $F \in \widetilde{\mathcal{F}}_1$. In the subsequent, we will disregard the specific choice of geodesics sharing both endpoints. Given a leaf $l\in \widetilde{\mathcal{G}}_F$ for any $F\in \widetilde{\mathcal{F}}_1$, we denote by $l^+, l^-\in\partial_{\infty}F$ its ideal points in the ideal circle. We refer to $l^*$ as an arbitrary geodesic in the same leaf of $\widetilde{\mathcal{F}}_1$ as $l$ joining $l^+$ and $l^-$. We use the same notation if $l$ is a segment or a ray in a leaf of $\widetilde{\mathcal{G}}$.

Given a point $x\in F$ for some leaf $F\in\widetilde{\mathcal{F}}_1$, any curve $c\in L$ determines a unique ideal interval $I_{x, c}$ in $\partial_{\infty}L$ by geodesic rays starting at $x$ through the points of $c$. We say the ideal interval $I_{x,c}$ defined in such a way is the \emph{shadow} of the curve $c$ at $x$.

We devote this section to reveal a relation between leaves of $\widetilde{\mathcal{G}}$ and their associated geodesics sharing the same endpoints. 

\subsection{Uniformly bounded distance}

For any two points $x, y \in \widetilde{M}$, we denote by $d(x,y)$ the distance of $x$ and $y$ in the universal cover $\widetilde{M}$, and $d_F(x, y)$ the induced distance in the leaf $F$ if $x, y \in F$. Let us introduce a lemma for general $\mathbb{R}$-covered foliations. In particular, it adapts to any sublamination of a $\mathbb{R}$-covered foliation with a possibly smaller constant $B$.
\begin{lem}\label{AB}\cite[Proposition 2.1]{Fenley92}
	If $\mathcal{F}$ is a $\mathbb{R}$-covered (branching) foliation in a closed manifold, then for any $A>0$ there is $B>0$ depending only on $A$ so that for any pair of points $x, y\in \widetilde{M}$ in a leaf $F\in\widetilde{\mathcal{F}}$ with $d(x, y)<A$, we have $d_F(x, y)<B$.
\end{lem}

We first show a result under the transverse continuity condition. It provides a neighborhood of uniform size for every curve of $\widetilde{\mathcal{G}}$ to contain the associated geodesic curve.

\begin{prop}\label{uniform_neighborhood_id}
	Assume that there is a leaf $l_0\in \widetilde{\mathcal{G}}_{L_0}$ for some $L_0\in\widetilde{\mathcal{F}}_1$ such that $l_0^+\neq l_0^-$ and $\widetilde{\mathcal{F}}_1$ has transverse continuity with respect to $\widetilde{\mathcal{F}}_2$. Then there is a uniform constant $C>0$ so that for any $L\in \widetilde{\mathcal{F}}_1$ and any $l\in \widetilde{\mathcal{G}}_L$, the geodesic $l^*$ is contained in the $C$-neighborhood of $l$. The same statement holds for any ray and any segment in a leaf of $\widetilde{\mathcal{G}}_L$.
\end{prop}
\begin{proof}
	We will only prove it for segments, since the same argument extends directly to leaves and rays. Suppose to the contrary that for any $n\in\mathbb{N}$, there is $L_n\in \widetilde{\mathcal{F}}_1$ and a segment $s_n$ in a leaf of $\widetilde{\mathcal{G}}_{L_n}$ satisfying that the associated geodesic segment $s_n^*$ is not entirely contained in the $2n$-neighborhood of $s_n$. It implies that there is a point $x_n\in L_n$ such that the disk $D(x_n, n)$ in $L_n$ centered at $x_n$ of radius $n$ is contained in an open connected region enclosed by $s_n$ and $s_n^*$. By the minimality of $\mathcal{F}_1$, up to a subsequence and the action of a sequence of deck transformations, we assume that the point $x_n\in L_n$ converges to a point $x\in F\in\widetilde{\mathcal{F}}_1$ as $n$ goes to infinity. The limit set of $\widetilde{\mathcal{G}}_F$ must be dense, since otherwise any leaf of $\widetilde{\mathcal{G}}_{L_0}$ should share a single common ideal point as shown in Lemma~\ref{nondense=single}. It follows that for any given $\xi\in\partial_{\infty}F$, there always is a leaf $l\in \widetilde{\mathcal{G}}_F$ with either $l^+\neq \xi$ or $l^-\neq \xi$. Note that $F$ might be distinct from $L_0$. But for our convenience, we will write $L_0=F$, which would not affect our result.
	
	\begin{figure}[htb]	
		\centering
		\includegraphics{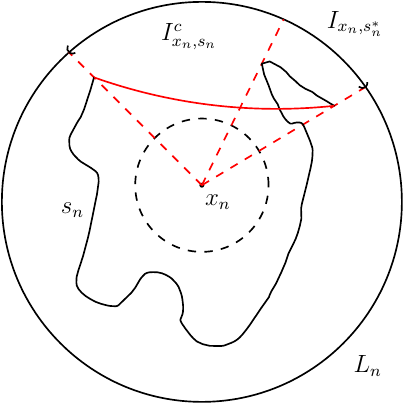}\\
		\caption{The shadows of the curve $s_n$ and its associated geodesic}
		
	\end{figure}
	
	Denote by $I_{x_n, s_n}\subset \partial_{\infty}L_n$ the shadow of the segment $s_n$ at the point $x_n$, and $I_{x_n, s_n}^c$ the complement ideal interval of $I_{x_n, s_n}$. Let $I_{x_n, s_n^*}\subset \partial_{\infty}L_n$ be the shadow of the geodesic segment $s_n^*$ at $x_n$. It is clear that $I_{x_n, s_n}^c$ is a sub-interval of $I_{x_n, s_n^*}$, but may be not proper. As $n$ tends to infinity, the visual measure $V_{x_n}(I_{x_n, s_n^*})$ shrinks to zero and the shadow $I_{x_n, s_n^*}$ converges to an ideal point, denoted by $\xi$, in the ideal boundary $\partial_{\infty}L_0$. Moreover, the wedge $W_{x_n}(I_{x_n, s_n^*})$ converges to a geodesic ray in $L_0$ from $x$ to $\xi$ as $n$ goes to infinity, denoted by $\alpha$. It turns out that the wedge $W_{x_n}(I_{x_n, s_n}^c)$ also converges to the geodesic ray $\alpha$.
	
	 Note that $l_0$ is a leaf of $\widetilde{\mathcal{G}}_{L_0}$ with $l_0^+\neq l_0^-$. Then $l_0$ has at least one ideal point distinct from $\xi$, say $l_0^+$. We can take a point $y\in l_0$ such that the ray of $l_0$ from $y$ to $l_0^+$, denoted by $r_0$, is disjoint with $\alpha$. In particular, the shadow $I_{x, r_0}\subset \partial_{\infty}L_0$ at the point $x$ does not contain $\xi$. Let $\theta:= V_x(I_{x, r_0})$ be the visual measure of the shadow $I_{x, r_0}$ at $x$. Changing the choice of $y$ if necessary, we can assume that $y$ is on the boundary of the wedge $W_x(I_{x, r_0})$ and $\theta\in(0, \pi)$. Let $d:= d_{L_0}(x, y)$ be the distance of $x$ and $y$ in the leaf $L_0$, which is finite. Denote by $T(x, a)$ a transversal to the foliation $\widetilde{\mathcal{F}}_1$, containing $x$ in the interior and having length $a>0$, that is contained in $\widetilde{\mathcal{F}}_2(x)$. By continuity, for any $\epsilon>0$, there is $\delta_d>0$ such that if a leaf $F\in \widetilde{\mathcal{F}}_1$ intersects a transversal $T(x, \delta_d)$, then it intersects $T(y, \epsilon)$. Since the leaves $L_n$ accumulate at the leaf $L_0$, there exists $N_{\delta_d}>0$ such that for any $n\geq N_{\delta_d}$, the leaf $L_n$ intersects a transversal $T(x, \delta_d)$. Therefore, the leaf $L_n$ intersects $T(y, \epsilon)$ in a point denoted by $y_n$.
	 
	 As indicated in Theorem \ref{Thurston}, for any $\rho>0$ and any leaf $F\in \widetilde{\mathcal{F}}_1$, there is a dense set of $\rho$ non-expanding directions at every point of $F$. For any $\rho>0$, there is $\epsilon\leq \rho$ such that if a leaf $F\in \widetilde{\mathcal{F}}_1$ intersects $T(y, \epsilon)$, then there is a dense set of $\rho$ non-expanding directions at $y$ in $L_0$ along which all points have distance less than $\rho$ to $F$. Pick a small $\omega>0$ so that any ideal interval of $\partial_{\infty}L_0$ containing $\xi$ in the interior of visual measure less than $5\omega$ at $x$ is disjoint with $I_{x, r_0}$. Let $\eta^1, \eta^2\in \partial_{\infty}L_0$ be two ideal points on the opposite sides of $l_0^+$ associated to two $\rho$ non-expanding directions at $y$. By the density of non-expanding directions and the fact that $x$ is at bounded distance from $y$, we can choose $\eta^1$ and $\eta^2$ to satisfy that the ideal interval between $\eta^1$ and $\eta^2$ containing $l_0^+$ has visual measure less than $\omega$ at the point $x$. The non-expanding directions from $y$ to $\eta^1$ and $\eta^2$ define two quasi-geodesics in $L_n$ starting at $y_n$ for each $n\geq N_{\delta_d}$. This uniquely determines two corresponding ideal points $\eta^i_n$, $i=1,2$, in the ideal boundary $\partial_{\infty}L_n$. The transverse continuity of ideal points shows that if a leaf of $\widetilde{\mathcal{G}}$ has an ideal point contained in the interior of an ideal interval that is bounded by two non-expanding ideal points, then this relation is persistent in nearby $\widetilde{\mathcal{F}}_1$-leaves joining the transversal associated to these non-expanding directions. It implies that the leaf of $\widetilde{\mathcal{G}}_{L_n}$ through $y_n$ has an ideal point $l(y_n)^+$ contained in the small interval bounded by $\eta_n^1$ and $\eta_n^2$. Denote by $r_n$ the ray of $l(y_n)\in \widetilde{\mathcal{G}}_{L_n}$ from $y_n$ to the ideal point $l(y_n)^+$.
	 
	 \begin{figure}[htb]	
	 	\centering
	 	\subcaptionbox{The choice of $y$ makes sure that the shadow of $r_0$ is separated from the limit ideal point $\xi$}
	 	{\includegraphics{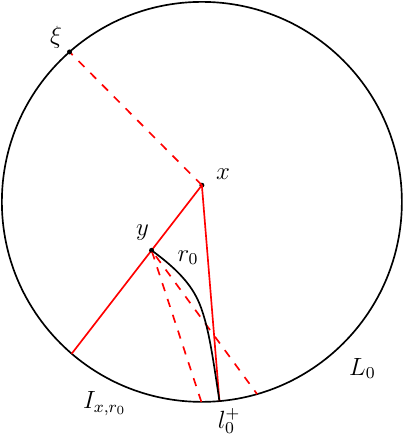}}
	 	\subcaptionbox{The shadow of $r_n$ is not contained in the shadow of $s_n$ if $r_n$ and $s_n$ are disjoint for $n$ large enough}
	 	{\includegraphics{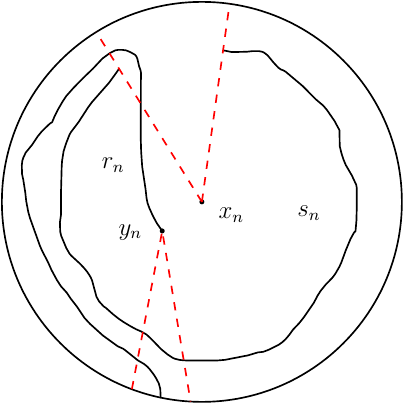}}
	 	\caption{Bounded distance from geodesics to curves under transverse continuity}
	 	
	 \end{figure}
	 
	 Let $I^i\subset \partial_{\infty}L_0$ be the ideal interval bounded by $\eta^i$ and the ideal point of the geodesic ray from $x$ through $y$ that intersects $I_{x, r_0}$, for $i=1, 2$. By the choice of $\eta^1$ and $\eta^2$, we can deduce that $|V_x(I^1)-V_x(I^2)|\leq \omega$. Denote by $I^i_n\subset \partial_{\infty}L_n$ the ideal interval bounded by $\eta^i_n$ and the ideal point of the geodesic ray starting at $x_n$ through $y_n$. For the given small $\omega>0$, there exists $N_{\omega}\in \mathbb{N}$ such that for any $n\geq N_{\omega}$, we have $|V_x(I^i)-V_{x_n}(I^i_n)|\leq \omega$ for $i=1, 2$. It turns out that for any $n\geq N_{\omega}$, we have $|V_{x_n}(I_{x_n, r_n})-V_x(I_{x, r_0})|\leq 3\omega$, since the ideal point $l(y_n)^+$ is bounded by $\eta^1_n$ and $\eta^2_n$. 
	 
	 As the wedge $W_{x_n}(I_{x_n, s_n}^c)$ converges to the geodesic ray $\alpha$ which has the ideal point $\xi$, there is $N_{\omega}'\in \mathbb{N}$ such that for any $n\geq N_{\omega}'$, the visual measure $V_{x_n}(I_{x_n, s_n}^c)$ is less than $\omega$. Also by the choice of $y$, for any $n\geq \max\{N_{\delta_d}, N_{\omega}, N_{\omega}'\}$, the ideal intervals $I_{x_n, r_n}$ and $I_{x_n, s_n}^c$ are disjoint. It concludes that the shadow $I_{x_n, r_n}$ is contained in the shadow $I_{x_n, s_n}$. However, for any $n\geq N_{\delta_d}$, we have 
	 \begin{equation}
	 	\begin{aligned}
	 		d(x_n, y_n) &\leq d(x_n, x)+ d(x, y)+ d(y, y_n)\\
	 		&\leq d(x_n, x)+ d_{L_0}(x, y)+ d(y, y_n)\\
	 		&\leq d+2\epsilon< 2d.
	 	\end{aligned}\nonumber
	 \end{equation}
	 By Lemma \ref{AB}, there is finite number $B>0$ such that $d_{L_n}(x_n, y_n)<B$ for any $n\geq N_{\delta_d}$. 
	 
	 Take $n>\max\{N_{\delta_d}, N_{\omega}, N_{\omega}', B\}$. Note that the curve $s_n$ is outside of the disk $D(x_n, n)$ in $L_n$, while $y_n$ is contained in $D(x_n, B)\subset D(x_n, n)$. Then, as $I_{x_n, r_n}$ is contained in $I_{x_n, s_n}$, the curve $r_n$ has to intersect $s_n$. This is impossible, which completes the proof.	
\end{proof}

\subsection{General leaf-dependent bounded distance}

The subsequent result addresses a broader framework where transverse continuity is not assumed. It presents a dichotomy where: The first case establishes a relationship analogous to Proposition~\ref{uniform_neighborhood_id} between curves and their associated geodesics, though permitting curve-dependent variation in constants; The second case demonstrates a foliation behavior exhibiting structural rigidity.

\begin{prop}\label{distance_geodesic to leaf_dichotomy}
	Assume that there is a leaf $l_0\in \widetilde{\mathcal{G}}_{L_0}$ for some leaf $L_0\in\widetilde{\mathcal{F}}_1$ such that $l_0^+\neq l_0^-$. Then either
	\begin{enumerate}
	    \item for any $L\in \widetilde{\mathcal{F}}_1$ and any $l\in \widetilde{\mathcal{G}}_L$, there is a constant $C_l>0$ such that the geodesic $l^*$ is contained in the $C_l$-neighborhood of $l$; or
        \item the limit set of $\widetilde{\mathcal{G}}_{L_0}$ is equal to $\{l_0^+, l_0^-\}$. 
	\end{enumerate}
The same statement holds for any ray and any segment in a leaf of $\widetilde{\mathcal{G}}_L$.
\end{prop}
\begin{proof}
	The first statement always holds if $l$ is chosen as a segment in a leaf of $\widetilde{\mathcal{G}}_L$. We only need to prove it for rays since the same argument also holds for leaves. 
	
	Suppose that the first possibility is false. Then there is $L\in \widetilde{\mathcal{F}}_1$ and a ray $l$ in a leaf of $\in \widetilde{\mathcal{G}}_L$ such that the geodesic ray $l^*$ is not contained in any bounded neighborhood of $l$ with respect to the induced distance $d_L$. It implies that there is a point $p_n\in L$ such that the disk $D(p_n, n)$ in $L$ is contained in a connected region bounded by $l$ and $l^*$ for each $n\in \mathbb{N}$. Denote by $l^{\infty}\in \partial_{\infty}L$ the ideal point of $l$. One can observe that $p_n$ converges to $l^{\infty}$ as $n$ goes to infinity.
	
	By the minimality of $\mathcal{F}_1$, there is a sequence of deck transformations $(\gamma_n)_{n\in \mathbb{N}}$ such that up to a subsequence, the point $\gamma_n(p_n)$ converges to a point $p_0\in L_0$. At the same time, the leaf $\gamma_n(L)$ converges to the leaf $L_0$ as $n$ goes to infinity. Let $L_n:= \gamma_n(L)$ and $I_{\gamma_n(p_n), \gamma_n(l^*)}\subset \partial_{\infty}L_n$ be the shadow of the geodesic ray $\gamma_n(l^*)$ at the point $\gamma_n(p_n)$. As indicated in the proof of Proposition \ref{uniform_neighborhood_id}, the shadow $I_{\gamma_n(p_n), \gamma_n(l^*)}$ converges to a single ideal point in $\partial_{\infty}L_0$, denoted by $\xi$. We are going to study the relation between $\xi$ and ideal points of $l_0$.
	
	Firstly, we assume that $\xi$ is not contained in $\{l_0^+, l_0^-\}$. Let $E\in \widetilde{\mathcal{F}}_2$ be the leaf containing $l_0$. As $L_0$ is accumulated by $L_n$, the leaf $l_0$ is also accumulated by a sequence of leaves in $L_n$, denoted by $l_n\in \widetilde{\mathcal{G}}_{L_n}$. Without loss of generality, up to a finite sequence, we can assume that $l_n$ is contained in the intersection of $E$ and $L_n$. Since the point $p_0$ and the leaf $l_0$ belong to the same leaf $L_0$, the distance $d_{L_0}(p_0, l_0)$ is bounded. For $n$ large enough, the leaf $l_n$ intersects the disk $D(\gamma_n(p_n), n)$ in $L_n$, see also the proof of Proposition \ref{uniform_neighborhood_id} for details. By the fact that $D(\gamma_n(p_n), n)$ is contained in a connected region, denoted by $R_n\subset L_n$, enclosed by $\gamma_n(l)$ and $\gamma_n(l^*)$, the leaf $l_n$ either escape the region $R_n$ from the geodesic ray $\gamma_n(l^*)$ or has a unique ideal point $\gamma_n(l^{\infty})\in \partial_{\infty}L_n$. 
	
	\begin{figure}[htb]	
		\centering
		\includegraphics{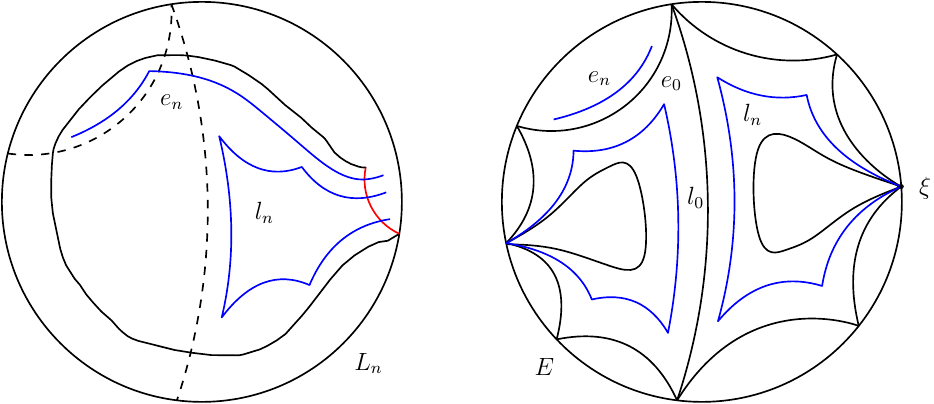}\\
		\caption{The leaf $l_0$ has two ideal points distinct from $\xi$}
		
	\end{figure}
	
	We claim that there are other leaves in $L_0\cap E$ distinct from $l_0$ accumulated by $l_n$. Indeed, if $l_n$ converges to a single leaf $l_0$, then the ideal points $l_n^{\pm}$ converge to $l_0^{\pm}$. Otherwise, by Proposition \ref{transverse-continuity-id}, there is at least one leaf non-separated from $l_0$ in $E$ accumualted by $l_n$. Using Proposition~\ref{intersection}, this leaf is also contained in $L_0$, which implies the claim. Now, the fact that the ideal points $l_n^{\pm}$ converge to $l_0^{\pm}$ implies the first statement of the proposition by the same argument as in the proof of Proposition \ref{uniform_neighborhood_id}. 
	
	Now, there is a collection of leaves including $l_0$ in $L_0$ that are accumulated by $l_n$ and non-separated from each other in $E$. It implies that all these leaves are contained in the half plane of $E$ bounded by $l_0$ that also contains $l_n$. Denote by $H_{l_0}^+\subset E$ this half plane and $H_{l_0}^-\subset E$ be the complementary half plane. We will show that there are also leaves in $H_{l_0}^-$ non-separated from $l_0$. Indeed, if not, then there are leaves $l_n'\in H_{l_0}^-$ converging to $l_0$ such that the ideal points $l_n'^{\pm}$ converge to $l_0^{\pm}$ by Proposition \ref{transverse-continuity-id}. We now have leaves $L_n'\in \widetilde{\mathcal{F}}_1$ containing $l_n'$ that converge to $L_0$ on the opposite side of $L_0$ from $L_n$. By minimality of $\mathcal{F}_1$, we can change the sequence of deck transformations $\gamma_n$ to $\gamma_n'$ so that $\gamma_n'(L)$ converges to $L_0$ on the same side as $L_n'$. Without loss of generality, we take $L_n':=\gamma_n'(L)$. Then the same argument as in the proof of Proposition \ref{uniform_neighborhood_id} shows that the first statement of the current proposition should hold, which provides a contradiction. Thus, the leaf $l_0$ is non-separated in $H_{l_0}^+$ and $H_{l_0}^-$. 
	
	\begin{figure}[htb]	
		\centering
		\includegraphics{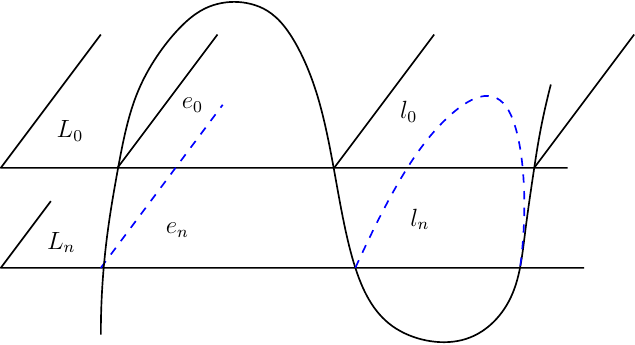}\\
		\caption{The leaf $e_n$ is contained in the half plane disjoint with $l_0$ bounded by $e_0$}
		
	\end{figure}
	
	Denote by $A$ the collection of leaves in $H_{l_0}^-$ non-separated from $l_0$ and $B$ the collection of leaves in $H_{l_0}^+$ non-separated from $l_0$. By Proposition~\ref{intersection}, we know that all leaves in $E$ non-separated from $l_0$ are also contained in $L_0$ since $\mathcal{F}_1$ is $\mathbb{R}$-covered. Moreover, the sets $A$ and $B$ belong to two distinct half planes in $L_0$ bounded by $l_0$, denoted by $K_{l_0}^-$ and $K_{l_0}^+$, respectively. It is not hard to see that $\xi$ is an ideal point in the boundary of the half plane $K_{l_0}^+$. Let $e_0\subset H_{l_0}^-$ be a leaf in $A$, which is also contained in $K_{l_0}^-$. Denote by $H_{e_0}^+$ the half plane of $E$ that is bounded by $e_0$ and contains $l_0$, and $H_{e_0}^-$ the complementary half plane in $E$. By transversality, for each $n$ large enough, the leaf $L_n$ intersects $H_{l_0}^-$ in a leaf, denoted by $e_n\subset L_n\cap E$, very close to $e_0\subset L_0\cap E$. We can always find a large $n$ such that $e_n$ intersects the disk $D(\gamma_n(p_n), n)$ in the leaf $L_n$. As discussed before, since $D(\gamma_n(p_n), n)$ is contained in the region $R_n$, the leaf $e_n$ either escapes from $R_n$ through the geodesic ray $\gamma_n(l^*)$ or accumulates at the ideal point $\gamma_n(l^{\infty})$. It turns out by the same argument for $l_n$ that $e_n$ converges to several leaves non-separated from $e_0$ in $E$. As the shadow $I_{\gamma_n(p_n), \gamma_n(l^*)}$ converges to $\xi$ and $\xi$ is in the boundary of $K_{l_0}^+$, there is at least one leaf in $K_{l_0}^+$ accumulated by $e_n$, which is a leaf in $B$. It implies that there is a leaf in the half plane $H_{l_0}^+$ non-separated from the leaf $e_0$ in the half plane $H_{l_0}^-$. This is absurd since two half planes are separated by $l_0$. 
	
	Therefore, we conclude that the ideal point $\xi$ is contained in $\{l_0^+, l_0^-\}$. We assume that $\xi$ coincides with $l_0^+$. 
	
	\begin{figure}[htb]	
		\centering
		\includegraphics{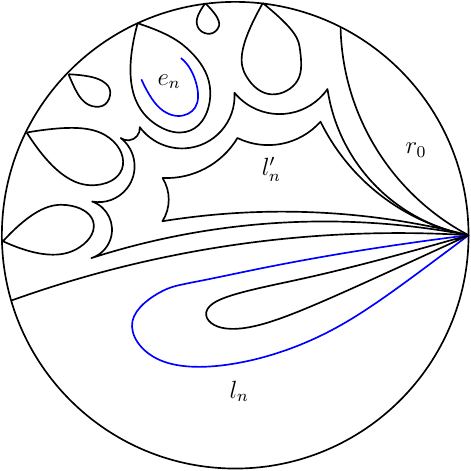}\\
		\caption{The leaf $L_n$ contains a leaf $e_n$ with a single ideal point distinct from $\xi$}
		
	\end{figure}
	
	For any $n$ large enough, there is a leaf $l_n\in\widetilde{\mathcal{G}}_{L_n}$ in the intersection of $L_n\cap E$ accumulating in the leaf $l_0$. We can choose $n$ large so that $l_n$ intersects the disk $D(\gamma_n(p_n), n)$. Note that the region $R_n$ contains the disk $D(\gamma_n(p_n), n)$ in $L_n$. Then the leaf $l_n$ either escapes from $R_n$ through $\gamma_n(l^*)$ or converges to $\gamma_n(l^{\infty})$. By the argument above, the leaf $l_n$ converges to several leaves non-separated from $l_0$ in $E$. 
	We first suppose that there is at least one leaf $r_0$ in $E$ non-separated from $l_0$ whose ideal points are not entirely contained in $\{l_0^+, l_0^-\}$. Recall that there are non-separated leaves with $l_0$ in the half plane $H_{l_0}^-$. Without loss of generality, we assume that $r_0$ is contained in the half plane $H_{l_0}^-$. The case where it is contained in $H_{l_0}^+$ can be discussed analogously. As shown above, every leaf non-separated from $l_0$ in $E$ that does not accumulate in the ideal point $\xi$ has a unique ideal point since any such leaf is also contained in $L_0$ by Proposition~\ref{intersection}. Let $l_n'\subset L_n' \cap E$ be a sequence of leaves converging to the leaves $l_0$ and $r_0$ in the half plane $H_{l_0}^-$. By the connectness of $l_n'$, there is at least one leaf, denoted by $e_0\subset H_{l_0}^-$, non-separated from $l_0$ and $r_0$ whose ideal point $e_0^+=e_0^-$ is distinct from $\xi$. For each $n$ large enough, the leaf $L_n$ intersects the half plane $H_{e_0}^-$ in a leaf $e_n\subset L_n\cap E$ very close to $e_0\subset L_0\cap E$. Note that the leaf $e_n$ intersects the disk $D(\gamma_n(p_n), n)$ for $n$ large enough. As the region $R_n$ contains the disk $D(\gamma_n(p_n), n)$ in $L_n$, the leaf $e_n$ either escapes from $R_n$ through $\gamma_n(l^*)$ or converges to $\gamma_n(l^{\infty})$. This is absurd since $e_n$ is contained in the half plane $H_{e_0}^-$ with a unique ideal point $e_0^+=e_0^-$ that is distinct from $\xi$. 
	
	Thus, we conclude that the ideal points of every leaf non-separated from $l_0$ in $E$ are contained in $\{l_0^+, l_0^-\}$. Now we have two possibilities for the half plane $H_{l_0}^+$: either (1) every leaf in $H_{l_0}^+$ has a single ideal point in $\{l_0^+, l_0^-\}$; or (2) there is a unique leaf $s_1\subset H_{l_0}^+$ accumulated by $l_n$ with two distinct ideal points $\{s_1^+, s_1^-\}=\{l_0^+, l_0^-\}$ such that every leaf accumulated by $l_n$ in $E$ is bounded by $l_0$ and $s_1$. As we mentioned above, there are leaves in $H_{l_0}^-$ non-separated from $l_0$. We also have two possibilities for $H_{l_0}^-$: either (1) every leaf in $H_{l_0}^-$ has a single ideal point in $\{l_0^+, l_0^-\}$; or (2) there is a unique leaf $s_{-1}\subset H_{l_0}^-$ non-separated from $l_0$ with two distinct ideal points $\{s_{-1}^+, s_{-1}^-\}=\{l_0^+, l_0^-\}$ such that every leaf non-separated from $l_0$ in $H_{l_0}^-$ is bounded by $l_0$ and $s_{-1}$.
	
	\begin{figure}[htb]	
		\centering
		\includegraphics{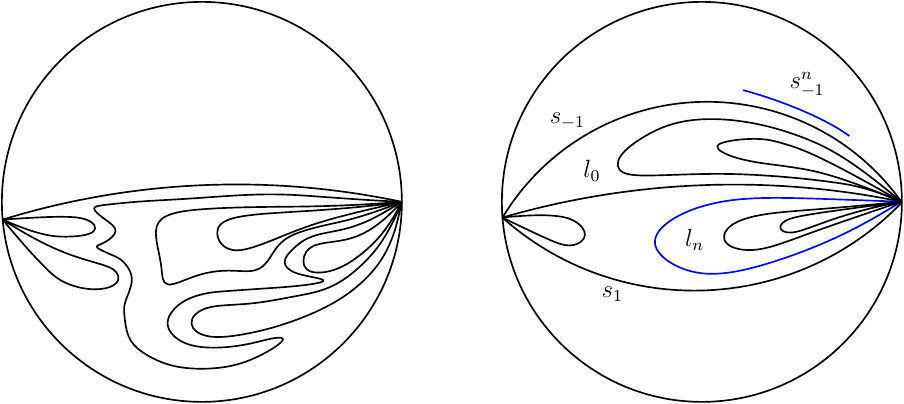}\\
		\caption{Two possibilities in each half plane: either every leaf has a single ideal point contained in $\{l_0^+, l_0^-\}$, or there is a unique leaf bounding a new half plane with the same ideal points}
		
	\end{figure}

	Assume the existence of $s_{-1}$. By transversality, for each $n$ large enough, the leaf $L_n$ intersects $H_{l_0}^-$ in a leaf $s_{-1}^n$ accumulating at $s_{-1}$ as $n$ goes to infinity. As discussed above, there are leaves in $H_{s_{-1}}^-\subset H_{l_0}^-$ non-separated from $s_{-1}$. Analogously, we also have two possibilities on the leaves in $H_{s_{-1}}^-$ non-separated from $s_{-1}$. 
	
	Assume the existence of $s_1$. We can discuss in a similar way on the opposite side of $L_0$ from $L_n$ by considering $L_n':=\gamma_n'(L)$ for another sequence of deck transformations $(\gamma_n')_{n\in \mathbb{N}}$. Then for $n$ large enough, the leaf $L_n'$ intersects $H_{l_0}^+$ in a leaf $s_1^n$ converging to $s_1$ as $n$ goes to infinity. The same argument applies to obtain two possibilities on the leaves in $H_{s_1}^+$ non-separated from $s_1$. 
	
	Inductively, we obtain a sequence of leaves 
	$$\cdots s_{-2}, s_{-1}, s_0=l_0, s_1, s_2, \cdots$$
	in $E\cap L_0$ so that each consecutive pair is non-separated in $E$, which may be finite in one side or in both sides. Here, for each $i\in \mathbb{Z}$, we have $\{s_i^+, s_i^-\}=\{l_0^+, l_0^-\}$. For each point $q\in L_0$, there is an integer $N_q\in \mathbb{Z}$ such that the leaf $l(q)\in \widetilde{\mathcal{G}}_{L_0}$ is bounded by $s_{N_q}$ and $s_{N_q+1}$. It turns out that $l(q)^{\pm}\in \{l_0^+, l_0^-\}$. Thus, the limit set of $\widetilde{\mathcal{G}}_{L_0}$ in $\partial_{\infty}L_0$ is equal to $\{l_0^+, l_0^-\}$.
	
	Hence, we finish the proof.
\end{proof}

Following Lemma~\ref{nondense=single} and the dichotomy in Proposition~\ref{distance_geodesic to leaf_dichotomy}, we obtain a consequence.
\begin{cor}\label{distance-corollary}
    Under the conditions of Proposition~\ref{distance_geodesic to leaf_dichotomy}, for any $L\in \widetilde{\mathcal{F}}_1$ and any $l\in \widetilde{\mathcal{G}}_L$, there is a constant $C_l>0$ such that the geodesic $l^*$ is contained in the $C_l$-neighborhood of $l$. The same also holds for any ray and any segment in a leaf of $\widetilde{\mathcal{G}}_L$.
\end{cor}

\section{Quasi-geodesic subfoliation}\label{section-Hausdorff}

This section investigates leafwise geometric structures of 2-dimensional foliations saturated by a one-dimensional subfoliation. We consider transverse minimal $\mathbb{R}$-covered foliations $\mathcal{F}_1$ and $\mathcal{F}_2$ comprising non-compact Gromov hyperbolic leaves in a closed 3-manifold, where $\mathcal{G} = \mathcal{F}_1 \cap \mathcal{F}_2$ denotes their intersecting one-dimensional subfoliation. By analyzing asymptotic properties of $\widetilde{\mathcal{G}}$ within leaves of $\widetilde{\mathcal{F}}_1$, we derive intrinsic geometric characterizations for $\widetilde{\mathcal{F}}_1$-leaves.

\subsection{Hausdorff leaf space}

Given a lifted foliation $\widetilde{\mathcal{F}}$ and a leaf $F\in \widetilde{\mathcal{F}}$, a curve $l\subset F$ is a \emph{quasi-geodesic} if there exist constants $\lambda>1$ and $c>0$ such that 
$$
\lambda^{-1}d_F(x, y)-c\leq d_l(x, y)\leq \lambda d_F(x,y)+c
$$
for any $x, y\in l$, where $d_F$ denotes the induced leaf metric in $F$ and $d_l$ denotes the path metric in $l$. A family of curves constitutes \emph{uniform quasi-geodesics} if the constants $\lambda, c$ can be chosen independent of curves. We mention that the proof of \cite[Lemma 10.20]{Calegari07book} shows that a one-dimensional foliation $\widetilde{\mathcal{G}}$ by quasi-geodesic leaves indeed comprises uniform quasi-geodesics if $\mathcal{G}$ subfoliates a two-dimensional foliation $\mathcal{F}$ by Gromov hyperbolic leaves in a closed manifold $M$. Even though, we will not use this fact in our subsequent arguments.

We will present equivalent statements for having Hausdorff leaf space of $\widetilde{\mathcal{G}}$ in leaves of $\widetilde{\mathcal{F}}_1$ under some technical condition. It is worth noting that this technical condition appeared in the following result is necessary for those equivalence since counter examples can be readily constructed. We mention here that the recent paper \cite{FP23intersection} shows a more general result. Our proof in this section is independent of \cite{FP23intersection}. 

\begin{prop}\label{equivalent}
	Assume that there is a leaf $l_0\in \widetilde{\mathcal{G}}_F$ for some $F\in \widetilde{\mathcal{F}}_1$ such that $l^+_0\neq l^-_0$. Then the following statements are equivalent: for any $L\in \widetilde{\mathcal{F}}_1$,
	\begin{enumerate}
	    \item every leaf $l\in\widetilde{\mathcal{G}}_L$ has two distinct ideal points $l^+$ and $l^-$;
        \item all leaves of $\widetilde{\mathcal{G}}_L$ are uniform quasi-geodesics;
        \item the leaf space of $\widetilde{\mathcal{G}}_L$ is Hausdorff.
	\end{enumerate}
\end{prop}
\begin{proof}
	It follows from \cite[Lemma 5.9]{FU1} that (1) implies (2). The implication from (2) to (3) can be obtained by the same argument as in \cite[Lemma 5.10]{FU1}. It suffices to show item (1) by assuming (3). Suppose there is a leaf $l\in \widetilde{\mathcal{G}}_L$ for some $L\in \widetilde{\mathcal{F}}_1$ such that $l^+=l^-$. By Proposition~\ref{intersection} and Proposition \ref{transverse-continuity-id}, the foliation $\widetilde{\mathcal{F}}_1$ has transverse continuity. Picking any point $z\in l$, we denote by $l_1$, $l_2$ two rays of $l$ from $z$ to the ideal point $l^+=l^-$. Applying Proposition \ref{uniform_neighborhood_id} for rays, there is a uniform constant $C>0$ such that the geodesic ray from $z$ to $l^+$ is entirely contained in the $C$-neighborhood of $l_1$ and $l_2$. We could find two sequences of points $(x_n)_{n\in \mathbb{N}}\subset l_1$, $(y_n)_{n\in \mathbb{N}}\subset l_2$ such that both $x_n$ and $y_n$ converge to the ideal point $l^+=l^-$ and the distance $d_L(x_n, y_n)$ is uniformly bounded by $2C$. By the compactness of $M$, up to taking subsequence of points, there is a sequence of deck transformations $(g_n)_{n\in \mathbb{N}}$ such that $g_n(x_n)$, $g_n(y_n)$ converge simultaneously as $n$ goes to infinity. Let $x_0$ and $y_0$ be two accumulation points of $g_n(x_n)$ and $g_n(y_n)$, respectively, and $L_0\in \widetilde{\mathcal{F}}_1$ be the leaf containing both $x_0$ and $y_0$. As any deck transformation preserves the foliation $\widetilde{\mathcal{G}}$ and the arc of $l$ connecting $x_n$ and $y_n$ has unbounded length as $n$ goes to infinity, the leaves of $\widetilde{\mathcal{G}}_{L_0}$ through $x_0$ and $y_0$ are distinct, denoted by $l(x_0), l(y_0)\in \widetilde{\mathcal{G}}_{L_0}$, respectively. 
	
	By transversality, there is $\epsilon>0$ such that for any point $p\in \widetilde{M}$ with $d(p, D_{L_0}(x_0, 2C))\leq \epsilon$, the leaf of $\widetilde{\mathcal{F}}_2$ through $p$ intersects $L_0$ in a $\widetilde{\mathcal{G}}$-curve. For this $\epsilon$, there exists an integer $N>0$ such that for any $n\geq N$, we have $d(g_n(x_n), x_0)\leq \epsilon$ and $d(g_n(y_n), y_0)\leq \epsilon$. Then there is a curve $l(x_n')\in \widetilde{\mathcal{G}}_{L_0}$ contained in the leaf of $\widetilde{\mathcal{F}}_2$ through $g_n(x_n)$, where $x_n'\in L_0$ is the point satisfying that $d(g_n(x_n), x_n')= \inf\{d(g_n(x_n), p); p\in l(x_n')\}$. The same leaf of $\widetilde{\mathcal{F}}_2$ also contains $g_n(y_n)$ and intersects $L_0$ in a curve $l(y_n')\in \widetilde{\mathcal{G}}_{L_0}$, where $y_n'\in L_0$ is a point satisfying that $d(g_n(y_n), y_n')= \inf\{d(g_n(y_n), p); p\in l(y_n')\}$. Since $g_n(x_n)$ converges to the point $x_0$, the point $x_n'$ converges to $x_0$ in the leaf $L_0$ as $n$ goes to infinity. Similarly, the point $y_n'$ converges to $y_0$ in the leaf $L_0$ as $n$ goes to infinity.
	
	Suppose that there are infinitely many integers $n\geq N$ such that $l(x_n')$ and $l(y_n')$ coincide. By the fact that $x_0$ and $y_0$ are accumulation points of $x_n'$ and $y_n'$, respectively, we deduce that $l(x_0)$ and $l(y_0)$ are two distinct leaves of $\widetilde{\mathcal{G}}_{L_0}$ accumulated by $l(x_n')=l(y_n')$ as $n$ goes to infinity. It implies that the leaf space of $\widetilde{\mathcal{G}}_{L_0}$ is not Hausdorff, which contradicts the assumption (3).
	
	Thus, there exists an integer $N_1\geq N$ so that for any $n\geq N_1$, the leaf $l(x_n')$ is distinct from $l(y_n')$. Note that the points $g_n(x_n)$ and $g_n(y_n)$ are contained in the same leaf of $\widetilde{\mathcal{G}}_{g_n(L)}$. Then there are two distinct leaves $l(x_n')$ and $l(y_n')$ in the intersection of a leaf of $\widetilde{\mathcal{F}}_2$ through $g_n(x_n)$ and $L_0$. As indicated by Theorem~\ref{taut} and Proposition~\ref{intersection}, the leaf space of $\widetilde{\mathcal{G}}_{L_0}$ cannot be Hausdorff. This is a contradiction.
	
	Hence, we finish the proof of implication from (3) to (1), and complete the proof.	
\end{proof}

\begin{rmk}
    Under the assumption and condition (3) of Proposition~\ref{equivalent}, we can actually obtain a more precise description on the limit leaf: for accumulation points $x_0, y_0\in L_0$, their leaves $l(x_0), l(y_0)\in \widetilde{\mathcal{G}}_{L_0}$ are indeed a non-separated pair in $L_0$ and accumulated simultaneously by $l(x_n')=l(y_n')\in \widetilde{\mathcal{G}}_{L_0}$ for $n$ large enough. Moreover, the leaves $l(x_n')$ have a single common ideal point $l(x_n')^+=l(x_n')^-\in \partial_{\infty}L_0$ for all $n$. As this fact is out of our use, we do not include its proof here.
\end{rmk}

\subsection{Leaf-wise quasi-geodesic subfoliation}

A leaf $ F \in \widetilde{\mathcal{F}} $ is called a \textit{weak quasi-geodesic fan} for a uniformly quasi-geodesic subfoliation \( \widetilde{\mathcal{G}}_F \) if there exists a unique ideal point $\eta_F \in \partial_\infty F$ serving as the common limit for all $\widetilde{\mathcal{G}}_F$-leaves. The uniqueness of $\eta_F$ is mandatory, while the existence of leaves with distinct ideal endpoints is necessary to maintain uniform quasi-geodesic property. We designate $\eta_F$ as the \emph{funnel point} of $F$.

By \cite[Lemma 5.3]{FU1}, uniformly quasi-geodesic foliations $\widetilde{\mathcal{G}}_F$ exhibit continuous variation of ideal points $l(x)^\pm \in \partial_\infty F$ as $x$ ranges over $F$. This ensures that every ideal point $\zeta \in \partial_\infty F \setminus \{\eta_F\}$ connects to $\eta_F$ via at least one $\widetilde{\mathcal{G}}_F$-leaf, enhancing the fan structure. A \emph{quasi-geodesic fan} strengthens this by requiring unique $\widetilde{\mathcal{G}}_F$-leaves between each $\zeta$ and $\eta_F$, unlike weak fans permitting multiple leaves with identical endpoints. See \cite{BFP23collapsed,Calegari06} for expanded discussions on (quasi-)geodesic foliation theories.

\begin{prop}\cite[Proposition 5.14]{FU1}\cite[Proposition 6.9]{BFP23collapsed}\label{fan1-weak}
    If some $L \in \widetilde{\mathcal{F}}_1$ has all $\widetilde{\mathcal{G}}_L$-leaves as uniform quasi-geodesics, then every $F \in \widetilde{\mathcal{L}}$ becomes a weak quasi-geodesic fan for $\widetilde{\mathcal{G}}_F$, where $\mathcal{L}$ denotes the minimal $\mathcal{F}_1$-sublamination containing $\pi(L)$.
\end{prop}

\begin{prop}\label{fan-weak}
    Let $\mathcal{L}$ be a minimal sublamination of $\mathcal{F}_1$. Assume that there is a leaf $l_0\in \widetilde{\mathcal{G}}_L$ for some $L\in \widetilde{\mathcal{L}}$ such that $l^+_0\neq l^-_0$, and $\widetilde{\mathcal{G}}_{L'}$ has Hausdorff leaf space for some $L'\in \widetilde{\mathcal{L}}$. Then every leaf $F\in \widetilde{\mathcal{L}}$ is a weak quasi-geodesic fan for $\widetilde{\mathcal{G}}_F$. In particular, $\mathcal{G}$ contains a closed leaf in $\mathcal{L}$.
\end{prop}
\begin{proof}
    Since the leaf space of $\widetilde{\mathcal{G}}_{L'}$ is Hausdorff, every leaf $F\in \widetilde{\mathcal{L}}$ admits Hausdorff $\widetilde{\mathcal{G}}_F$-leaf space by Proposition~\ref{intersection}. Proposition~\ref{equivalent} implies that all leaves of $\widetilde{\mathcal{G}}_F$ are uniform quasi-geodesics. Then, each leaf $F\in \widetilde{\mathcal{L}}$ is a weak quasi-geodesic fan by Proposition~\ref{fan1-weak}. It turns out from \cite[Corollary 5.15]{FU1} that every leaf of $\widetilde{\mathcal{L}}$ stabilized by some non-trivial deck transformation projects to a cylinder leaf in $M$, which exits due to Theorem~\ref{rosenberg}. Moreover, such a deck transformation fixes a leaf of $\widetilde{\mathcal{G}}$, which projects to a closed leaf in $M$.
\end{proof}

We point out that we did not require the uniformness of foliations in our preceding results. Here, we impose this restriction for the rest of this section, even though they are not necessary for our proof of accessibility.

\begin{lem}\label{fan-weak-strong}
    Let $\mathcal{F}_1, \mathcal{F}_2$ be two transverse foliations by non-compact Gromov hyperbolic leaves, and $\mathcal{G}$ be the intersected foliation. Assume that $\mathcal{F}_2$ is minimal and uniform. If $\widetilde{\mathcal{G}}$ forms a weak quasi-geodesic fan in each leaf of $\widetilde{\mathcal{F}}_1$ and $\widetilde{\mathcal{F}}_2$, then all leaves of $\widetilde{\mathcal{F}}_1$ are quasi-geodesic fans for $\widetilde{\mathcal{G}}$.
\end{lem}
\begin{proof}
    Suppose that there is a leaf $F\in \widetilde{\mathcal{F}}_1$ that is weak quasi-geodesic fan but not a quasi-geodesic fan. Then, there exist two leaves $l_1,l_2\in \widetilde{\mathcal{G}}_F$ with two common ideal points $l_1^+=l_2^+, l_1^-=l_2^-\in \partial_{\infty}F$. As $\widetilde{\mathcal{G}}_F$ is a uniformly quasi-geodesic foliation, there exists a constant $K>0$ such that both $l_1$ and $l_2$ are contained in the $K$-neighborhood of a geodesic connecting ideal points $l_1^+$ and $l_1^-$. It implies that $l_1$ and $l_2$ are at $2K$-distance apart in the leaf $F$ and, in particular, in $\widetilde{M}$. 

    Denote by $\mathcal{L}_2$ the leaf space of $\widetilde{\mathcal{F}}_2$, which is homeomorphic to $\mathbb{R}$. By the uniformness of $\mathcal{F}_2$, there is a quasi-isometric map producing a universal circle by identifying all ideal circles of $\widetilde{\mathcal{F}}_2$. Let $\mathcal{U}_2$ be the universal circle defined by the identification map $\tau_L:\mathcal{U}_2\rightarrow \partial_{\infty}L$ for any $L\in \widetilde{\mathcal{F}}_2$, where $\tau_{L'}\circ(\tau_L)^{-1}$ maps a point $x\in L$ to a point $y\in L'$ such that $d((x, y)$ is less than the Hausdorff distance $d_H(L, L')$. Define a map $\sigma:\mathcal{L}_2\rightarrow\mathcal{U}_2$ as $\sigma(L)=(\tau_L)^{-1}\circ\eta_L$, where $\eta_L\in \partial_{\infty}L$ denotes the funnel point of $L$.
    
    Let $L_1, L_2\in\widetilde{\mathcal{F}}_2$ be the leaves through $l_1, l_2$, respectively. These two leaves are distinct since $\widetilde{\mathcal{G}}_F$ has Hausdorff leaf space and $\mathcal{F}_2$ is Reebless (see Proposition~\ref{intersection}). The leaves $L_1$ and $L_2$ bound a closed leaf interval $I\subset \mathcal{L}_2$ in which every leaf intersects $F$ in a $\widetilde{\mathcal{G}}$-leaf bounded by $l_1$ and $l_2$. Note that the leaves in $I$ intersect $F$ in $\widetilde{\mathcal{G}}$-leaves of Hausdorff distance uniformly bounded by $2K\leq d_H(L_1,L_2)$. Then, for any leaves $L, L'\in I$, the map $\tau_{L'}\circ(\tau_L)^{-1}$ sends $\eta_L$ to $\eta_{L'}$. It turns out that $\sigma$ collapses the interval $I$ to a single point in $\mathcal{U}_2$.

    Analogously, for any deck transformation $\gamma$, the leaf interval $\gamma(I)\subset \mathcal{L}_2$ is also mapped by $\sigma$ to a single point of $\mathcal{U}_2$. Since $\mathcal{F}_2$ is a minimal foliation, the leaf space $\mathcal{L}_2$ is covered by the $\pi_1(M)$-orbit of $I$. Thus, the map $\sigma$ collapses all funnel points of $\widetilde{\mathcal{F}}_2$ to a single point of $\mathcal{U}_2$. This gives a global fixed point under $\pi_1(M)$ in the universal circle $\mathcal{U}_2$, which contradicts with \cite[Theorem 1.2]{FP20minimal}. We finish the proof.
\end{proof}

\begin{prop}\label{fan-strong}
    Let $\mathcal{F}_1, \mathcal{F}_2$ be two transverse minimal uniform foliations by non-compact Gromov hyperbolic leaves. Assume that each leaf $F\in\widetilde{\mathcal{F}}_i$ admits dense limit set in $\partial_{\infty}F$ for $\widetilde{\mathcal{G}}_F$, $i=1, 2$. If there is a leaf $L$ of $\widetilde{\mathcal{F}}_1$ or $\widetilde{\mathcal{F}}_2$ such that $\widetilde{\mathcal{G}}_{L}$ has Hausdorff leaf space, then all leaves of $\widetilde{\mathcal{F}}_1$ and $\widetilde{\mathcal{F}}_2$ are quasi-geodesic fans for $\widetilde{\mathcal{G}}$.
\end{prop}
\begin{proof}
    Assume that $\widetilde{\mathcal{G}}_L$ has Hausdorff leaf space for a leaf $L\in \widetilde{\mathcal{F}}_2$. Proposition~\ref{intersection} implies the Hausdorff leaf space of $\widetilde{\mathcal{G}}$ in all leaves of $\widetilde{\mathcal{F}}_1$ and $\widetilde{\mathcal{F}}_2$. As each leaf $F\in \widetilde{\mathcal{F}}_1$ admits dense limit set in $\partial_{\infty}F$ for $\widetilde{\mathcal{G}}_F$, there must be a leaf $l\in \widetilde{\mathcal{G}}_F$ with distinct ideal points $l^+\neq l^-\in \partial_{\infty}F$. Otherwise, if any leaf of $\widetilde{\mathcal{G}}_F$ possesses a single ideal point in $\partial_{\infty}F$, then the denseness of limit set implies that the leaf space of $\widetilde{\mathcal{G}}_F$ cannot be Hausdorff. By Proposition~\ref{equivalent}, the foliation $\widetilde{\mathcal{G}}_F$ constitutes uniform quasi-geodesics for each $F\in \widetilde{\mathcal{F}}_2$. Then, Proposition~\ref{fan-weak} implies that each leaf $F\in \widetilde{\mathcal{F}}_1$ is a weak quasi-geodesic fan. The same argument applies to show that every leaf of $\widetilde{\mathcal{F}}_2$ is also a weak quasi-geodesic fan. Lemma~\ref{fan-weak-strong} applies for both foliations to obtain both quasi-geodesic fans, which finishes the proof.
\end{proof}

\subsection{Topological Anosov flow}

Note that for topological Anosov flows, the lifted weak-stable and weak-unstable leaves form quasi-geodesic fans relative to the lifted orbit foliation. Our construction produces a topological Anosov flow via the foliation $\widetilde{\mathcal{F}}_1$ endowed with leafwise quasi-geodesic fans. As these results remain ancillary to subsequent developments, readers may safely omit this subsection without loss of continuity.

Now, let us present our main result of this subsection.

\begin{prop}\label{topAnosovflow}
    Let $\mathcal{F}_1, \mathcal{F}_2$ be two transverse foliations such that all leaves of $\widetilde{\mathcal{F}}_1$ and $\widetilde{\mathcal{F}}_2$ are quasi-geodesic fans for $\widetilde{\mathcal{G}}$. Then, there exists a topological Anosov flow $\phi_t: M\rightarrow M$ such that 
    \begin{itemize}
        \item the $\phi_t$-orbits coincide with the leaves of $\mathcal{G}$;
        \item the leaves of $\mathcal{F}_1$ (resp. $\mathcal{F}_2$) constitute the weak-stable (resp. weak-unstable) foliation of $\phi_t$.
    \end{itemize}
\end{prop}
\begin{proof}
    Since the foliation $\mathcal{G}$ has been assumed to be orientable, it determines a non-singular vector field $X_{\mathcal{G}}$ by fixing an orientation. This vector field generates a non-singular continuous flow $\phi_t$. By definition, it is clear that $\phi_t$ preserves two transverse foliations $\mathcal{F}_1$ and $\mathcal{F}_2$. We proceed to show that $\phi_t$ is an expansive flow.

    Denote by $\widetilde{\phi}_t:\widetilde{M}\rightarrow \widetilde{M}$ be the lifted flow of $\phi_t$. The expansiveness of $\widetilde{\phi}_t$ will imply that of $\phi_t$. Since $\widetilde{\phi}_t$ has not self accumulation point, it is sufficient to find an expansive constant separating $\widetilde{\phi}_t$-oribits in the sense of Hausdorff distance. There exists a constant $\epsilon>0$ such that each set of diameter less than $\epsilon$ is contained in a foliation chart for both $\widetilde{\mathcal{F}}_1$ and $\widetilde{\mathcal{F}}_2$. By transversality, there is a constant $\delta>0$ such that for any points $x, y\in \widetilde{M}$ of distance $d(x, y)\leq \delta$, the leaves $\widetilde{\mathcal{F}}_1(x)$ and $\widetilde{\mathcal{F}}_2(y)$ intersect in a unique connected curve in the $\epsilon$-neighborhood of $x$ and $y$.

    For each leaf $F\in \widetilde{\mathcal{F}}_1$, any pair of leaves of $\widetilde{\mathcal{G}}_F$ cannot be at bounded Hausdorff distance apart. Indeed, if there were two leaves of $\widetilde{\mathcal{G}}_F$ of bounded Hausdorff distance, then they must share two ideal points in $\partial_{\infty}F$. This is impossible since the leaf $F$ is a quasi-geodesic fan for $\widetilde{\mathcal{G}}_F$. Analogously, since all leaves of $\widetilde{\mathcal{F}}_2$ are quasi-geodesic fans for $\widetilde{\mathcal{G}}$, the Hausdorff distance between any pair of $\widetilde{\mathcal{G}}$-leaves in the same leaf of $\widetilde{\mathcal{F}}_2$ is unbounded. Suppose that there are two leaves $l, l'\in \widetilde{\mathcal{G}}$ of Hausdorff distance $d_H(l,l')\leq\delta$. By the choice of $\delta$, the intersection of the leaf of $\widetilde{\mathcal{F}}_1$ through $l$ and the leaf of $\widetilde{\mathcal{F}}_2$ through $l'$ determines a leaf $r\in \widetilde{\mathcal{G}}$ satisfying $d_H(l, r)\leq \epsilon$ and $d_H(l', r)\leq\epsilon$. Notice that $l$ and $r$ are two leaves of $\widetilde{\mathcal{G}}$ in the same leaf of $\widetilde{\mathcal{F}}_1$ of bounded Hausdorff distance. This is impossible by our argument above.

    Therefore, any pair of $\widetilde{\mathcal{G}}$-leaves has Hausdorff distance no less than $\delta$. It follows that the flow $\widetilde{\phi}_t$ is expansive with an expansivity constant $\delta$. Thus, $\phi_t$ is an expansive flow preserving two transverse foliations $\mathcal{F}_1$ and $\mathcal{F}_2$. Theorem~\ref{flow-defn} implies that $\phi_t$ is a topological Anosov flow. By \cite[Proposition 5.5]{BFP23collapsed}, up to reversing the flow, we have that $\mathcal{F}_1$ is the weak-stable foliation of $\phi_t$ and $\mathcal{F}_2$ is the weak-unstale foliation.
\end{proof}

The leafwise quasi-geodesic fan structure proves more than sufficient the preceding argument. The result remains valid under the weaker constraint that distinct $\widetilde{\mathcal{G}}$-leaves within any common leaf of $\widetilde{\mathcal{F}}_1$ or $\widetilde{\mathcal{F}}_2$ never share both ideal boundary points.

\begin{prop}\label{topAnosovflow-mixing}
    Under the conditions of Proposition~\ref{topAnosovflow}, if either $\mathcal{F}_1$ or $\mathcal{F}_2$ is minimal, then the flow $\phi_t$ is topological mixing (in particular, it is transitive).
\end{prop}
\begin{proof}
    For an arbitrarily small $\theta>0$, it suffices to show that for any two $\theta$-balls $U$ and $V$ in $M$, there exists $T>0$ such that $\phi_t(U)\cap V\neq \emptyset$ or $\phi_{-t}(U)\cap V\neq \emptyset$ for any $t\geq T$. Assume that $\mathcal{F}_1$ is minimal. Then, there exists a uniform constant $R_{\theta}>0$ such that each disk of radius $R_{\theta}$ in any leaf of $\mathcal{F}_1$ intersects every $\theta$-ball in $M$. Let $U_0\subset \widetilde{M}$ be a lift of $U$ and $F\in \widetilde{\mathcal{F}}_1$ be a leaf intersecting $U_0$ in a disk $D\subset F$. Pick two points $p, q\in D$ in distinct leaves of $\widetilde{\mathcal{G}}_F$. As $F$ is a quasi-geodesic fan and a weak-stable leaf of $\widetilde{\phi}_t$, the $\widetilde{\phi}_t$-orbits accumulate at the funnel point as $t\rightarrow+\infty$ and accumulate at distinct ideal points of $\partial_{\infty}F$ as $t\rightarrow -\infty$. It turns out that the distance $d_F(\widetilde{\phi}_t(p), \widetilde{\phi}_t(q))$ is unbounded as $t\rightarrow -\infty$. The region $\bigcup_{t>0}\widetilde{\phi}_{-t}(D)$ contains a disk of arbitrarily large radius. In particular, there exists $T>0$ such that $\bigcup_{t\geq T}\widetilde{\phi}_{-t}(D)$ contains a disk of radius $R_{\theta}$. Projecting to $M$, the half orbit $\phi_{-t}(U)$ always intersects $V$ for $t\geq T$. Thus, we finish the proof.
\end{proof}

\section{Degenerate limit set}\label{section-degenerate}
In this section, let $\mathcal{F}_1$ and $\mathcal{F}_2$ be transverse minimal $\mathbb{R}$-covered foliations with non-compact Gromov hyperbolic leaves in a closed 3-manifold $M$, and let $\mathcal{G} = \mathcal{F}_1 \cap \mathcal{F}_2$ denote their one-dimensional intersection. We assume orientability for $\mathcal{F}_1$, $\mathcal{F}_2$, $\mathcal{G}$, and $M$ for simplicity. Our focus centers on the case where every leaf $F \in \widetilde{\mathcal{F}}_1$ exhibits a degenerate limit set for $\widetilde{\mathcal{G}}$ (see Section~\ref{subsection-coherent}). By Lemma~\ref{nondense=single}, this degeneracy propagates universally across $\widetilde{\mathcal{F}}_1$-leaves if one leaf has a non-dense limit set.

The main result of this section is the following:
\begin{thm}\label{degenerate-limit-set}
    Let $\mathcal{F}_1, \mathcal{F}_2$ be two transverse uniform minimal foliations by non-compact Gromov hyperbolic leaves in a closed 3-manifold $M$ and $\mathcal{G}$ be their one-dimensional intersection. Assume that every leaf of $\widetilde{\mathcal{F}}_1$ has degenerate limit set for $\widetilde{\mathcal{G}}$. Then, $\mathcal{F}_2$ has trivial holonomy, and all leaves of $\mathcal{F}_2$ share $\pi_1(M)$ as their fundamental group. Moreover, for any transverse one-dimensional foliation $\mathcal{T}$ to $\mathcal{F}_2$, the lifted foliations $\widetilde{\mathcal{F}}_2$ and $\widetilde{\mathcal{T}}$ constitute a trivial foliated $\mathbb{R}$-bundle.
\end{thm}

We stress that the foliations $\mathcal{F}_1$ and $\mathcal{F}_2$ are indeed $\mathbb{R}$-covered by Corollary~\ref{uniform-Gromov-Rcovered}. 

By definition, a leaf $F \in \widetilde{\mathcal{F}}_1$ with a degenerate limit set for $\widetilde{\mathcal{G}}$ admits a unique accumulation point $\xi_F \in \partial_\infty F$ for all $\widetilde{\mathcal{G}}_F$-leaves. We designate $\xi_F$ as the \emph{distinguished ideal point} of the leaf $F\in\widetilde{\mathcal{F}}_1$.

\subsection{Weak-stable foliation of a topological Anosov flow}

In this subsection, we construct a topological Anosov flow whose weak-stable foliation aligns precisely with $\widetilde{\mathcal{F}}_1$.

\begin{lem}
    Let $\mathcal{F}_1$ be a foliation with non-compact Gromov hyperbolic leaves and $\mathcal{G}$ a one-dimensional subfoliation. If every leaf of $\widetilde{\mathcal{F}}_1$ has a degenerate limit set for $\widetilde{\mathcal{G}}$, then every leaf of $\mathcal{F}_1$ is either a cylinder or a plane.
\end{lem}

\begin{proof}
    For each $F \in \widetilde{\mathcal{F}}_1$, deck transformations fixing $F$ act isometrically and preserve $\widetilde{\mathcal{G}}_F$. Since the limit set of $\widetilde{\mathcal{G}}_F$ coincides with the distinguished ideal point $\xi_F$, these transformations must fix $\xi_F$. Such isometries cannot be elliptic. Parabolic isometries are excluded due to the uniformly bounded radius of $\mathcal{F}_1$'s foliation charts. Thus, all such deck transformations are hyperbolic isometries.  

    Let $F_0 = \pi(F) \in \mathcal{F}_1$ be the projected leaf in $M$. Each generator of $\pi_1(F_0)$ corresponds to an invariant geodesic (its axis) in $F$. Distinct generators yield axes with disjoint ideal endpoints, forcing $\pi_1(F_0)$ to have at most one generator and to be abelian. Since $\mathcal{F}_1$ has no compact leaves, $F_0$ must be a cylinder or a plane.
\end{proof}

A (topological) Anosov flow in a closed 3-manifold is \emph{skewed} if it is not orbit-equivalent to a suspension Anosov flow.

The following result is indicated by \cite[Theorem 5.5.8]{Calegari06}, see also \cite{Calegari07book, BFP23collapsed, FP23intersection}.
\begin{prop}\label{construct_skew_flow}
	Let $\mathcal{F}_1$ be a minimal foliation by non-compact Gromov hyperbolic leaves and $\mathcal{G}$ be a one-dimensional subfoliation. Assume that every leaf of $\widetilde{\mathcal{F}}_1$ has degenerate limit set for $\widetilde{\mathcal{G}}$. Then, there exists a topological Anosov flow preserving $\mathcal{F}_1$ as its weak-stable foliation. If moreover $\mathcal{F}_1$ is $\mathbb{R}$-covered and $\pi_1(M)$ is not virtually solvable, then $\mathcal{F}_1$ is the weak-stable foliation of a skewed $\mathbb{R}$-covered topological Anosov flow.
\end{prop}
\begin{proof}
    For each leaf \( F \in \widetilde{\mathcal{F}}_1 \), we construct a geodesic foliation in \( F \) directed toward the distinguished ideal point \( \xi_F \). This foliation induces a non-vanishing unit vector field on \( \widetilde{M} \), aligned with geodesics terminating at distinguished ideal points across \( \widetilde{\mathcal{F}}_1 \)-leaves. The resultant vector field generates a non-singular flow \( \widetilde{\phi}_t \). The union \( \bigcup_{F\in \widetilde{\mathcal{F}}_1} \xi_F \) forms a \( \pi_1(M) \)-invariant continuous curve in the cylinder at infinity \( \mathcal{A}_{\mathcal{F}_1} \), ensuring the \( \pi_1(M) \)-equivariance of \( \widetilde{\phi}_t \). Projection yields a non-singular flow \( \phi_t: M \to M \), which is expansive by \cite[Theorem 5.5.8]{Calegari06}.  
    
    Theorem~\ref{flow-defn} establishes \( \phi_t \) as a topological Anosov flow preserving \( \mathcal{F}_1 \), with \( \mathcal{F}_1 \) functioning as its weak-stable foliation by construction. The \( \mathbb{R} \)-covered property of \( \mathcal{F}_1 \) ensures \( \phi_t \) is a \( \mathbb{R} \)-covered topological Anosov flow and thus it is transitive by \cite{Fenley94}. Following \cite{Shannon}, \( \phi_t \) becomes orbit-equivalent to a smooth \( \mathbb{R} \)-covered Anosov flow \( \psi_t \). As \( \pi_1(M) \) is not virtually solvable, \( \psi_t \) cannot be orbit equivalent to a suspension Anosov flow, thereby classifying it as a skewed \( \mathbb{R} \)-covered Anosov flow. Consequently, \( \phi_t \) constitutes a skewed \( \mathbb{R} \)-covered topological Anosov flow with \( \mathcal{F}_1 \) as its weak-stable foliation.
\end{proof}

\subsection{Uniformly equivalent foliations}

Throughout this subsection, we will assume $\mathcal{F}_1$, $\mathcal{F}_2$, and $\mathcal{G}$ are given as in the assumption of Theorem~\ref{degenerate-limit-set}.

We will make use of the following property of skewed Anosov flows, see \cite{Fenley94, BM24, FP23intersection}.

\begin{prop}\label{skew_prop}
	Let $\psi$ be a skewed $\mathbb{R}$-covered (topological) Anosov flow and $\mathcal{W}$ be its weak-stable foliation. Consider any non-trivial deck transformation $\gamma$ that has at least one fixed leaf of $\widetilde{\mathcal{W}}$. Then, it has a countable number of fixed leaves $\{L_n\}_{n\in \mathbb{Z}}\subset \widetilde{\mathcal{W}}$ ordered by the transverse orientation. Moreover, for each $n\in \mathbb{Z}$, there are unique $\gamma$-invariant orbits $\alpha\subset L_n$ and $\beta\subset L_{n+1}$, whose projections in $M$ are freely homotopic to each other with reversed orientation.
\end{prop}

\begin{lem}\label{L-E}
    If some leaf of $\mathcal{F}_2$ has fundamental group distinct from $\pi_1(M)$, then there exists a non-trivial $\gamma\in \pi_1(M)$ and a pair of leaves $L\in \widetilde{\mathcal{F}}_1$ and $E\in \widetilde{\mathcal{F}}_2$ such that $\gamma(L)=L$ and $\gamma(E)\neq E$.
\end{lem}
\begin{proof}
    By assumption, we assume the existence of a leaf $E_0\in \mathcal{F}_2$ such that $\pi_1(E_0)\neq \pi_1(M)$. Denote by $E\in \widetilde{\mathcal{F}}_2$ the lifted leaf of $E_0$ in the universal cover $\widetilde{M}$. Then, there exists a non-trivial deck transformation $\gamma\in \pi_1(M)$ with $\gamma(E)\neq E$. As shown in the proof of Lemma~\ref{construct_skew_flow}, $\mathcal{F}_1$ is the weak-stable foliation of a transitive topological Anosov flow $\phi_t$. By \cite{Shannon}, $\phi_t$ is orbit equivalent to a transitive smooth Anosov flow $\psi_t$. In particular, the closed orbits collections of $\phi_t$ and $\psi_t$ correspond to the same subgroup of $\pi_1(M)$. The homotopy classes of closed orbits of $\psi_t$ generate $\pi_1(M)$ \cite{Adachi87}, so do the closed orbits of $\phi_t$. Then, there exists a closed orbit $o$ of $\phi_t$ whose homotopy class contains the element $\gamma\in \pi_1(M)$. For a leaf $L\in \widetilde{\mathcal{F}}_1$ containing a lift of $o$, we have that $\gamma(L)=L$ and thus finish the proof.
\end{proof}

While we invoke Lemma~\ref{construct_skew_flow} in the preceding argument, the non-virtually solvable fundamental group assumption serves exclusively to establish the skewed property of constructed flow for subsequent applications. This constraint is non-essential to the immediate result but critical for the following one.

We say that two foliations $\mathcal{F}_1, \mathcal{F}_2$ are \emph{uniformly equivalent} if each pair of lifted leaves $L_1\in \widetilde{\mathcal{F}}_1$ and $L_2\in \widetilde{\mathcal{F}}_2$ has bounded Hausdorff distance.

\begin{lem}\label{uniform-equivalent}
	Let $\pi_1(M)$ be not virtually solvable. Assume that there exists a non-trivial $\gamma\in \pi_1(M)$ and a pair of leaves $L\in \widetilde{\mathcal{F}}_1$ and $E\in \widetilde{\mathcal{F}}_2$ such that $\gamma(L)=L$ and $\gamma(E)\neq E$. Then, $\mathcal{F}_1$ is uniformly equivalent to $\mathcal{F}_2$.
\end{lem}
\begin{proof}
	By the minimal and $\mathbb{R}$-covered properties of $\mathcal{F}_1$ and $\mathcal{F}_2$, we are sufficient to assume that every leaf of $\widetilde{\mathcal{F}}_2$ intersects every leaf of $\widetilde{\mathcal{F}}_1$ as shown in \cite[Theorem 1.1]{BFP25transverse} and \cite{Thurston97}. Since $\gamma$ induces an orientation preserving action on the leaf space of $\widetilde{\mathcal{F}}_2$, there exists another leaf $F\in \widetilde{\mathcal{F}}_2$ distinct from $E$ satisfying that $\gamma(F)\neq F$. By Lemma~\ref{construct_skew_flow} and Proposition~\ref{skew_prop}, there exists a $\gamma$-invariant leaf $L'\in \widetilde{\mathcal{F}}_1$ distinct from $L$ such that the leaf interval of $\widetilde{\mathcal{F}}_1$ between $L$ and $L'$ has no other $\gamma$-invariant leaves. As shown in Proposition~\ref{skew_prop}, the $\gamma$-axes $\alpha\in L$ and $\beta\in L'$ project to two freely homotopic loops with reversed orientations. Since the distinguished ideal point $\xi_L$ is fixed by all deck transformations leaving $L$ invariant, it is one of two endpoints of the geodesic $\alpha$. Analogously, the geodesic $\beta$ has an ideal point $\xi_{L'}\in \partial_{\infty}L'$. The fact $\pi(\alpha)$ is freely homotopic to $-\pi(\beta)$ implies that the actions of $\gamma$ on $\alpha$ and $\beta$ have reversed orientations. We assume that $\xi_L$ is an attractor of $\gamma$ and $\xi_{L'}$ is a repeller of $\gamma$. 
	
	Notice by assumption that $L\cap E\neq \emptyset$ and $L'\cap F\neq \emptyset$. Up to changing curves in the same homotopy classes as $\pi(\alpha)$ and $\pi(\beta)$, we assume that $E\cap \alpha\neq\emptyset$ and $F\cap \beta\neq \emptyset$. As $\mathcal{F}_2$ is $\mathbb{R}$-covered, the leaf space of $\widetilde{\mathcal{F}}_2$ between $E$ and $F$, denoted by $I_{[E,F]}$, is a connected interval. Up to taking $\gamma^{-1}$, we can assume that $\gamma(I_{[E,F]})$ is a strict sub-interval of $I_{[E, F]}$, since both boundary leaves $E$ and $F$ are not $\gamma$-invariant. It implies that there exist two $\gamma$-invariant leaves $E_0, F_0\in I_{[E,F]}$ that are accumulated by $\gamma^i(E), \gamma^i(F)$, respectively, as $i\rightarrow +\infty$, where the case $E_0=F_0$ might happen. Then, the leaf $E_0$ is disjoint with $L$. Otherwise, the intersection $E_0\cap L$ would contain a curve homotopic to $\alpha$. However, each leaf of $\widetilde{\mathcal{G}}_L$ has a unique ideal point $\xi_L$, which means that no leaf in $E_0\cap L$ has bounded distance from $\alpha$. Thus, we conclude that $E_0\cap L=\phi$, and similarly $F_0\cap L'=\phi$. Note that $I_{[E, F]}$ is a bounded interval, so is $I_{[E_0,F_0]}$. Moreover, since $\mathcal{F}_1$ is $\mathbb{R}$-covered, the leaf space between $L$ and $L'$ is also a bounded interval. As the leaves of $\widetilde{\mathcal{F}}_1$ are properly embedded, $L$ splits $\widetilde{M}$ into two half spaces, denoted by $R_L^+$ and $R_L^-$. Let $R_L^+$ be the one containing $L'$. Analogously, $L'$ splits two half spaces $R_{L'}^+$ and $R_{L'}^-$, where $R_{L'}^-$ denotes the one containing $L$. By the fact that $E\cap L\neq \emptyset$ and $F\cap L'\neq \emptyset$, both $E_0$ and $F_0$ are contained in the region $R_L^+\cap R_{L'}^-$. It turns out that the Hausdorff distance $d_H(E_0, L)$ is bounded by $d_H(L, L')$. Thus by uniformness of $\mathcal{F}_1$ and $\mathcal{F}_2$, for any leaves $P\in \widetilde{\mathcal{F}}_1$ and $Q\in \widetilde{\mathcal{F}}_2$, we have that
	$$
	d_H(P, Q)\leq d_H(P, L)+d_H(L,E_0)+d_H(E_0, Q)<+\infty.
	$$
	Therefore, $\mathcal{F}_1$ and $\mathcal{F}_2$ are uniformly equivalent. We finish the proof.
\end{proof}

The following proposition is independent of our preceding results in this section and it may have its own interest.

\begin{prop}\label{degenerate-distinct}
    Let $\mathcal{F}_1$ and $\mathcal{F}_2$ be transverse uniform minimal foliations that are uniformly equivalent to each other. Then, both $\widetilde{\mathcal{F}}_1$ and $\widetilde{\mathcal{F}}_2$ have leaves with non-degenerate limit set for $\widetilde{\mathcal{G}}$.
\end{prop}

\begin{proof}
	We will only show it for $\widetilde{\mathcal{F}}_1$ since the symmetric argument applies to $\widetilde{\mathcal{F}}_2$. Suppose that every leaf of $\widetilde{\mathcal{F}}_1$ has degenerate limit set for $\widetilde{\mathcal{G}}$. Then, each $F\in \widetilde{\mathcal{F}}_1$ admits a distinguished ideal point $\xi_F\in \partial_{\infty}F$ shared by all leaves of $\widetilde{\mathcal{G}}_F$. Moreover, the point $\xi_F$ is fixed by all deck transformations preserving the leaf $F$.
	
	For any pair of leaves $L, F\in \widetilde{\mathcal{F}}_1$, let $a_0:= d_H(L,F)>0$ be the Hausdorff distance between $L$ and $F$, which is finite by the uniformness of $\mathcal{F}_1$. We define a map $\tau_{L,F}: L\rightarrow F$ as $\tau_{L,F}(x)=y$, where the point $y\in F$ is chosen to satisfy $d(x, y)\leq a_0$. This map is well-defined up to an error by the $\mathbb{R}$-covered property of $\mathcal{F}_1$. Indeed, if $y_1, y_2\in F$ are two points with $d(x, y_i)\leq a_0$ for $i=1,2$, then the triangle inequality implies $d(y_1, y_2)\leq 2a_0$. Since the foliation $\mathcal{F}_1$ is $\mathbb{R}$-covered, for any $\epsilon>0$, there is $\delta(\epsilon)>0$ so that any pair of points $p,q\in \widetilde{M}$ in the same leaf $F\in \widetilde{\mathcal{F}}_1$ with $d(p,q)<\epsilon$ satisfy $d_F(p,q)<\delta(\epsilon)$, see Proposition \ref{AB}. Then there is a constant $\delta(2a_0)>0$ depending only on $a_0$ such that $d_F(y_1, y_2)\leq \delta(2a_0)$. It means that all points in $F$ of distance not larger than $a_0$ from $x$ are contained in the $\delta(2a_0)$-neighborhood of $y_1$. So up to a $\delta(2a_0)$-neighborhood, the image of $x$ under $\tau_{L,F}$ is uniquely determined. Thus, in this sense, the map $\tau_{L,F}$ is coarsely well-defined. Moreover, this map is a quasi-isometry, see \cite[Proposition 3.4]{Fenley02} for a complete account. 
	
	The map $\tau_{L,F}$ can be extended to the ideal boundary $\partial_{\infty}L$ as a homeomorphism between $\partial_{\infty}L$ and $\partial_{\infty}F$, denoted also by $\tau_{L,F}$. We only show that the extension map on the boundaries of any pair of leaves is well-defined. Indeed, given any ideal point $z\in \partial_{\infty}L$, we consider the geodesic ray starting at any point $x\in L$ with ideal point $z$. The image of this geodesic under $\tau_{L,F}$ is a ray contained in the $2a_0$-neighborhood of a geodesic starting at the point $\tau_{L,F}(x)$ in $F$ by definition. This implies that $\tau_{L,F}$ maps any geodesic in $L$ to a quasi-geodesic in $F$, and thus it maps any ideal point in $\partial_{\infty}L$ to a unique ideal point in $\partial_{\infty}F$. Therefore, the extension map $\tau_{L,F}: L\cup \partial_{\infty}L\rightarrow F\cup \partial_{\infty}F$ is well-defined.
	
	We are going to produce a vertical foliation on the cylinder at infinity $\mathcal{A}$ by following \cite{Fenley02}. Given any leaf $L\in \widetilde{\mathcal{F}}_1$ and any ideal point $z\in \partial_{\infty}L$, the union $\sigma_z:=\bigcup\limits_{F\in \widetilde{\mathcal{F}}_1} \tau_{L,F}(z)$ is a continuous curve on $\mathcal{A}$. We first give the interpretation of the continuity of $\sigma_z$. Let $F\in \widetilde{\mathcal{F}}_1$ be any leaf and $T$ be any transversal to $\widetilde{\mathcal{F}}_1$ through a point $y\in F$. Let $L_i\in \widetilde{\mathcal{F}}_1$, $i\in \mathbb{N}$, be a sequence of leaves intersecting $T$ such that intersection points $x_i\in L_i\cap T$ converge to $y$ as $i$ goes to infinity. The continuity of $\sigma_z$ means that the ideal points $\tau_{L,L_i}(z)$ converge to $\tau_{L,F}(z)$ as $x_i$ converges to $y$. Now let us prove the continuity of $\sigma_z$. We can assume the leaves $L_i$ are in the interval $[L,F]$ of the leaf space. Let $\alpha_0$ be any geodesic ray in $L$ with ideal point $z$, $\alpha_i$ be the geodesic rays in $L_i$ from $x_i$ to $\tau_{L,L_i}(z)$, and $\beta$ be the geodesic ray in $F$ accumulated by $\alpha_i$. By the uniformness of $\mathcal{F}_1$ and the quasi-isometric property of $\tau_{L,L_i}$ and $\tau_{L,F}$, the Hausdorff distance $d_H(\tau_{L,L_i}(\alpha_0),\tau_{L,F}(\alpha_0))$ between two quasi-geodesic rays is finite for any $i \in\mathbb{N}$. As the quasi-geodesic ray $\tau_{L,L_i}(\alpha_0)$ has the same ideal point as $\alpha_i$, the Hausdorff distance $d_H(\tau_{L,L_i}(\alpha_0), \alpha_i)$ is finite. Since $\beta$ is accumulated by $\alpha_i$, the Hausdorff distance $d_H(\beta, \alpha_i)$ is finite. It turns out that $$d_H(\beta, \tau_{L,F}(\alpha_0))\leq d_H(\beta, \alpha_i)+d_H(\alpha_i, \tau_{L,L_i}(\alpha_0))+d_H(\tau_{L,L_i}(\alpha_0),\tau_{L,F}(\alpha_0))$$ is also finite. Therefore, the quasi-geodesic ray $\tau_{L,F}(\alpha_0)$ and the geodesic ray $\beta$ have the common ideal point $\tau_{L,F}(z)$, and thus the continuity of $\sigma_z$ follows. Note that each pair of curves $\sigma_z, \sigma_{z'}$ are disjoint unless $z=z'\in \partial_{\infty}L$. Then we obtain a vertical foliation $\bigcup\limits_{z\in \partial_{\infty}L}\sigma_z$ of $\mathcal{A}_{\mathcal{F}_1}$. We can construct a universal circle as the quotient $\mathcal{U}=\mathcal{A}/\sim$, where two ideal points $z_1\in \partial_{\infty}L$ and $z_2\in \partial_{\infty}F$ are equivalent if $z_2=\tau_{L,F}(z_1)$. This gives a canonical way to identify all the ideal circle of $\widetilde{\mathcal{F}}_1$.
	
	Since $\mathcal{F}_1$ and $\mathcal{F}_2$ are uniformly equivalent, we can define analogously a coarsely well-defined map $\tau_{L, E}: L\rightarrow E$ for any leaves $L\in \widetilde{\mathcal{F}}_1$ and $E\in \widetilde{\mathcal{F}}_2$ such that $d(x, \tau_{L, E}(x))\leq d_H(L, E)<+\infty$ for any point $x\in L$. As explained above, this map is quasi-isometric and extends to a homeomorphism between the boundaries $\tau_{L, E}:  \partial_{\infty}L\rightarrow \partial_{\infty}E$. Picking leaves $L$ and $E$ with $L\cap E\neq\emptyset$, we denote by $l\subset L\cap E$ a leaf (non-necessarily unique) of $\widetilde{\mathcal{G}}$ in the intersection. Let $l_L^{\pm}\in \partial_{\infty}L$ and $l_E^{\pm}\in \partial_{\infty}E$ be the ideal points of $l$ in two ideal circles. By assumption, we know $l_L^+=l_L^-=\xi_L$. The map $\tau_{L,E}$ can be chosen as the identity map form $l\subset L$ to $l\subset E$. The fact that $\tau_{L,E}:\partial_{\infty}L\rightarrow \partial_{\infty}E$ is a homeomorphism implies $l_E^+=l_E^-\in \partial_{\infty}E$. The leaf $l$ splits $E$ into two connected half planes, one of which has a single ideal point $\theta:=l_E^+=l_E^-$, denoted by $E_l^+$. Then, every leaf of $\widetilde{\mathcal{G}}_E$ contained in $E_l^+$ has to accumulate at $\theta$. By the transversality of $\mathcal{F}_1$ and $\mathcal{F}_2$, there exists a leaf interval $I$ in the leaf space of $\widetilde{\mathcal{F}}_1$ bounded by $L$ such that for any leaf $L'\in I$, the intersection $L'\cap E$ contains a $\widetilde{\mathcal{G}}$-leaf $l'\in E_l^+$ with $l'^+_E=l'^-_E=\theta\in \partial_{\infty}E$. Since $l'\in \widetilde{\mathcal{G}}_{L'}$, we have that $l'^+_{L'}=l'^-_{L'}=\xi_{L'}$, where $\xi_{L'}\in \partial_{\infty}L'$ is the distinguished ideal point of $\widetilde{\mathcal{G}}_{L'}$. 
	
	Notice that $\tau_{L, E}(\xi_{L})=\theta$. By the choice of $I$, for any leaf $L'\in I$, we have that $\tau_{L', E}(\xi_{L'})=\theta$ by making $\tau_{L', E}$ as the identity on $l'\subset L'\cap E$. It follows that $\tau_{L,L'}(\xi_L)=\xi_{L'}$ for any $L'\in I$.

	Recall that each leaf $F\in \widetilde{\mathcal{F}}_1$ has a distinguished ideal point $\xi_F\in \partial_{\infty}F$ which is the common ideal point of all leaves of $\widetilde{\mathcal{G}}_F$. The union of all distinguished ideal points $\xi:=\bigcup\limits_{F\in \widetilde{\mathcal{F}}}\xi_F$ is a $\pi_1(M)$-invariant continuous curve in $\mathcal{A}_{\mathcal{F}_1}$. Its continuity follows from the continuity of $\widetilde{\mathcal{G}}$. 
	By minimality of $\mathcal{F}_1$, the leaf interval $I$ extends to the whole leaf space of $\mathcal{F}_1$ under all deck transformations. Since $\xi$ is $\pi_1(M)$-invariant, we have that $\tau_{L,F}(\xi_L)=\xi_F$ for any leaf $F\in \widetilde{\mathcal{F}}_1$. This means that the curve $\xi$ is a leaf of the vertial foliation of $\mathcal{A}_{\mathcal{F}_1}$ constructed above, which corresponds to a single point in the universal circle $\mathcal{U}$. This single point is a global fixed point of $\pi_1(M)$ since $\xi$ is a $\pi_1(M)$-invariant curve. However, by \cite{FP20minimal}, the fundamental group $\pi_1(M)$ is a minimal action on the universal circle $\mathcal{U}$. We arrive to a contradiction, and thus finish the proof.	
\end{proof}

In the preceding argument, the universal circle construction for $\widetilde{\mathcal{F}}_1$ plays a central role. This framework extends naturally to all uniform foliations with Gromov hyperbolic leaves (see \cite{Thurston97, Calegari00, Fenley02}). We explicitly implement this construction here to establish a correspondence between ideal boundaries of uniformly equivalent foliations.

\subsection{Trivial foliated bundle}

\begin{lem}\label{trivialholonomy-lem}
    Let $\mathcal{F}$ be a foliation whose leaves have the same fundamental group. Then $\mathcal{F}$ has trivial holonomy.
\end{lem}
\begin{proof}
    Let $E\in \mathcal{F}$ be any leaf and $\alpha$ be any loop expressing a non-trivial homotopy class in $E$. For any transversal $\tau$ to $\mathcal{F}$ joining the loop $\alpha$, any $\mathcal{F}$-holonomy map sending $\tau$ to itself must be the identity. Otherwise, there would be a leaf $F\in\mathcal{F}$ intersecting $\tau$ so that the homotopy class given by $\alpha$ is not contained in $\pi_1(F)$. It contradicts the assumption that $\pi_1(F)=\pi_1(E)$. This implies that the holonomy group of $\mathcal{F}$ is trivial.
\end{proof}

\begin{thm}\cite[Chapter VIII, Theorem 2.2.1]{HectorHirsch}\label{product}
    Let $\mathcal{F}$ be a codimension-one foliation without holonomy on a closed manifold. Then for any transverse one-dimensional foliation $\mathcal{T}$, the lifted foliations $\widetilde{\mathcal{F}}$ and $\widetilde{\mathcal{T}}$ constitute a trivial foliated $\mathbb{R}$-bundle.
\end{thm}

\medskip
Now, we can complete the proof of Theorem~\ref{degenerate-limit-set}.

\begin{proof}[Proof of Theorem~\ref{degenerate-limit-set}]
    Suppose that there exists a leaf of $\mathcal{F}_2$ whose fundamental group is not equal to $\pi_1(M)$. Lemma~\ref{L-E} provides the existence of leaves $L\in \widetilde{\mathcal{F}}_1$ and $E\in \widetilde{\mathcal{F}}_2$ such that $\gamma(L)=L$ and $\gamma(E)\neq E$ for some non-trivial $\gamma\in \pi_1(M)$. Then, under the condition that $\pi_1(M)$ is not virtually solvable, the foliations $\mathcal{F}_1$ and $\mathcal{F}_2$ are uniformly equivalent by Lemma~\ref{uniform-equivalent}. As shown in Proposition~\ref{degenerate-distinct}, $\widetilde{\mathcal{F}}_1$ must have a leaf with non-degenerate limit set for $\widetilde{\mathcal{G}}$. This contradicts to the assumption. Thus, the holonomy group of $\mathcal{F}_2$ is trivial by Lemma~\ref{trivialholonomy-lem}. Under the trivial holonomy condition, the trivial foliated $\mathbb{R}$-bundle for $\widetilde{\mathcal{F}}_2$ and any transverse foliation $\widetilde{\mathcal{T}}$ follows directly from Proposition~\ref{product}. Hence, we finish the proof.
\end{proof}

\section{Analysis on non-separated leaves}\label{section-nonHausdorff}

For a 2-dimensional foliation $\mathcal{F}$ with Gromov hyperbolic leaves and subfoliation $\mathcal{G}$, each $F \in \widetilde{\mathcal{F}}$ determines a simply connected, separable one-dimensional leaf space from $\widetilde{\mathcal{G}}_F$. Recall from Section~\ref{subsection-coherent} that non-Hausdorff leaf spaces of $\widetilde{\mathcal{G}}_F$ correspond to pairs of non-separated leaves lacking Hausdorff neighborhood bases. Specifically, such leaf pairs are simultaneously accumulated by a common sequence of approximating leaves. This section investigates obstructions to Hausdorffness in $\widetilde{\mathcal{G}}_F$'s leaf space structure.

Our principal result of this section is:

\begin{thm}\label{closed-leaf-nonHausdorff}
    Let $\mathcal{F}_1$, $\mathcal{F}_2$ be transverse minimal $\mathbb{R}$-covered foliations with non-compact Gromov hyperbolic leaves in a closed 3-manifold $M \neq \mathbb{T}^3$, and let $\mathcal{G} = \mathcal{F}_1 \cap \mathcal{F}_2$. Suppose there exists $l_0 \in \widetilde{\mathcal{G}}_{L_0}$ with distinct ideal endpoints $l_0^+ \neq l_0^- \in \partial_\infty L_0$ for some $L_0 \in \widetilde{\mathcal{F}}_1$. If there is a leaf $L \in \widetilde{\mathcal{F}}_1$ exhibiting a non-Hausdorff leaf space for $\widetilde{\mathcal{G}}_L$, then $\mathcal{G}$ contains a closed leaf.
\end{thm}

Note that the minimality condition serves exclusively to apply Proposition~\ref{distance_geodesic to leaf_dichotomy} universally across $\widetilde{\mathcal{F}}_1$-leaves. By this proposition, the existence of a leaf with distinct ideal points implies that for every $F \in \widetilde{\mathcal{F}}_1$ and $l \in \widetilde{\mathcal{G}}_F$, the associated geodesic $l^* \subset F$ connecting $l$'s endpoints lies within a $C_l$-neighborhood of $l$ for some $C_l > 0$. Corollary~\ref{distance-corollary} explicitly precludes the second case in Proposition~\ref{distance_geodesic to leaf_dichotomy}.

We recall contracting directions from Definition~\ref{def_contracting}. An ideal point $\xi \in \partial_\infty F$ is \emph{contracting} for $F \in \widetilde{\mathcal{F}}$ if $F$ contains a contracting direction toward $\xi$. By contrast, $\xi$ is \emph{non-expanding} in $\partial_\infty F$ if for every $\epsilon > 0$ and $x \in F$, there exist $\delta > 0$ and a $\widetilde{\mathcal{F}}$-transversal $\beta$ through $x$ of length less than $\delta$ such that for all $L \in \widetilde{\mathcal{F}}$ intersecting $\beta$ and any point $y$ along the geodesic ray from $x$ to $\xi$, the distance $d(y, L)$ is bounded by $\epsilon$. An ideal point is \emph{expanding} if it fails to be non-expanding.

The transversal $\beta$ in this definition can be selected to pass through $x$ in its interior when $\mathcal{F}$ lacks compact leaves, as indicated by Theorem~\ref{Thurston}. This constitutes the sole application of the compact-leaf exclusion condition in our main theorem in this section.

\subsection{Transverse expanding direction}

Given a ray $l$ contained in $L\cap E$ for $L\in \widetilde{\mathcal{F}}_1$ and $E\in \widetilde{\mathcal{F}}_2$, the ideal point of $l$ in $L$, say $l^+_L\in\partial_{\infty}L$, uniquely defines an ideal point $l^+_E\in\partial_{\infty}E$ since the leaves of both $\widetilde{\mathcal{F}}_1$ and $\widetilde{\mathcal{F}}_2$ are properly embedded in the universal cover $\widetilde{M}$. Initially, no canonical identification exists between the ideal circles of  $L$ and $E$ — a property typically arising when $\mathcal{F}_1$ and $\mathcal{F}_2$ are uniformly equivalent foliations.  

The subsequent lemma characterizes transverse expansion in non-separated directions. This generalizes \cite[Proposition 3.9]{FP23_transverse}, which proves non-separated rays exclusively accumulate at non-marker points in their transverse leaves in case where the manifold is the unit tangent bundle of a closed hyperbolic surface.

\begin{lem}\label{expanding_transverse}
	Let $s\in \widetilde{\mathcal{G}}_{L}$ be a non-separated leaf in $L\in \widetilde{\mathcal{F}}_1$, and $s_L^+\in \partial_{\infty}L$ be the ideal point of a non-separated ray of $s$. Then $s_E^+\in \partial_{\infty}E$ is an expanding ideal point in $\partial_{\infty}E$, where $E\in \widetilde{\mathcal{F}}_2$ is the leaf through $s$.
\end{lem}
\begin{proof} 
	Given a point $x_0\in s$, we denote by $l$ the ray starting at $x_0$ towards $s_L^+$ non-separated from another leaf $s'\in \widetilde{\mathcal{G}}_L$. Let $r_n\in \widetilde{\mathcal{G}}_L$ be a sequence of leaves converging to $s$ and $s'$ simultaneously as $n\rightarrow +\infty$.
	
	Suppose that $s_E^+$ is a non-expanding ideal point in $\partial_{\infty}E$. By the definition of non-expanding ideal point, for any $\epsilon>0$, there is $\delta>0$ such that if a leaf $E'\in \widetilde{\mathcal{F}}_2$ intersects a transversal through $x_0$ of length less than $\delta$, then every point in the geodesic ray in $E$ from $x_0$ to $s_E^+$ has distance less than $\epsilon$ from the leaf $E'$. By Proposition \ref{distance_geodesic to leaf_dichotomy}, there is a constant $C>0$ such that the associated geodesic ray $l^*\subset E$ with the same endpoints of $l$ is contained in the $C$-neighborhood of $l$. It implies that
	for any sequence of points $z_m\in l^*$ converging to $s_E^+$, there is a sequence of points $x_m\in l$ of distance $d_E(x_m, z_m)\leq C$. Up to choosing a smaller $\epsilon$, by continuity of $\widetilde{\mathcal{F}}_2$, there exists a constant $\rho=\rho(\epsilon, C)>0$ such that $d(x_m, E')\leq \rho$ for each $m\in \mathbb{N}$. As $\epsilon$ and thus $\rho$ are small enough, by the transversality of $\widetilde{\mathcal{F}}_1$ and $\widetilde{\mathcal{F}}_2$, there is a sequence of points $y_m\in L\cap E'$ so that $d_L(x_m, y_m)\leq \omega$ for $m\in \mathbb{N}$ and some constant $\omega>0$. One can see that $\omega$ becomes smaller and smaller as $\epsilon\rightarrow0$. Thus, without loss of generality, we can assume that the $\widetilde{\mathcal{G}}_L$-leaf through $y_m$ intersects a transversal to $\widetilde{\mathcal{G}}$ through $x_m$ in $L$.
    
    Denote by $E_n\in \widetilde{\mathcal{F}}_2$ the leaf containing $r_n$ for each $n\in \mathbb{N}$. There is $N\in \mathbb{N}$ such that for any $n\geq N$, $E_n$ contains a point $y_{m, n}\in L\cap E_n$ given as above: $d_L(x_m, y_{m, n})\leq \omega$ and the leaf $l_{m,n}\in \widetilde{\mathcal{G}}_L$ through $y_{m,n}$ intersects a transversal to $\widetilde{\mathcal{G}}_L$ through $x_m$ in the leaf $L$ for all $m\in \mathbb{N}$. As $r_n$ does not accumulate at $s_L^+$, the distance $d_L(x_m, r_n)$ goes to infinity as $m$ goes to infinity. Thus, we have that $l_{m,n}$ is distinct from $r_n$ for $m$ large enough, while both of them are contained in $L\cap E_n$.

	Let $u, u'\subset L$ be two transversals to $\widetilde{\mathcal{G}}_L$ that intersects $s, s'$, respectively. Since $s$ and $s'$ are two leaves accumulated by $r_n$ simultaneously, there is $N'\in \mathbb{N}$ such that for any $n\geq N'$, the leaf $r_n$ intersects both $u$ and $u'$. Denote by $R\subset L$ the region bounded by leaves $s$, $s'$, $r_n$ and tranversals $u$, $u'$. Then, the leaf $l_{m,n}$ is entirely contained in $R$. Otherwise, it could only escape from $R$ through either $u$ or $u'$. However, as both $l_{m,n}$ and $r_n$ are contained in $E_n$, either $u$ or $u'$ serves as a transversal to $\widetilde{\mathcal{F}}_2$ intersecting $E_n$ in two distinct points. This contradicts Theorem \ref{taut}.

    Recall that $l_{m,n}$ intersects a transversal, denoted by $\tau_m\subset L$, to $\widetilde{\mathcal{G}}_L$ through $x_m\in s$. If $l_{m, n}$ is non-separated from $s$, then $r_n$ intersects $\tau_m$ at least twice for $n$ large enough. This is impossible. If $l_{m, n}$ is separated from $s$, then there is a leaf $s''\in \widetilde{\mathcal{G}}_L$ non-separated from $s$ that separates $s$ from $l_{m,n}$. This is due to the fact that $l_{m,n}$ is contained in $R$ whose boundaries $s$ and $s'$ are non-separated leaves. It implies that $s''$ also intersects $\tau_m$. Then, as $r_n$ accumulate at $s$ and $s''$, there exists a large $n$ so that $r_n$ intersects $\tau_m$ at least twice for $n$ large. This gives a contradiction and thus finishes the proof.
\end{proof}

We stress that the lemma above holds in a more generality even if $\mathcal{F}_1$ and $\mathcal{F}_2$ are not $\mathbb{R}$-covered by assuming the conclusion of Corollary~\ref{distance-corollary}. Analogous argument also applies to provide the $\widetilde{\mathcal{F}}_2$-holonomy expansion along each non-separated ray.

\subsection{Deck invariant leaf}

Before going into discussion on pair of leaves, we first establish ideal points correspondence between ideal circles of $\widetilde{\mathcal{F}}_1$ and $\widetilde{\mathcal{F}}_2$. Regardless of $\widetilde{\mathcal{G}}$-leaves, for any $L\in \widetilde{\mathcal{F}}_1$, each ideal point in $\partial_{\infty}L$ defines an ideal point of $\partial_{\infty}E$ for some leaf $E\in \widetilde{\mathcal{F}}_2$ since the leaves of $\widetilde{\mathcal{F}}_1$ and $\widetilde{\mathcal{F}}_2$ are all properly embedded in $\widetilde{M}$. For ideal points lying in the limit set of $\widetilde{\mathcal{G}}_L$, each ideal point $l_L^+\in \partial_{\infty}L$ associated with a ray $l$ of $\widetilde{\mathcal{G}}_L$ uniquely determines an ideal point $l_E^+\in \partial_{\infty}E$ for the leaf $E\in \widetilde{\mathcal{F}}_2$ containing $l$. 

The following lemma establishes a bijective correspondence between ideal points of paired leaves. This will provide us convenience to transfer an ideal point from one circle to another. However, this structural compatibility alone remains insufficient to conclude uniform equivalence of $\mathcal{F}_1$ and $\mathcal{F}_2$.

\begin{lem}
    Let $l, l'\subset L\cap E$ be two rays of $\widetilde{\mathcal{G}}$ with the same ideal point $l_L^+=l'^+_L\in \partial_{\infty}L$, for $L\in \widetilde{\mathcal{F}}_1$ and $E\in \widetilde{\mathcal{F}}_2$. Then, they possess the same ideal point $l_E^+=l'^+_E\in \partial_{\infty}E$ in the leaf $E$.
\end{lem}
\begin{proof}
    Denote by $l_L^*$ and $l_L'^*$ be associated geodesic rays in $L$ joining endpoints of $l$ and $l'$, respectively. The fact $l_L^+=l_L'^+$ implies that $l_L^*$ and $l_L'^*$ share the same ideal point in $\partial_{\infty}L$. Then, they have bounded Hausdorff distance $d_H(l_L^*, l_L'^*)\leq C$ for some $C>0$. Corollary~\ref{distance-corollary} provides a constant $K>0$ such that $l_L^*$ and $l_L'^*$ are contained in $K$-neighborhoods of $l$ and $l'$, respectively, with respect to the induced distance $d_L$. We deduce that there exist two sequence of points $x_i\in l$, $y_i\in l'$ such that $x_i\rightarrow l_L^+$, $y_i\rightarrow l_L'^+$, and $d_L(x_i, y_i)\leq C+2K$ for all $i$. Note that $x_i$ and $y_i$ are all contained in $E$. As $\mathcal{F}_2$ is $\mathbb{R}$-covered, by Lemma~\ref{AB}, we have that $d_E(x_i, y_i)\leq C'$ for some constant $C'=C'(C, K)>0$. Since both $x_i$ and $y_i$ are chosen to escape from any compact set as $i\rightarrow$, the leaves $l$ and $l'$ possess the same ideal point $l_E^+=l'^+_E\in \partial_{\infty}E$. Thus, we finish the proof.
\end{proof}

\begin{prop}\label{expanding-both}
	Let $l\subset L\cap E$ be a non-separated ray of $\widetilde{\mathcal{G}}_L$ for $L\in \widetilde{\mathcal{F}}_1$ and $E\in \widetilde{\mathcal{F}}_2$, and $l_L^+\in \partial_{\infty}L$ be its ideal point. Assume that $l_L^+$ is an expanding (or contracting) ideal point in $\partial_{\infty}L$. Then there exists a deck transformation $\gamma$ fixing a leaf $F\in \widetilde{\mathcal{F}}_1$ and a leaf $s\in \widetilde{\mathcal{G}}_F$ with $\gamma(s)=s$. The same conclusion holds if we exchange the roles of $\widetilde{\mathcal{F}}_1$ and $\widetilde{\mathcal{F}}_2$.
\end{prop}
\begin{proof}
	Let $l^*_L$ be the associated geodesic ray of $l$ in $L$ joining the endpoints of $l$. By Corollary~\ref{distance-corollary}, there is a constant $C>0$ such that $l^*_L$ is contained in the $C$-neighborhood of $l$. In particular, for any sequence of points $z_n\in l^*_L$ converging to $l_L^+$, there is a sequence of points $x_n\in l$ of distance $d_L(z_n, x_n)\leq C$. We take $z_n$ so that for any $n\in \mathbb{N}$, the concatenate points $z_n$ and $z_{n+1}$ have distance at least $3C$. If the projection of $l$ by the covering map $\pi: \widetilde{M}\rightarrow M$ is closed in $M$, then there exists $\gamma\in \pi_1(M)$ such that $\gamma(l)=l$ and $\gamma(L)=L$, which finishes the proof. Now, we can assume that $l$ projects to a non-closed curve in $M$. Up to a subsequence, we have no loss of generality by assuming that the projections of $x_n$ are all distinct. By compactness of $M$, up to a subsequence again, we can assume that $\pi(x_n)$ is convergent. Then there is a sequence of deck transformations $\gamma_n$ such that $\gamma_n(x_n)$ is convergent. 
	
	Denote by $B_p$ a compact foliation box around a point $p\in \widetilde{M}$ of small diameter $\delta_0>0$. 
	Pick $\delta_0>0$ sufficiently small such that for any $p\in l$, each compact foliation box $B_p$ around $p$ of diameter smaller than $\delta_0$ satisfies that (i) every plaque of $\widetilde{\mathcal{F}}_1$ intersects every plaque of $\widetilde{\mathcal{F}}_2$ in a unique arc of $\widetilde{\mathcal{G}}$ in $B_p$; (ii) if there is a non-trivial deck transformation $\gamma$ fixing a leaf $F$ of $\widetilde{\mathcal{F}}_1$ intersecting $B_p$, then either $F$ is the unique $\gamma$-invariant leaf intersecting $B_p$, or all leaves intersecting $B_p$ are invariant by $\gamma$; (iii) if there is a non-trivial deck transformation $\gamma$ fixing a leaf $W$ of $\widetilde{\mathcal{F}}_2$ intersecting $B_p$, then either $W$ is the unique $\gamma$-invariant leaf intersecting $B_p$, or all leaves intersecting $B_p$ are invariant by $\gamma$.

	Pick a sufficiently small $\delta>0$ that is smaller than $\delta_0/6$. There is $N>0$ such that for any $n\geq N$, the point $\gamma_N^{-1}\circ\gamma_n(x_n)$ is contained in the $\delta$-neighborhood of $x_N$ and thus contained in a foliation box $B_{x_N}$ of diameter $\delta_0$. We claim that each deck transformation $\gamma_N^{-1}\circ\gamma_n$ is non-tirvial for any $n>N$. Otherwise, there is $n>N$ such that $\gamma_n=\gamma_N$. By the choice of $N$, we know that $\gamma_n(x_n)=\gamma_N(x_n)$ is contained in a small neighborhood of $\gamma_N(x_N)$. As the distance between $z_N$ and $z_n$ is larger than $3C$, the distance between $x_N$ and $x_n$ is larger than $C$. Then $\gamma_N(x_N)$ and $\gamma_N(x_n)$ have distance larger than $C>\delta_0$, which implies the claim.
	
	Suppose $l_L^+$ is an expanding ideal point in $\partial_{\infty}L$. Given an orientation in the leaf space of $\widetilde{\mathcal{F}}_1$, we denote by $L^+$ and $L^-$ two connected components of $\widetilde{M}$ split by $L$. Let $L_1\in L^+$, $L_2\in L^-$ be two leaves of $\widetilde{\mathcal{F}}_1$ satisfying that $d(x_N, L_i)<\delta$ for $i=1, 2$. Given a transversal segment $c_{x_N}$ to $\widetilde{\mathcal{F}}_1$ through $x_N$ so that it connects $L_1$ and $L_2$ and has length less than $2\delta$. As $l_L^+$ is an expanding ideal point in $\partial_{\infty}L$ accumulated by $l$, there is $n>N$ satisfying that $d(x_n, L_i)>3\delta$ for $i=1,2$. The transversal $c_{x_N}$ produces a transversal segment, denoted by $c_{x_n}$, through the point $x_n$ by the holonomy of $\widetilde{\mathcal{F}}_1$. The transversal $c_{x_n}$ has endpoints contained in $L_1$ and $L_2$, so its length is greater than $6\delta$. Since the foliation $\widetilde{\mathcal{F}}_1$ has Hausdorff leaf space, each point in $c_{x_N}$ and $c_{x_n}$ corresponds a leaf of $\widetilde{\mathcal{F}}_1$ in the interval $[L_1, L_2]$ of the leaf space. The fact that $\gamma_N^{-1}\circ\gamma_n(x_n)$ is contained in the $\delta$-neighborhood of $x_N$ implies that any point in the segment $c_{x_N}$ is contained in the $2\delta$-neighborhood of $\gamma_N^{-1}\circ\gamma_n(x_n)$. It turns out the existence of a $\gamma_N^{-1}\circ\gamma_n$-invariant leaf of $\widetilde{\mathcal{F}}_1$ which intersects $c_{x_N}$ and thus intersects $B_{x_N}$. By the choice of $\delta_0$, such an invariant leaf is unique, denoted by $L'$, since $\gamma_N^{-1}\circ\gamma_n$ expands the leaf interval $[L_1, L_2]$.
		
	Now, we are going to show a similar result for $\widetilde{\mathcal{F}}_2$. By Lemma \ref{expanding_transverse}, the ideal point $l_E^+$ of $l$ is an expanding ideal point in $\partial_{\infty}E$. Let $E_1\in E^+$, $E_2\in E^-$ be two leaves of $\widetilde{\mathcal{F}}_2$ satisfying that $d(x_N, E_i)<\delta$ for $i=1, 2$. There is a large $n>N$, which can be taken to be the same one as in the last paragraph, such that $d(x_n, E_i)>3\delta$ for $i=1, 2$. Denote by $u_{x_N}$ a transversal to $\widetilde{\mathcal{F}}_2$ through $x_N$ of length less than $2\delta$,  and $u_{x_n}$ a transversal to $\widetilde{\mathcal{F}}_2$ through $x_n$ of length greater than $6\delta$, both of whose endpoints are contained in $E_1$ and $E_2$. As $\mathcal{F}_2$ is a $R$-covered foliation, each point in $u_{x_N}$ and $u_{x_n}$ corresponds a leaf of $\widetilde{\mathcal{F}}_2$. Using the same argument as above, we conclude that the deck transformation $\gamma_N^{-1}\circ\gamma_n$ expands the leaf interval $[E_1, E_2]$. Therefore, there is a unique leaf $E'$ of $\widetilde{\mathcal{F}}_2$ invariant by $\gamma_N^{-1}\circ\gamma_n$. 
	
	Therefore, there are two leaves $L'\in \widetilde{\mathcal{F}}_1$ and $E'\in \widetilde{\mathcal{F}}_2$ invariant by the non-trivial deck transformation $\gamma:= \gamma_N^{-1}\circ\gamma_n$, both of which intersect the foliation box $B_{x_N}$ of diameter $\delta_0$. By the choice of $\delta_0$, the leaves $L'$ and $E'$ intersect in a unique leaf of $\widetilde{\mathcal{G}}$ through $B_{x_N}$, which is also invariant by $\gamma$. It turns out the existence of a closed leaf of $\mathcal{G}$ in $M$, which completes the proof. 
	
	Similarly, if $l_L^+$ is a contracting ideal point in $\partial_{\infty}L$, then we could find a leaf interval $I_L$ containing $L$ and a non-trivial deck transformation $\gamma$ such that $\gamma$ strictly contracts $I_L$. Thus, there is a $\gamma$-invariant leaf $L'\in \widetilde{\mathcal{F}}_1$ intersecting a foliation box around a point $x_N\in L$. The same argument turns out that there is a leaf in the intersection $L'\cap E'$ invariant by $\gamma$, which projects to a closed leaf of $\mathcal{G}$ in $M$. Hence, we finish the proof.
\end{proof}

\medskip

Now, we present the proof of Theorem~\ref{closed-leaf-nonHausdorff}.

\begin{proof}[Proof of Theorem~\ref{closed-leaf-nonHausdorff}]
    By Proposition~\ref{intersection}, for every leaf $F\in \widetilde{\mathcal{F}}_1$, the foliation $\widetilde{\mathcal{G}}_F$ has non-Hausdorff leaf space. Let $l\subset L\cap E$ be any non-separated ray of $\widetilde{\mathcal{G}}_L$ for leaves $L\in\widetilde{\mathcal{F}}_1$, $E\in \widetilde{\mathcal{F}}_2$. Denote by $l_L^+\in\partial_{\infty}L$ and $l_E^+\in \partial_{\infty}E$ be the ideal points of $l$ in corresponding ideal circles. Lemma~\ref{expanding_transverse} implies that $l_E^+$ is an expanding ideal point in $\partial_{\infty}E$. By changing another leaf $L\in \widetilde{\mathcal{F}}_1$ if necessary, we can assume that $l$ is also a non-separated ray in $\widetilde{\mathcal{G}}_E$ thanks to \cite[Theorem 8.1]{BFP25transverse}. Then, we apply Proposition~\ref{expanding-both} by exchanging $L\in \widetilde{\mathcal{F}}_1$ and $E\in\widetilde{\mathcal{F}}_2$ to conclude that there exists a leaf $s\in \widetilde{\mathcal{G}}$ such that $\gamma(s)=s$ for some non-trivial deck transformation $\gamma$. Therefore, we obtain a closed leaf of $\mathcal{G}$ in $M$ by projecting $s$ through the universal covering map. Hence, the proof is complete.
\end{proof}

\section{Accessibility}\label{section-accessible}

In this section, we consider $f:M \rightarrow M$ as a partially hyperbolic diffeomorphism of a closed 3-manifold with non-virtually solvable fundamental group. Assume that $f$ is homotopic to the identity and $NW(f)=M$. Throughout this section, we consistently assume the orientability of manifold $M$ and its subbundles $E^s$, $E^c$, and $E^u$. Furthermore, we require that the partially hyperbolic diffeomorphism $f$ preserves these orientations up to an appropriate finite lift and iterate when necessary. This assumption entails no loss of generality, as explained in \cite{FP_hyperbolic,FU1}.

By Theorem~\ref{maptori}, there cannot be any 2-dimensional embedded torus tangent to $E^s\oplus E^u$. Suppose that $f$ is not accessible. Then, as shown in Theorem~\ref{su-foliation}, either there is an $f$-invariant minimal foliation $\mathcal{F}^{su}$, or there exists an $f$-invariant minimal lamination, denoted by $\Lambda^{su}$, tangent to $E^s\oplus E^u$ whose complementary regions are $I$-bundles. Moreover, the lamination $\Lambda^{su}$ can extend to a foliation without compact leaves.

\subsection{Absence of invariant leaves}

We start considering the action of a good lift of $f$ on the universal cover $\widetilde{M}$. In particular, we are going to show that there is no leaf of $\widetilde{\Lambda}^{su}$ fixed by a good lift. Following \cite{BFFP1}, we define:
\begin{defn}
	We say that a lift $\widetilde{f}: \widetilde{M}\rightarrow \widetilde{M}$ of a homeomorphism $f: M\rightarrow M$ is a \emph{good lift} if it satisfies the following properties:
	\begin{enumerate}
	    \item $\widetilde{f}$ has uniformly bounded distance from the identity (i.e., there exists $K>0$ such that $d_{\widetilde{M}}(x, \widetilde{f}(x))<K$ for any $x\in \widetilde{M}$);
        \item $\widetilde{f}$ commutes with every deck transformation.
	\end{enumerate}
\end{defn}

The definition above does not require the diffeomorphism $f$ to be partially hyperbolic. There always exists a good lift if $f$ is homotopic to the identity. In general, a good lift might not be unique, since one can construct a new good lift by composing a good lift with any element in the centralizer of the fundamental group. Subsequently, we fix $\widetilde{f}$ a good lift of $f$.

Recall that the fundamental group of $M$ is not (virtually) solvable, so the lamination $\Lambda^{su}$ has no compact leaves, see Theorem \ref{maptori}. We will make use of the following dichotomy.

\begin{thm}\cite[Corollary 3.10]{FP_hyperbolic}\label{dichotomy}
	
	Let $f:M\rightarrow M$ be a homeomorphism homotopic to the identity that preserves a lamination $\Lambda$ with $C^1$ non-compact leaves on a closed 3-manifold, and $\widetilde{f}$ be a good lift. Assume that each completion of a complementary region of $\Lambda$ is an $I$-bundle. Then either
	\begin{itemize}
	    \item there is a minimal sublamination of $\Lambda$ whose lifted leaves are all invariant by $\widetilde{f}$; or
        \item $\Lambda$ extends to a uniform $\mathbb{R}$-covered foliation $\mathcal{F}$, and $\widetilde{f}$ acts as a translation on the leaf space of $\widetilde{\Lambda}$ as a subset of the leaf space of $\widetilde{\mathcal{F}}$.
	\end{itemize}
\end{thm}

We need the following theorem to show the main result of this subsection.
\begin{thm}\cite{Mendes77}\label{Mendes}
	Any Anosov map on a plane has at most one fixed point.
\end{thm}

Recall that we are considering a minimal lamination $\Lambda^{su}$. In this subsection, we will discard the first case in Theorem \ref{dichotomy}.

The main purpose of this subsection is to discard the first case in Theorem \ref{dichotomy} in the partially hyperbolic setting. We mention that there is another proof of the proposition below provided in \cite[Section 6]{FP_hyperbolic}. However, the proof we present here is much simpler as a consequence of \cite[Theorem 1.4]{FU1}. Here, we use the unified notation $\Lambda^{su}$ to denote either an $f$-invariant minimal foliation or an invariant minimal proper lamination tangent to $E^s\oplus E^u$ for the given partially hyperbolic diffeomorphism $f$.

\begin{prop}\label{translation}
	Let $f$ be a non-accessible partially hyperbolic diffeomorphism homotopic to the identity in a closed 3-manifold and $\widetilde{f}$ be a good lift. Assume $\pi_1(M)$ is not virtually solvable and $NW(f)=M$. Then, the $su$-lamination $\Lambda^{su}$ is uniform and $\mathbb{R}$-covered, and $\widetilde{f}$ acts as a translation on the leaf space of $\widetilde{\Lambda}^{su}$. In particular, there is no points in $\widetilde{M}$ fixed by $\widetilde{f}$.
\end{prop}

\begin{proof}
	By Theorem~\ref{su-foliation}, $f$ preserves a minimal lamination $\Lambda^{su}$ without compact leaves. Here, $\Lambda^{su}$ is either a minimal foliation or a minimal proper lamination with $I$-bundles in each complementary region. Then, Theorem \ref{dichotomy} provides us a dichotomy: either the good lift $\widetilde{f}$ fixes all leaves of $\widetilde{\Lambda}^{su}$, or the lamination $\Lambda^{su}$ satisfies the desired property. Now we suppose that all leaves of $\widetilde{\Lambda}^{su}$ are invariant by $\widetilde{f}$. 
	
	In case $\Lambda^{su}$ is a minimal proper lamination, there are periodic points on its boundary leaves by Theorem~\ref{su-foliation}. If $\Lambda^{su}$ is a minimal foliation, then $f$ also admits a periodic point as an immediate corollary of \cite[Theorem 1.4]{FU1}. Up to a finite iterate, let $p\in \Lambda^{su}$ be a fixed point of $f$. We apply an argument in \cite{2020Seifert} in the following proof. 
	
	Lifting to the universal cover $\widetilde{M}$, we denote by $\tilde{p}$ a lifted point of $p$ and $L\in \widetilde{\Lambda}^{su}$ be the leaf through $\tilde{p}$. There is a deck transformation $\rho$ satisfying that $\rho\circ \widetilde{f}(\tilde{p})=\tilde{p}$. Note that $\widetilde{f}$ commutes with any deck transformation since it is a good lift. Then we have $\rho\circ \widetilde{f}\circ \rho(\tilde{p})=\rho \circ \rho\circ \widetilde{f}(\tilde{p})=\rho(\tilde{p})$, which implies that $\rho(\tilde{p})$ is also a fixed point of $\rho\circ \widetilde{f}$. The lift $\widetilde{f}$ fixes the leaf $L$, so does the deck transformation $\rho$. It follows that $\rho\circ \widetilde{f}$ is an Anosov map acting on the properly embedded plane $L$. Then this map has at most one fixed point by Theorem \ref{Mendes}, which means that $\rho(\tilde{p})$ and $\tilde{p}$ coincide and $\rho$ is the identity map. 
	
	Now the point $\tilde{p}$ is a fixed point of $\tilde{f}$. Also by the property of a good lift, each point $\tau(\tilde{p})$ is a fixed point of $\tilde{f}$ in the leaf $\tau(L)$ for any deck transformation $\tau$. As the lamination $\Lambda^{su}$ is minimal, there is a non-trivial deck transformation $\tau$ so that the point $\tau(\tilde{p})$ is arbitrarily close to the leaf $L$. The deck transformation $\tau$ is trivial if the leaf $\tau(L)$ coincides with $L$. We assume that $\tau(L)$ is distinct from $L$. By transversality, there is one of the center curves through $\tau(\tilde{p})$, denoted by $c$, that is $\widetilde{f}$-invariant and intersects the leaf $L$ in a point distinct from $\tilde{p}$. There could be more than one center curves through $\tau(\tilde{p})$ since we did not assume the dynamical coherence. As the lamination $\Lambda^{su}$ has no Reeb component, the center curve $c$ can intersect $L$ in at most one point by Theorem \ref{taut}. Denote by $\tilde{q}=L\cap \gamma$ the unique intersection point. Both the center curve $\gamma$ and the leaf $L$ are invariant under the action $\widetilde{f}$. It turns out that $\tilde{q}$ and $\tilde{p}$ are two distinct point in $L$ fixed by $\widetilde{f}$, which is a contradiction.
\end{proof}

\subsection{A technical lemma}

Here, we provide an independent technical lemma that may have it own interest and may be applicable to other contexts.

We say that a one-dimensional foliation $\mathcal{T}$ is \emph{regulating} for a transverse codimension-one foliation $\mathcal{F}$ if each lifted leaf in $\widetilde{\mathcal{T}}$ intersects every lifted leaf of $\widetilde{\mathcal{F}}$.

\begin{lem}\label{technical-lemma}
    Let $\mathcal{F}$ be a codimension-one uniform Reebless foliation in a compact manifold $M$, and $\mathcal{T}$ be a one-dimensional foliation transverse to $\mathcal{F}$. Assume that both foliations are invariant under a homeomorphism $f$. Then, either
    \begin{itemize}
        \item $\mathcal{T}$ is not regulating for $\mathcal{F}$; or
        \item for any $L, F\in \widetilde{\mathcal{F}}$, the curves of $\widetilde{\mathcal{T}}$ between $L$ and $F$ have uniformly bounded lengths. Moreover, if $f$ uniformly expands (or contracts) each curve of $\mathcal{T}$, then $\widetilde{f}$ has an invariant leaf of $\widetilde{\mathcal{F}}$.
    \end{itemize}
\end{lem}
\begin{proof}
    For any $L, F\in \widetilde{\mathcal{F}}$, we denote by $d:=d_H(L, F)$ their Hausdorff leaf space, which is bounded by the uniformness of $\mathcal{F}$. As $M$ is compact, we choose a fundamental domain $K$ in $\widetilde{M}$, which is a connected compact set with $\pi(K)=M$ for the covering map $\pi:\widetilde{M}\rightarrow M$.

    Assume $\mathcal{T}$ is regulating for $\mathcal{F}$, i.e., each leaf of $\widetilde{\mathcal{T}}$ intersects every leaf of $\widetilde{\mathcal{F}}$. Then, each leaf of $\widetilde{\mathcal{T}}$ intersects every leaf of $\widetilde{\mathcal{F}}$ in exactly one point by Theorem~\ref{taut}. It implies that the leaf space of $\widetilde{\mathcal{F}}$ is homeomorphic to each leaf of $\widetilde{\mathcal{T}}$, and thus it is homeomorphic to $\mathbb{R}$. For any leaf $\alpha\in \widetilde{\mathcal{T}}$, we define a total order $<$ according to the transverse orientation of $\widetilde{\mathcal{F}}$, which determines an order in the leaf space of $\widetilde{\mathcal{F}}$. Define a projection $\tau: \widetilde{M}\rightarrow \alpha$ such that $\tau(x):=\widetilde{\mathcal{F}}(x)\cap \alpha$. Then, the compact set $K$ is projected by $\tau$ to a compact interval on $\alpha$. Notice that $\tau$ is monotone, that is, it preserves the order $<$. 

    Denote by $B(K, d)$ the closed set of points with distance at most $d$ from $K$. The projection $\tau(B(K, d))$ is a closed interval $I$ of $\alpha$. Let $I^+, I^-$ be respectively the uppermost and lowermost points of $I$ according to the order $<$. For any point $z\in \alpha$ with $I^+<z$, the leaf $\widetilde{\mathcal{F}}(z)$ is disjoint with $B(K, d)$ and thus its Hausdorff distance from any leaf through $K$ is larger than $d$. Analogously, the leaf through any point $z<I^-$ is $d$-apart from $K$. It implies that for any leaf $P\in \widetilde{\mathcal{F}}$, if $d_H(P, Q)\leq d$ for some leaf $Q\in \widetilde{\mathcal{F}}$ with $Q\cap K\neq \emptyset$, then we have $\tau(P)\in I$. We mention that, a priori, a leaf through a point in $I$ does not necessarily have $d$-Hausdorff distance from a leaf intersecting $K$.

    For any $p\in L$, there exists an element $\gamma\in \pi_1(M)$ with $\gamma(p)\in K$. Then, the leaf $\gamma(F)$ must intersect $B(K, d)$ and we have $\tau(\gamma(F))\in I$. Now, the leaf $\gamma(F)$ is bounded by the leaves of $\widetilde{\mathcal{F}}$ through $I^+$ and $I^-$. Without loss of generality, we assume the relation $\gamma(L)<\gamma(F)<\widetilde{\mathcal{F}}(I^+)$. Denote by $a:=\widetilde{\mathcal{T}}(\gamma(p))\cap\gamma(F)$ and $b:=\widetilde{\mathcal{T}}(\gamma(p))\cap \widetilde{\mathcal{F}}(I^+)$. One can deduce that the length of $\widetilde{\mathcal{T}}$-curve between $\gamma(p)$ and $a$ is smaller than that between $\gamma(p)$ and $b$. We can define a continuous function $\eta: K\rightarrow \mathbb{R}^+$ such that $\eta(x)$ is the length of $\widetilde{\mathcal{T}}$-curve between $x$ and the point $\widetilde{\mathcal{T}}(x)\cap\widetilde{\mathcal{F}}(I^+)$. By the compactness of $K$, $\eta$ has an upper bound $C$ only depending on $d$. It turns out that the length of $\widetilde{\mathcal{T}}$-curve from $\gamma(p)$ to $\gamma(F)$ is bounded by $C$ and so is the length of the curve from $p$ to $F$. Since $p$ is arbitrarily chosen, we conclude that any curve of $\widetilde{\mathcal{T}}$ between $L$ and $F$ has uniformly bounded length $C$.

    If $f$ uniformly expands every curve of $\mathcal{T}$, then each interval in the leaf space of $\widetilde{\mathcal{F}}$ is expanded under $\widetilde{f}$ since any leaf of $\widetilde{\mathcal{T}}$ is homeomorphic to the leaf space of $\widetilde{\mathcal{F}}$. We can find a bounded interval $J$ in this leaf space so that $J\subset\widetilde{f}(J)$. This implies the existence of an $\widetilde{f}$-invariant leaf of $\widetilde{\mathcal{F}}$. The same holds when $f$ uniformly contracts curves of $\mathcal{T}$.
\end{proof}

\subsection{Leaf-wise Gromov hyperbolicity} 

We employ the following theorem showing leafwise Gromov hyperbolicity for our lamination $\Lambda^{su}$. The result holds more generally for any minimal lamination with no compact leaves satisfying that the completions of its complementary regions are $I$-bundles. 

\begin{thm}\cite[Corollary 5.8]{FP_hyperbolic}\label{Gromov-su}
	Let $f$ be a partially hyperbolic diffeomorphism on a 3-manifold $M$ so that it is not accessible and $NW(f)=M$. If the fundamental group $\pi_1(M)$ is not (virtually) solvable, then the leaves of $\Lambda^{su}$ are uniformly Gromov hyperbolic.
\end{thm}

\begin{rmk}
	The proofs of this theorem relies on Candel's Uniformization Theorem \cite{Candel93} (see also \cite[Section 12.6]{CC00I}). The original Candel's Uniformization Theorem states for leaf-wise smooth foliations. In the light of \cite{Calegari01AGT}, it could be applied for surface laminations or foliations with no regularity requirement for leaves.
\end{rmk}

We utilize the following theorem, compiled from results in \cite{BI08, 2011TORI, BFFP2}. For background on branching foliation leaf spaces, see \cite[Section 3]{BFFP2}. A branching foliation is $\mathbb{R}$-covered precisely when its $\epsilon$-approximating foliations $\mathcal{W}_\epsilon$ are $\mathbb{R}$-covered for sufficiently small $\epsilon > 0$, as established in \cite[Proposition 3.16]{BFFP2}.  

\begin{thm}\label{cs_branching}
	Let $f: M\rightarrow M$ be a partially hyperbolic diffeomorphism on a closed 3-manifold whose fundamental group is not (virtually) solvable. Then up to a finite lift and iterate, there is an $f$-invariant branching foliation $\mathcal{W}^{cs}$ tangent to $E^c\oplus E^s$ without compact leaves, which is well-approximated by foliations $\mathcal{W}^{cs}_{\epsilon}$. Moreover, if $f$ is homotopic to the identity and $NW(f)=M$, then either
	\begin{itemize}
	    \item some iterate of $f$ is a discretized Anosov flow; or
        \item $\mathcal{W}^{cs}$ is $\mathbb{R}$-covered and uniform, and $\widetilde{f}$ acts as a translation on the leaf space of $\widetilde{\mathcal{W}}^{cs}$.
	\end{itemize}
\end{thm}

Using minimality of $\mathcal{W}^{cs}$ (Lemma~\ref{minimal-branching}), we can also obtain the leaf-wise Gromov hyperbolicity for $\mathcal{W}^{cs}$. 

\begin{prop}\cite[Theorem 5.1]{FP_hyperbolic}\cite[Appendix A.3]{BFP23collapsed}\label{Gromov-cs}
    Let $\mathcal{W}^{cs}$ be an $f$-invariant branching foliation for a partially hyperbolic diffeomorphism with $NW(f)=M$. Assume $\pi_1(M)$ is not virtually solvable. Then, the leaves of $\mathcal{W}^{cs}$ are Gromov hyperbolic. Analogously, the leaves of an well-approximated foliation $\mathcal{W}^{cs}_{\epsilon}$ are Gromov hyperbolic for $\epsilon$ small enough.
\end{prop}

\subsection{Proof of Theorem \ref{accessible-id}}

Let $f:M\rightarrow M$ be a partially hyperbolic diffeomorphism homotopic to the identity and its non-wandering set is the whole manifold. There is no loss of generality by considering a finite cover and a finite iterate. We assume that the manifold $M$ and the bundles $E^s$, $E^c$ and $E^u$ are all orientable and their orientations are preserved by $f$.

Suppose that $f$ is not accessible and the fundamental group of the ambient manifold is not (virtually) solvable. As shown in Theorem \ref{maptori}, there cannot be any torus tangent to $E^s\oplus E^u$. By Theorem \ref{su-foliation}, there is a unique $f$-invariant minimal lamination $\Lambda^{su}$ tangent to $E^s\oplus E^u$. This lamination is either an invariant $su$-foliation or a proper lamination that can trivially extend to a foliation. In either case, the diffeomorphism $f$ admits a periodic point in $\Lambda^{su}$. Indeed, it can be deduced directly by \cite[Theorem 1.4]{FU1} if $\Lambda^{su}$ is a foliation. In the proper lamination case, $\Lambda^{su}$ contains dense periodic points in its boundary leaves by Theorem~\ref{su-foliation} (alternatively, we can also obtain periodic points in $\Lambda^{su}$ using $I$-bundle structure and \cite[Theorem 1.4]{FU1}).

Let $\tilde{f}$ be a good lift of $f$ on the universal cover $\widetilde{M}$. Proposition \ref{translation} shows that $\Lambda^{su}$ is a uniform $\mathbb{R}$-covered lamination, and the good lift $\tilde{f}$ acts as a translation on the leaf space of $\widetilde{\Lambda}^{su}$. Moreover, the leaves of $\Lambda^{su}$ are uniformly Gromov hyperbolic. Since the complementary regions of $\Lambda^{su}$ are $I$-bundles, we can collapse those complementary regions to obtain a minimal foliation, denoted by $\mathcal{F}^{su}$. After collapsing, the manifold would not lose any differentiability since it is still a topological 3-manifold, by a classical result of Moise \cite{Moise77}. Moreover, the fundamental group keeps to be non-virtually solvable. We still use the notations as before. Thus, $\mathcal{F}^{su}$ is an $f$-invariant uniform $\mathbb{R}$-covered minimal foliation by non-compact Gromov hyperbolic leaves. Together Theorem~\ref{cs_branching} and Theorem~\ref{DAF}, the branching foliation $\mathcal{W}^{cs}$ and its well-approximated foliation $\mathcal{W}^{cs}_{\epsilon}$ are uniform and $\mathbb{R}$-covered for $\epsilon$ small enough. Moreover, the good lift $\widetilde{f}$ acts as a translation on the leaf space of $\widetilde{\mathcal{W}}^{cs}_{\epsilon}$. Using Theorem~\ref{cs_branching}, Proposition~\ref{Gromov-cs} and Lemma~\ref{minimal-branching}, $\mathcal{W}^{cs}_{\epsilon}$ is a minimal foliation by non-compact Gromov hyperbolic leaves.

We utilize the results in preceding sections by applying $\mathcal{F}^{su}=\mathcal{F}_1$ and $\mathcal{W}^{cs}_{\epsilon}=\mathcal{F}_2$. Let $\mathcal{G}= \mathcal{F}^{su}\cap\mathcal{W}^{cs}_{\epsilon}$ be the intersected one-dimensional subfoliation. Denote by $\widetilde{\mathcal{F}}^{su}$, $\widetilde{\mathcal{W}}^{cs}_{\epsilon}$ and $\widetilde{\mathcal{G}}$ be the corresponding lifts to the universal cover $\widetilde{M}$. As in Lemma~\ref{singlelimit}, each ray of $\widetilde{\mathcal{G}}_F$ accumulates in a single ideal point of $\partial_{\infty}F$ for any $F\in \widetilde{\mathcal{F}}^{su}$. The ideal limit set of $\widetilde{\mathcal{G}}_F$ is defined as the union of all ideal points in $\partial_{\infty}F$. With respect to the topology of ideal boundaries given in Section~\ref{subsection-topology}, we will discuss separately the cases where the ideal limit set of $\widetilde{\mathcal{G}}_F$ is dense or not.

Suppose the ideal limit set of $\widetilde{\mathcal{G}}_L$ is dense for one leaf $L\in \widetilde{\mathcal{F}}^{su}$, then the denseness holds for any leaf $F\in \widetilde{\mathcal{F}}^{su}$ (see Lemma~\ref{nondense=single}). Without loss of generality, we can assume that there exists a leaf $l_0\in \widetilde{\mathcal{G}}_{L_0}$ with distinct ideal points $l_0^+\neq l_0^-\in \partial_{\infty}L_0$ for some $L_0\in \widetilde{\mathcal{F}}^{su}$. Otherwise, if every leaf of $\widetilde{\mathcal{G}}_F$ admits a single ideal point in its ideal boundary, then by the denseness of limit set, we can find a leaf in the intersected foliation $\widetilde{\mathcal{F}}^{su}\cap \widetilde{\mathcal{W}}^{cu}_{\epsilon}$ with two distinct ideal points in $\partial_{\infty}F$. Then, we can proceed our argument by replacing $\widetilde{\mathcal{W}}^{cs}_{\epsilon}$ by $\widetilde{\mathcal{W}}^{cu}_{\epsilon}$.

In case the leaf space of $\widetilde{\mathcal{G}}_F$ is non-Hausdorff in $F\in \widetilde{\mathcal{F}}^{su}$, Theorem~\ref{closed-leaf-nonHausdorff} implies the existence of a leaf $L\in \widetilde{\mathcal{F}}^{su}$, a non-trivial deck transformation $\gamma$, and a leaf $s\in\widetilde{\mathcal{G}}_L$ satisfying $\gamma(s)=s$ and $\gamma(L)=L$. Note that the well-approximating map $h_{\epsilon}: \mathcal{W}^{cs}_{\epsilon}\rightarrow \mathcal{W}^{cs}$ is homotopic to the identity and lifts to a homeomorphism in $\widetilde{M}$ sending leaves of $\widetilde{\mathcal{W}}^{cs}_{\epsilon}$ to leaves of $\widetilde{\mathcal{W}}^{cs}$. As $\widetilde{h}_{\epsilon}(s)$ is at bounded Hausdorff distance from $s$, we can find a leaf of $\widetilde{\mathcal{F}}^s$ fixed by $\gamma$. This produces a closed stable leaf of $f$, which is absurd by partial hyperbolicity.

Now, the leaf space of $\widetilde{\mathcal{G}}_F$ is Hausdorff in $F\in \widetilde{\mathcal{F}}_1$. Then, as shown in Proposition~\ref{fan-weak}, the foliation $\mathcal{G}$ admits a closed leaf in $M$. The same argument as above applies to obtain a closed stable leaf of $f$, which is a contradiction. 

Thus, the ideal limit set of $\widetilde{\mathcal{G}}_F$ is not dense for every $F\in \widetilde{\mathcal{F}}^{su}$. Lemma~\ref{nondense=single} implies that every leaf $F\in \widetilde{\mathcal{F}}^{su}$ admits a degenerate limit set for $\widetilde{\mathcal{G}}$. As in Theorem~\ref{degenerate-limit-set}, the foliation $\mathcal{F}^{su}$ has trivial holonomy. Moreover, for the unstable foliation $\mathcal{F}^u$, the lifted foliations $\widetilde{\mathcal{W}}^{cs}_{\epsilon}$ and $\widetilde{\mathcal{F}}^u$ form a trivial foliated $\mathbb{R}$-bundle, so do $\widetilde{\mathcal{W}}^{cs}$ and $\widetilde{\mathcal{F}}^u$. It means that $\widetilde{\mathcal{F}}^u$ is regulating for $\widetilde{\mathcal{W}}^{cs}$. As $f$ uniformly expands each unstable curve, by Lemma~\ref{technical-lemma}, there exists at least one $\widetilde{f}$-invariant leaf of $\widetilde{\mathcal{W}}^{cs}$. This contradicts the fact that $\widetilde{f}$ acts as a translation on the leaf space of $\widetilde{\mathcal{W}}^{cs}$ given in Theorem~\ref{cs_branching}.

Hence, we complete the proof of Theorem~\ref{accessible-id}

\subsection{Other proofs}

Here, we present the proofs of rest results including Theorem~\ref{ergodic_id}, Theorem~\ref{ergodic=transitive}, and Corollary~\ref{cor-id-nonDAF}.

\begin{proof}[Proof of Theorem~\ref{ergodic_id}]
    In case where the fundamental group of manifold $M$ is not virtually solvable, we obtain K-property and ergodicity directly using Theorem~\ref{accessible-id} and \cite{08invent,BW10annals}. In the other case, we can apply \cite[Corollary 1.5]{HS21DA} to conclude our result using \cite{08invent,BW10annals} again unless the ambient manifold is the 3-torus. However, since the partially hyperbolic diffeomorphism we are considering is homotopic to the identity, it cannot exist on 3-torus by \cite{BI08, Parwani10}. Thus, we complete our proof.
\end{proof}

\begin{proof}[Proof of Theorem~\ref{ergodic=transitive}]
    The implication from ergodicity to transitivity is immediate due to the conservative property. To obtain ergodicity, we are sufficient to prove it in case where the fundamental group is virtually solvable thanks to results in \cite{HP15,Ham17CMH}. Thus, the conclusion follows directly from Theorem~\ref{accessible-id} and results in \cite{08invent,BW10annals}.
\end{proof}

\begin{proof}[Proof of Corollary~\ref{cor-id-nonDAF}]
    For a $C^1$ partially hyperbolic diffeomorphism $f$ homotopic to the identity, the induced action on the first homology group necessarily becomes trivial. This observation leads to the exclusion of (virtually) nilpotent fundamental groups for such manifolds, as established in \cite{BI08,Parwani10}. 

    Consider the case where $\pi_1(M)$ is (virtually) solvable. Through appropriate finite lifts and iterates, Theorem \ref{sol-nil} guarantees that the dynamically coherent diffeomorphism $\hat{f}$ can be realized as a discretized suspension Anosov flow. Descending through the covering map, this structural property persists for $f$ itself \cite{BW05}. The conclusion then follows directly from Theorem \ref{accessible-id}.
\end{proof}



\bibliographystyle{alpha}
\bibliography{ref}

\end{document}